\documentclass{amsart}[10pt]
\usepackage{amssymb,amsfonts, color,epsf}
\usepackage [cmtip,arrow]{xy}
\xyoption{all}
\usepackage {pb-diagram,pb-xy}
\usepackage{enumerate}
\usepackage{accents}
\usepackage[noautoscale]{youngtab}

\ifx\pdfpageheight\undefined
   \usepackage[dvips,colorlinks=true,linkcolor=blue,citecolor=red,%
      urlcolor=green]{hyperref}
   \usepackage[dvips]{graphicx}
   \makeatletter
   \edef\Gin@extensions{\Gin@extensions,.mps}
   \makeatother
\else
 \usepackage[pdftex]{graphicx}
   \usepackage[bookmarksopen=false,pdftex=true,breaklinks=true,%
      backref=page,pagebackref=true,plainpages=false,%
      hyperindex=true,pdfstartview=FitH,colorlinks=true,%
      pdfpagelabels=true,colorlinks=true,linkcolor=blue,%
      citecolor=red,urlcolor=green,hypertexnames=false%
      ]%
   {hyperref}
\fi
\usepackage{bm}

\usepackage{amsmath,amssymb,amsfonts}

\newtheorem{theorem}{Theorem}
\newtheorem{lemma}{Lemma}
\newtheorem{corollary}{Corollary}
\newtheorem{proposition}{Proposition}
\newtheorem{conjecture}{Conjecture}

\theoremstyle{definition}
\newtheorem{definition}{Definition}
\newtheorem{example}{Example}

\newtheorem{notation}{Notation}

\newtheorem{property}{Property}

\theoremstyle{remark}
\newtheorem{remark}{Remark}

\definecolor{DarkBlue}{rgb}{0,0.1,0.55}

\numberwithin{equation}{section}

\newcommand {\hide}[1]{}

 \newcommand {\sign} {\mbox{\bf sign}}

\newcommand {\junk}[1]{}

\newcommand {\R} {\mathrm{R}}

\newcommand {\C}     {\mathrm{C}}

\newcommand {\Real}[1]   {\mbox{$\mathbb{R}^{#1}$}}
\newcommand {\Sphere}{\mbox{${\bf S}$}}     % Sphere
 \newcommand {\re}         {\Real{}}

\newcommand {\Z}  {\mathbb{Z}}
 \newcommand {\N}         {\mathbb{N}}
\newcommand {\Q}         {\mathbb{Q}}
\newcommand{\kb}{\mathbf{k}}

\newcommand {\kk}         {\mathbf{k}}

\newcommand{\dd} {\mathbf{d}}

\newcommand{\F}{\mathbb{F}}

\newcommand {\ZZ} {\mathrm{Zer}}
\newcommand {\RR} {\mathrm{Reali}}

\newcommand {\la}   {{\langle}}
\newcommand {\ra}   {{\rangle}}
\newcommand {\eps} {{\varepsilon}}

\newcommand {\PP}     {\mathbb{P}} %projective space

\newcommand{\card}{\mathrm{card}}

\def\addots{\mathinner{\mkern1mu
\raise1pt\vbox{\kern7pt\hbox{.}}
\mkern2mu\raise4pt\hbox{.}\mkern2mu
\raise7pt\hbox{.}\mkern1mu}}

\newcommand{\HH}  {\mbox{\rm H}}

\DeclareMathOperator{\trace}{Tr}

\newcommand{\X}{\mathbf{X}}

\newcommand{\Y}{\mathbf{Y}}

\newcommand{\supp}{\mathrm{supp}}

\newcommand{\Ext}{\mathrm{Ext}}

\newcommand{\grad}{\mathrm{grad}}

\newcommand{\MV}{\mathrm{MV}}

\newcommand{\jj}{\mathbf{j}}

\newcommand{\Gr}{\mathrm{Gr}}
\newcommand{\AGr}{\mathrm{AffGr}}

\newcommand{\vol}{\mathrm{vol}}

\newcommand{\transversal}{\mathrm{Transversal}}
\newcommand{\Cont}{\mathrm{Cont}}
\newcommand{\gen}{\mathrm{gen}}
\newcommand{\Fix}{\mathrm{Fix}}

\begin{document}
\title[Multi-degree bounds on Betti numbers]
{
Multi-degree bounds on the Betti numbers of  real varieties and semi-algebraic sets and applications
}

\author{Saugata Basu}
\address{Department of Mathematics,
Purdue University, West Lafayette, IN 47906, U.S.A.}
\email{sbasu@math.purdue.edu}

\author{Anthony Rizzie}
\address{
Department of Mathematics, University of Connecticut, Storrs, CT 06269, U.S.A.}
\email{anthony.rizzie@uconn.edu}

%    General info
\date{\textbf{\today}}
%%\date{January 1, 1994 and, in revised form, June 22, 1994.}

%%\dedicatory{This paper is dedicated to our authors.}
\thanks{
Basu  was partially supported by NSF grants
DMS-1161629, CCF-1319080, CCF-1618918 and DMS-1620271.
Rizzie was partially supported by NSF grant
CCF-1319080.
}

\begin{abstract}
We prove new bounds on the Betti numbers of real varieties and semi-algebraic sets that have a more refined dependence on the degrees of the polynomials defining them than results known before. Our method also unifies several different types of results under a single framework, such as bounds depending on the total degrees, on multi-degrees, as well as in the case of quadratic and partially quadratic polynomials. 
The bounds we present in the case of partially quadratic polynomials offer a significant improvement over what was previously known.
Finally, we extend a result of Barone and Basu on bounding the number of connected components of real varieties defined by two polynomials of differing degrees to the sum of all Betti numbers, thus making progress on an open problem posed in that paper.
\end{abstract}

\maketitle
\tableofcontents

\section{Introduction}
\label{sec:intro}
Throughout this paper $\R$ will denote a fixed real closed field and $\C$ the algebraic closure of $\R$. 
For any
semi-algebraic subset $S \subset \R^k$ we denote by $b_i(S,\Z_2)$ the dimension of the $i$-th homology group,
$\HH_i(S,\Z_2)$, and by $b(S,\F) = \sum_{i \geq 0} b_i(S,\F)$ (we refer the reader to \cite[Chapter 6]{BPRbook2} for definition of homology groups of semi-algebraic sets defined over  
arbitrary real closed fields). 

\begin{remark}
\label{rem:universal}
Notice that by the universal coefficient theorem \cite[page 222]{Spanier}, $b_i(S,\Z_2) \geq b_i(S)$ (where $b_i(S)$ denotes the rank of the $i$-th homology group of
$S$ with integer coefficients), and thus an upper bound on $b_i(S,\Z_2)$ is automatically an upper bound
on $b_i(S)$.
\end{remark}

\subsection{Background}
\label{subsec:background}
The problem of bounding the Betti numbers of real algebraic varieties as well as semi-algebraic subsets
of $\R^k$, in terms of the format of their defining formulas, has been an active topic of investigation
for a long time starting from the first results bounding the Betti numbers of real varieties proved by Ole{\u\i}nik and Petrovski{\u\i} \cite{OP},
Thom \cite{T} and Milnor \cite{Milnor2}.  Later these results were extended to more general semi-algebraic sets \cite{B99,GaV, GV07}. These results were based on Morse-theoretic arguments involving bounding the number of critical points  of a  Morse function on a bounded, non-singular real algebraic hypersurface, using Bezout's theorem, arguments involving infinitesimal perturbations, and use of inequalities coming from the Mayer-Vietoris exact sequence. The bounds were singly exponential in the
dimension of the ambient space and polynomial in the number of polynomials used in the definition of the given semi-algebraic set, and also in the \emph{maximum of the total degrees of these polynomials}  (see Theorems \ref{thm:O-P-T-M}, \ref{thm:B99}, \ref{thm:GV07},  below for precise statements).

In another direction, bounds which are \emph{polynomial} in the dimension were proved for a restricted class of semi-algebraic sets -- namely, semi-algebraic sets defined by few (i.e,. a constant number of) quadratic equalities and inequalities. Barvinok \cite{Bar97}
proved a polynomial bound on the Betti numbers of semi-algebraic sets (see Theorem \ref{thm:barvinok} for a precise statement), which were sharpened in \cite{Bas05-first-Kettner, Lerario2012},
and also extended to a more general setting in \cite{BP'R07jems}  (see Theorem \ref{thm:BP'R07} for a precise statement). 
Some of these results were proved using different techniques than a simple counting of critical points -- for example,
a spectral sequence argument first proposed by Agrachev \cite{Agrachev, Ag, AgGa} plays an important role in the results proved in \cite{BP'R07jems,Lerario2012}.

Much more recently, because of certain new techniques developed in incidence geometry, more refined bounds than those mentioned above were needed. In particular, it was not enough to prove bounds which depended on the maximum of the degrees of the polynomials, 
and it was necessary to  prove bounds with
a more refined dependence on the sequence of degrees. Nearly optimal bounds on the zero-th Betti numbers (i.e. the number of connected components) of semi-algebraic sets was proved later in 
\cite{Barone-Basu11a, Barone-Basu13}  which has proved useful in applications. However, the techniques used to prove the
results in \cite{Barone-Basu11a, Barone-Basu13}  are not sufficient for bounding the higher Betti numbers.
Extending the bounds proved in \cite{Barone-Basu11a, Barone-Basu13} to the sum of all the Betti numbers (i.e. not just the zero-th Betti number) remains a challenging open problem in real algebraic geometry. 

The main contributions of the current paper are as follows.
The first contribution is to develop a single framework which allows one to prove 
the bounds on general semi-algebraic sets, as well as those defined by quadratic or even partially quadratic polynomials
(Theorems \ref{thm:algebraic-total-degree}, 
\ref{thm:partly-quadratic-total-algebraic}, 
\ref{thm:Lerario-new},
\ref{thm:partly-quadratic-total-semi}, \ref{thm:BP'R-new}).
Moreover, we  improve the known bounds in all of these cases. In the process, we also answer an open question of Lerario \cite{Lerario2012} on the asymptotic behavior of the Betti numbers of 
complete intersections of  projective quadrics  over $\C$ and $\R$
(Theorems \ref{thm:Lerario-new} and  \ref{thm:quadratic-projective}).

Additionally, the framework allows us to prove bounds in terms of the \emph{multi-degrees} of the polynomials
instead of the total degrees (Theorems \ref{thm:multi-algebraic}, \ref{thm:multi-semi}, \ref{thm:different-boxes-algebraic},
\ref{thm:different-boxes-semi}, 
\ref{thm:partly-quadratic-multi-algebraic}, \ref{thm:partly-quadratic-multi-semi}). 
We give several applications in which this new flexibility proves to be important (Theorems \ref{thm:pull-back}, 
\ref{thm:push-forward}, \ref{thm:Fourier-Mukai}, \ref{thm:transversals}).
Note that there have been some  other applications of multi-degree bounds in special cases  (see  for instance \cite{Riener-et-al} for a recent algorithmic application). 

As mentioned above,
extending the bounds proved in \cite{Barone-Basu11a, Barone-Basu13}  to the sum of all the Betti numbers  remains an open problem.
The second contribution of the current paper is extending the result \cite{Barone-Basu11a} to the sum of all the Betti numbers to the case of degree sequence of length bounded by $2$ (Theorems
\ref{thm:B-B-new} and \ref{thm:refined-betti-main}), thereby making progress on the open problem
posed in \cite{Barone-Basu13}.

\subsection{Prior Results}
\label{subsec:prior}

In this section we state more precisely the prior results mentioned in the previous section.
We first fix some notation that we will use for the rest of the paper.
\subsubsection{Basic notation and definition}
\label{subsec:basic-notation-definition}
\begin{notation}
\label{not:zeros}
  For $P \in \R [X_{1} , \ldots ,X_{k}]$ (resp. $P \in \C [ X_{1} ,
  \ldots ,X_{k} ]$) we denote by $\ZZ(P, \R^{k})$ (resp.
  $\ZZ (P, \C^{k})$) the set of zeros of $P$ in
  $\R^{k}$(resp. $\C^{k}$). More generally, for any finite set
  $\mathcal{P} \subset \R [ X_{1} , \ldots ,X_{k} ]$ (resp.
  $\mathcal{P} \subset \C [ X_{1} , \ldots ,X_{k} ]$), we denote by $\ZZ
  (\mathcal{P}, \R^{k})$ (resp. $\ZZ (\mathcal{P},
  \C^{k})$) the set of common zeros of $\mathcal{P}$ in
  $\R^{k}$  (resp. $\C^{k}$).  
  For a homogeneous polynomial $P \in \R [X_{0} , \ldots ,X_{k}]$,
  (resp. $P \in \C [ X_{0} ,\ldots ,X_{k} ]$)   
  we denote by $\ZZ(P, \PP_\R^{k})$ 
  (resp. $\ZZ (P, \PP_\C^{k})$) 
  the set of zeros of $P$ in $\PP_\R^{k}$.
  (resp. $\PP_\C^{k}$). 
  And, more generally, for any finite set of 
  homogeneous polynomials
  $\mathcal{P} \subset \R [ X_{0} , \ldots ,X_{k} ]$, 
  (resp.  $\mathcal{P} \subset \C [ X_{0} , \ldots ,X_{k} ]$),
   we denote by $\ZZ(\mathcal{P}, \PP_\R^{k})$ 
   (resp. $\ZZ (\mathcal{P},\PP_\C^{k})$) 
   the set of common zeros of $\mathcal{P}$ in $\PP_\R^{k}$.
  (resp. $\PP_\C^{k}$). 
  \end{notation}

\begin{notation}
  \label{not:sign-condition} For any finite family of polynomials $\mathcal{P}
  \subset \R [ X_{1} , \ldots ,X_{k} ]$, we call an element $\sigma \in \{
  0,1,-1 \}^{\mathcal{P}}$, a \emph{sign condition} on $\mathcal{P}$. For
  any semi-algebraic set $Z \subset \R^{k}$, and a sign condition $\sigma \in
  \{ 0,1,-1 \}^{\mathcal{P}}$, we denote by $\RR (\sigma ,Z)$ the
  semi-algebraic set defined by $$\{ \mathbf{x} \in Z \mid \sign (P (
  \mathbf{x})) = \sigma (P)  ,P \in \mathcal{P} \},$$ and call it the
  \emph{realization} of $\sigma$ on $Z$. More generally, we call any
  Boolean formula $\Phi$ with atoms, $P \{ =,>,< \} 0, P \in \mathcal{P}$, to
  be a \emph{$\mathcal{P}$-formula}. We call the realization of $\Phi$,
  namely the semi-algebraic set
  \begin{eqnarray*}
    \RR \left(\Phi , \R^{k} \right) & = & \left\{ \mathbf{x} \in \R^{k} \mid
    \Phi (\mathbf{x}) \right\}
  \end{eqnarray*}
  a \emph{$\mathcal{P}$-semi-algebraic set}. Finally, we call a Boolean
  formula without negations, and with atoms $P \{\geq, \leq \} 0$, $P\in \mathcal{P}$, to be a 
  \emph{$\mathcal{P}$-closed formula}, and we call
  the realization, $\RR \left(\Phi , \R^{k} \right)$, a \emph{$\mathcal{P}$-closed
  semi-algebraic set}.
\end{notation}

\subsubsection{General Bounds}
\label{subsubsec:general}
The first results on bounding the Betti numbers of real varieties were proved by
Ole{\u\i}nik and Petrovski{\u\i} \cite{OP}, Thom \cite{T} and Milnor \cite{Milnor2}.
Using a Morse-theoretic argument and Bezout's theorem they proved:

\begin{theorem}\cite{OP, T, Milnor2}
\label{thm:O-P-T-M}
Let $\mathcal{Q} \subset \R[X_1,\ldots,X_k]$ with $\deg(Q) \leq d, Q \in \mathcal{Q}$. Then, 
\begin{equation}
\label{eqn:O-P-T-M}
b(\ZZ(\mathcal{Q},\R^k),\Z_2) \leq d(2d-1)^{k-1}.
\end{equation}
\end{theorem}

\begin{remark}
\label{rem:O-P-T-M}
Theorem \ref{thm:O-P-T-M} is proved (see for example proof of Theorem 11.5.3 in \cite{BCR} for an exposition)  by first replacing the given variety by
a bounded basic, closed semi-algebraic set having the same homotopy type as $\ZZ(\mathcal{Q},\R^k)$
defined by a single polynomial inequality of total degree at most twice the
maximum of the degrees of the polynomials in $\mathcal{Q}$. Moreover, the critical points of the coordinate function $X_1$ are shown to be non-degenerate and their number can be bounded using Bezout's theorem.  
\end{remark}

\begin{remark}
\label{rem:O-P-T-M-2}
Also, note that the bound in Theorem \ref{thm:O-P-T-M} holds for dimensions of the homology groups with coefficients in any field and was proved in that generality. The same is true for some of the other results surveyed below. However, the new results in this paper give bounds only for Betti numbers over the field $\Z_2$ (because our technique for proving them  involves using Smith inequalities cf. Theorem \ref{thm:Smith}),  and are thus correspondingly weaker. On the other hand they do imply via the universal coefficients theorem (see Remark \ref{rem:universal}) the same bounds on the 
ranks of the homology groups with integer coefficients. Moreover, $\Z_2$-homology is very  natural in the context of real algebraic geometry. We will state all bounds for the $\Z_2$-Betti numbers from now on without comment (except in \S \ref{sec:refined-betti} below).
\end{remark}

Theorem \ref{thm:O-P-T-M} was later generalized to arbitrary semi-algebraic sets defined by quantifier-free
formulas in two steps. In the first step, Theorem \ref{thm:O-P-T-M} was extended to a particular
-- namely $\mathcal{P}$-closed semi-algebraic sets, where $\mathcal{P} \subset \R[X_1,\ldots,X_k]$ is a finite family of polynomials. The following theorem (which makes more precise an earlier result appearing in
\cite{B99}) appears in \cite{BPR02}.
\begin{theorem} \cite{BPR02}
\label{thm:B99}
If $S \subset \R^k$ is a $\mathcal{P}$-closed semi-algebraic set, then
\begin{eqnarray}
\label{eqn:B99}
b(S,\Z_2) \leq \sum_{i=0}^{k} \sum_{j=0}^{k-i} \binom{s+1}{j} 6^j d(2d-1)^{k-1},
\end{eqnarray}
where $s = \card(\mathcal{P}) > 0$ and $d = \max_{P \in \mathcal{P}} \deg(P)$.
\end{theorem}

Using an additional ingredient (namely, a technique to replace an arbitrary semi-algebraic set by a locally
closed one with a very controlled increase in the number of polynomials used to describe the given set),
Gabrielov and Vorobjov \cite{GaV} extended Theorem \ref{thm:B99} to arbitrary 
$\mathcal{P}$-semi-algebraic sets with only a small increase in the bound. Their result 
in conjunction with Theorem \ref{thm:B99} gives the following theorem.

\begin{theorem} \cite{GV07}
\label{thm:GV07}
If $S \subset \R^k$ is a $\mathcal{P}$-semi-algebraic set, then
\begin{eqnarray}
\label{eqn:GaV}
b(S,\Z_2) \leq \sum_{i=0}^{k} \sum_{j=0}^{k-i} \binom{2ks+1}{j} 6^j d(2d-1)^{k-1},
\end{eqnarray}
where $s = \card(\mathcal{P})$ and $d = \max_{P \in \mathcal{P}} \deg(P)$.
\end{theorem}

\subsubsection{Quadratic and partially quadratic case}
\label{subsubsec:quadratic}
Semi-algebraic sets defined by few quadratic inequalities are topologically simpler. This was first noticed
by Agrachev \cite{Agrachev, Ag, AgGa} who proved a bound which is polynomial in the number of
variables and exponential in the number of inequalities for \emph{generic} quadratic inequalities. 
The technique introduced by Agrachev was very important in later developments as well.
Independently, using a different technique (closer to the spirit of Morse theoretic arguments) Barvinok \cite{Bar97}  proved the following theorem (no genericity assumption is required).
 
\begin{theorem} \cite{Bar97}
\label{thm:barvinok}
Let $S \subset \R^k$ be defined by $P_1 \geq 0, \ldots, P_s \geq 0, \; \deg (P_i) \leq 2, \; 1 \leq i \leq s$. Then, $$ b(S,\Z_2) \leq k^{O(s)}.$$
\end{theorem}

This bound was later sharpened in \cite{Bas05-first-Kettner} and  further sharpened  in the case of 
algebraic sets by Lerario in \cite[Theorem 15]{Lerario2012},
where the following nearly optimal result was proved. 
\begin{theorem} \cite{Lerario2012}
\label{thm:lerario}
Let $\mathcal{Q} \subset \R[X_0,\ldots,X_k]$ be a set of $\ell$ quadratic forms, and $V =\ZZ(\mathcal{Q},\PP_\R^k)$ be 
the \emph{projective} variety defined by $\mathcal{Q}$.
Then, 
\[
b(V,\Z_2) \leq (O(k))^{\ell-1}.
\]
\end{theorem}

Theorem \ref{thm:barvinok} was later extended in 
\cite{B00} where the following theorem was proved. Notice that this bound is polynomial
even in the number of inequalities (for fixed $\ell$).

\begin{theorem} \cite{B00}
\label{thm:B00}
Let $\ell$ be any fixed number and let $\R$ be a real closed field. Let $S \subset \R^k$ be defined by $P_1 \geq 0,\ldots, P_s \geq 0, \; \deg(P_i) \leq 2$. Then, $$ b_{k-\ell }(S,\Z_2) \leq \binom{s}{\ell }k^{O(\ell )}.$$
\end{theorem}

Theorem \ref{thm:B00} was further improved using bounds on the Betti numbers of non-singular complete intersections and the Smith inequality (Theorem \ref{thm:Smith}) in
\cite{Bas05-first-Kettner} where the following theorem is proved.

\begin{theorem}\cite{Bas05-first-Kettner}
\label{thm:basukettner}
Let $\mathcal{P}=\{P_1,\ldots,P_s\} \subset \R[X_1,\ldots,X_k], \; s \leq k$. Let $S \subset \R^k$ be defined by $$P_1 \geq 0,\ldots,P_s \geq 0$$
with $\deg(P_i)\leq 2$. Then, for $0\leq i \leq k-1$,
$$b_i(S,\Z_2)\leq \frac{1}{2} \left( \sum_{j=0}^{min\{s,k-i\}} \binom{s}{j} \binom{k+1}{j} 2^j \right).$$
In particular, for $2\leq s \leq k/2$, we have
$$b_i(S,\Z_2) \leq \frac{1}{2} 3^s \binom{k+1}{s} \leq \frac{1}{2} \left( \frac{3e(k+1)}{s} \right) ^2.$$
\end{theorem}

Finally, in \cite{BP'R07jems} the authors also prove a 
result that generalizes the bounds on Betti numbers of general semi-algebraic sets (defined by $s$ polynomials having degrees bounded by $d$, cf. Theorem \ref{thm:B99}), as well as the bounds in the
quadratic case (cf. Theorems \ref{thm:barvinok}, \ref{thm:B00}, and \ref{thm:basukettner}).
More precisely they prove:

\begin{theorem}\cite{BP'R07jems}
\label{thm:BP'R07}
Let
$\mathcal {P}_1 \subset \R[X_1,\ldots,X_{k_1}]$,
a finite set of polynomials 
with
\[
\deg_{X}(P) \leq d, P \in \mathcal{P}_1, \card({\mathcal P}_1)=s,
\]
and let 
$\mathcal{P}_2\subset  \R[X_1,\ldots,X_{k_1},Y_1,\ldots,Y_{k_2}]$,
a finite set of polynomials
with 
\[
\deg_{X}(P) \leq d, 
\deg_{Y}(P) \leq 2,  P\in \mathcal{P}_2, \card(\mathcal{P}_2)=m,
\]

Let $S \subset \R^{k_1+k_2}$ 
be a $(\mathcal{P}_1 \cup \mathcal{P}_2)$-closed semi-algebraic set. Then
\begin{eqnarray}
\label{eqn:BP'R07-1}
b(S,\Z_2) &\leq& 
k_2^2 (O(k_2 +s+m)k_2 d)^{k_1+2m}. 
\end{eqnarray}

In particular, for $m \leq k_2$,
$
\displaystyle{
b(S,\Z_2) \leq k_2^2 (O(s+k_2)k_2 d)^{k_1+2m}. 
}
$
\end{theorem}

\begin{remark}
\label{rem:BP'R07-1}
In particular, if in Theorem \ref{thm:BP'R07}, $\mathcal{P}_1 = \emptyset$ (and hence $s=0$), and
$m, k_1 < k_2$, we get
\begin{eqnarray}
\label{eqn:BP'R07}
b(S,\Z_2) &\leq&  k_2^2 (O(m+k_2)k_2 d)^{k_1+2m}. 
\end{eqnarray}
\end{remark}

\begin{remark}
\label{rem:BP'R07-2}
The main tool used in the proof of Theorem \ref{thm:BP'R07} was a technique introduced by Agrachev in
\cite{AgGa,Ag,Agrachev} and later exploited by several authors \cite{BP'R07jems, Agrachev-Lerario, Lerario2014, Lerario2012} 
for bounding the Betti numbers of semi-algebraic sets defined by quadratic polynomials.  

The techniques used in the proof of the  theorems  corresponding to Theorem \ref{thm:BP'R07} in the current paper  
(namely, Theorems \ref{thm:partly-quadratic-total-algebraic},
\ref {thm:partly-quadratic-total-semi}, and \ref{thm:BP'R-new}) are quite different -- involving the method of infinitesimal perturbations, 
Mayer-Vietoris inequalities as explained in \cite[Chapter 7]{BPRbook2}, and bounds on the Betti numbers
of real affine varieties defined by partially quadratic polynomials proved in Proposition \ref{prop:many-total-mixed} below.  
\end{remark}

\subsubsection{Generic vs special}
\label{subsubsect:generic-vs-special}

\begin{definition}[Generic]
\label{def:generic}
In this paper while considering tuples  of polynomials
$(P_1,\ldots,P_\ell)$ whose supports are contained in some fixed tuple of convex polytopes
$\pmb{\Delta} = (\Delta_1,\ldots,\Delta_\ell)$ we will often make the assumption that the polynomials in the system are \emph{generic}. This means that the vector of coefficients of
$(P_1,\ldots,P_\ell)$ 
lie outside of  some  Zariski-closed subset, $V _{\pmb{\Delta}}$,  of strictly positive codimension in the vector space of coefficients of the polynomial system (the Zariski-closed
subset $V_{\pmb{\Delta}}$ will depend on $\pmb{\Delta}$). 

%%Note that this notion of being generic applies in both real and complex setting. 
Since the intersection of a Zariski-closed subset of $\C^N$ of strictly positive codimension with $\R^N \subset \C^N$ is a Zariski-closed subset of $\R^N$ of strictly positive (real)
codimension, while considering systems of generic complex polynomials, we can always assume that the coefficients belong to $\R$. 
However, note that in the complex case 
since the discriminant hypersurface does not disconnect the space of coefficients, the varieties defined by generic systems are all diffeomorphic to each other, 
the same is clearly not always true for the real parts of such varieties. 
\end{definition}

In many quantitative results in algebraic geometry, one assumes that the given system of polynomials is generic. Bounding the topological complexity of varieties defined by generic systems of polynomials 
(over $\R$ as well as $\C$) is often easy. However, such a result does not imply a bound in the non-generic situation. The following example which appears in \cite{Fulton}  is very well-known and shows that even for
the zero-th Betti number, the ``generic''  bound might not hold for all special systems.

\begin{example} \cite{Fulton}
  \label{ex:fulton} 
  Let $k = 3 $ and let
  \begin{eqnarray*}
    Q_1 & = & X_{3},\\
    Q_2 & = & X_{3},\\
    Q_3 & = & \sum_{i = 1}^2 \left( \prod^d_{j = 1} (X_i - j)^2 \right).
 \end{eqnarray*}

  The real variety defined by $\mathcal{Q} = \{ Q_1, Q_2, Q_3 \}$ is $0$-dimensional,
  and has $d^2$ isolated (in $\R^3$) points. However, a ``generic'' system of three
  polynomials in $\R[X_1,X_2,X_3]$ having degrees $1,1,2d$ will have by Bezout's theorem at most $2d$
  isolated points as its real zeros. Observe that even though the real variety $\ZZ(\mathcal{Q},\R^3)$ is zero-dimensional, the complex variety $\ZZ(\mathcal{Q},\C^3)$ is not, which accounts for this discrepancy. We refer  the reader to \cite{Barone-Basu13} for a  Bezout-type inequality that works over $\R$ as well.
\end{example}

There has been some work on bounding the number of connected components of real algebraic varieties
defined by systems of polynomials satisfying certain genericity conditions. For example, the following
theorem appears in \cite{Safey}.

\begin{theorem}\cite{Safey}
Let $(P_1,\ldots,P_s)\subset \Q[X_1,\ldots,X_k]$ (with $s \leq k-1$) generate a radical ideal and define a smooth algebraic variety $V \subset \C^k$ of dimension $k'$. Denote by $d_1,\ldots,d_s$ the respective degrees of $P_1,\ldots,P_s$ and by $d$ the maximum of $d_1,\ldots,d_s$. The number of connected components of $V \cap \R^k$ is bounded by
$$ d_1\cdots d_s\sum_{i=0}^{k'}(d-1)^{k-s-i}\binom{k-i}{k-i-s}.$$
Moreover, if $(P_1,\ldots,P_s)$ is a regular sequence, the number of connected components of $V \cap \R^k$ is bounded by
$$ d_1\cdots d_s\sum_{i=0}^{k'}(d-1)^{k-s-i}\binom{k-i-1}{k-i-s}.$$
\end{theorem}

Even though, bounds on generic systems do not immediately produce a bound on the Betti
numbers of general semi-algebraic sets, but with extra effort such bounds can be used to prove
(possibly worse) bounds for general semi-algebraic sets. This is in fact the approach taken in this
paper, but the approach already appears in the paper by Benedetti, Loeser and Risler \cite{Benedetti-Loeser},
which is the starting point of the results presented in the current paper. 
Using a clever reduction from the general case to the generic case they
 prove the following theorem.

\begin{theorem}\cite[Proposition 2.6]{Benedetti-Loeser}
\label{thm:BLR}
Let $\mathcal{P}=\{P_1,\ldots,P_\ell\}$ with $\deg(P_i) \leq d, \; 1 \leq i \leq \ell, \; P_i \in \R[X_1,\ldots,X_k]$.
 Then,
\[
b_0(\ZZ(\mathcal{P},\R^k),\Z_2) \leq \lambda (d,k,\ell ),
\]
where 
\[
\lambda(d,k,\ell)=Q_1(d,k)+2Q_2(d,k)+\cdots+2^{k-\ell-1}Q_{k-\ell}(d,k)+2^{k-\ell}\mu_d(\ell),
\] 
and each  $Q_i$ is a polynomial in $d$ of degree $k-i+1$, the leading coefficient of $Q_1$ is a polynomial 
in $k$ of degree $\ell-1$ with leading coefficient $(\ell+1)/2$, and the other terms polynomials in $k$ of degree $\max\{\ell-1,1\}$,
and $\mu_d(\ell)=d(2d-1)^{\ell-1}$.
\end{theorem}

\begin{remark}
\label{rem:BLR}
Of special interest here is that for every fixed $\ell$, and $k$ large enough (depending on $\ell$), and for $d$ tending to infinity,   
$\lambda(d,k,\ell)$ is asymptotically equal to 
%%$\frac{1}{2}(\ell+1)k^{\ell-1}d^k$ 
$\left(\frac{1}{2}(\ell+1)k^{\ell-1} + O_\ell(k^{\ell-2}) \right) d^k + O_{k,\ell}(d^{k-1})$
\cite[Corollary 2.7]{Benedetti-Loeser}
where the implied constant in the notation  $O_\ell$ (resp. $O_{k,\ell}$) depends only on $\ell$ (resp. $k,\ell$)
(compare with Theorem \ref{thm:algebraic-total-degree} and Remark \ref{rem:algebraic-total-degree}  below).
\end{remark}

In this paper, we consider the problem of bounding  the sum of all the Betti numbers
of real varieties and
semi-algebraic sets with a more refined dependence on the degrees of the polynomials. These refinements
are of two kinds. First,  we allow different blocks of variables to have different degrees
(see Theorems \ref{thm:multi-algebraic} and \ref{thm:multi-semi}).
Second, we allow
different polynomials to have different degrees (see Theorems \ref{thm:different-boxes-algebraic} and
\ref{thm:different-boxes-semi}).
Using our techniques we also improve existing bounds on the sum of the Betti numbers of real varieties and
semi-algebraic sets in terms of the number and total degrees of polynomials defining them
(Theorem \ref{thm:algebraic-total-degree}) as well as 
in the partially quadratic case (Theorems \ref{thm:partly-quadratic-total-algebraic}, 
\ref{thm:Lerario-new},
\ref{thm:partly-quadratic-total-semi}, 
\ref{thm:BP'R-new},
\ref{thm:partly-quadratic-multi-algebraic}, and
\ref{thm:partly-quadratic-multi-semi}).

The results mentioned above are all proved using a \emph{single framework}, which can be summarized
as follows. Using infinitesimal perturbation and inequalities derived from the Mayer-Vietoris exact
sequence, we reduce the problem of bounding the sum of the $\Z_2$-Betti numbers of a particular class of semi-algebraic sets to  bounding the same for a set of real algebraic varieties, which
are non-singular complete intersections in affine space. The perturbations need to be chosen carefully so that the degree dependencies of the various blocks of variables in the original set of polynomials are preserved.  
We then use Smith inequalities and 
a result of Khovanski{\u\i} \cite{Khovansky78} to bound the $\Z_2$-Betti numbers of these varieties.

We apply the results mentioned above to prove refined bounds on
the Betti numbers of pull-backs and direct images under polynomial maps 
(Theorems \ref{thm:pull-back}, \ref{thm:push-forward}, \ref{thm:Fourier-Mukai}), and as an application
of the last result (i.e. Theorem \ref{thm:Fourier-Mukai}) give a better bound (than possible before) on the Betti numbers of the space of affine subspaces of a fixed dimension that  meet a given semi-algebraic subset
of $\R^k$ (related to an important problem studied in discrete geometry). 
Finally, we make some progress on extending the theorem on refined bounds on
the number of connected components of real varieties defined by polynomials having different degrees,
to a bound on the sum of the Betti numbers using results proved in this paper and some other ingredients
(namely stratified Morse theory and Lefschetz duality from topology of manifolds). This is reported in
Theorem  \ref{thm:refined-betti-main}.

The rest of the paper is organized as follows. In \S \ref{sec:main} we state the main results
proved in this paper. In \S \ref{sec:preliminaries}, we state some preliminary results that are needed in the proofs of the main theorems.
In \S \ref{sec:main-results}, we prove the main theorems of the paper. In \S \ref{sec:applications}, we
prove bounds on the Betti numbers of pull-backs, direct images, and the space of transversals of semi-algebraic sets. In \S \ref{sec:refined-betti},
we prove a refined bound on the Betti numbers of varieties defined by polynomials having two different degree bounds. 
Finally, in \S \ref{sec:open} we list some open problems and possible future research directions.  

\section{Main Results}
\label{sec:main}
\subsection{Betti numbers of sets defined by polynomials of bounded total degree}
\label{subsec:total-degree}
We begin with the classical case of bounding the sum of the Betti numbers of varieties and semi-algebraic sets in terms of the total degrees of the polynomials defining them. This is the classical situation already considered by many authors and already surveyed in
\S \ref{subsec:prior},
but our methods produce slight improvements which we record here.
We prove the following theorems.
\begin{theorem}
\label{thm:algebraic-total-degree}
Let $\mathcal{Q}=\{Q_1,\ldots,Q_\ell\} \subset \R[X_1,\ldots,X_k]$ be a finite set of polynomials whose (total) degrees
are bounded by $d$ with $\ell>0$. Let $V$ denote $\ZZ(\mathcal{Q},\R^k)$. Then, 
$b(V,\Z_2)$ is bounded by 
\begin{equation}
\label{eqn:algebraic-total-degree}
\min\left(\sum_{j=1}^{k-1}\binom{\ell}{j}2^j (F_1(d',k,j) + F_2(d',k,j)) + \binom{\ell}{k}2^k d'^k+3 ,  \frac{1}{2}(1 + (2d-1)^k) \right),
\end{equation}
where 
\begin{eqnarray*}
F_1(d',k,j) &=& 
%%1 + (-1)^{k- \ell'+1} + 2d^{\ell'-1}(2+(i-1)d)^{k-i} \cdot \binom{k-1}{i-1} \\
1 + (-1)^{k- j+1} + \\
&&2d'^{j-1} \cdot\left( \sum_{h=0}^{k-j } \sum_{i=0}^{h}(-1)^{k-j+h}\binom{k}{h+j}\binom{j+i-2}{i}2^{h-i}d'^i\right) \\
& \leq &  2 \binom{k-2}{j -2} d^{k-1} + (O(d)) ^{k-2},\\
F_2(d,k,j) &=&  1 + (-1)^{k-j+1} + \binom{k-1}{j-1}(d'^k+k -1) \\
&\leq & \binom{k-1}{j-1}d^k + O(1)^k,
\end{eqnarray*}
and $d'$ is the least even integer $\geq d$.

In particular, 
if $k$ is fixed and $\ell \leq  k$, then for large $d$ we have
\begin{eqnarray}
\label{eqn:algebraic-total-degree-asymptotic}
b(V,\Z_2) &\leq & \sum_{j =1}^{\ell} 2^j  \binom{\ell}{j} \binom{k-1}{j-1}d^k 
 + (O(d))^{k-1}.
\end{eqnarray}
\end{theorem}

\begin{remark}
\label{rem:algebraic-total-degree}
Writing the bound in Theorem \ref{thm:algebraic-total-degree} as a polynomial in $d$, 
the leading coefficient is
\begin{eqnarray*}
\label{eqn:algebraic-total-degree-leading}
\sum_{j=1}^{\ell}2^j \binom{\ell}{j} \binom{k-1}{j-1} &\leq& 
 (3^\ell -1) \left(\sum_{j=1}^{\ell} \binom{k-1}{j-1}\right) \\
&\leq & \ell (3^\ell -1)\binom{k-1}{\ell-1} \mbox{ for } \ell <k/2.
\end{eqnarray*}
Thus, for every fixed $\ell$ and every $k$ sufficiently large (depending on $\ell$),  and as $d$ tends to infinity,  the bound in 
\eqref{eqn:algebraic-total-degree-asymptotic}  is asymptotically equal to 
\[
\left(\frac{\ell(3^\ell-1)}{(\ell-1)!} k^{\ell-1} + O_\ell (k^{\ell-2})\right) d^k + O_{k,\ell}(d^{k-1}),
\]
where the implied constant in the notation  $O_\ell$ (resp. $O_{k,\ell}$) depends only on $\ell$ (resp. $k,\ell$).
Notice that for $\ell >8$,
\[
\frac{\ell(3^\ell-1)}{(\ell-1)!} 
< \frac{1}{2}(\ell+1)
\]
(cf. Remark \ref{rem:BLR} following Theorem \ref{thm:BLR}).
Thus, for fixed $\ell$ and $k$ (sufficiently large) the bound in Theorem \ref{thm:algebraic-total-degree} is 
asymptotically better (as $d$ tends to infinity) than the bound in Theorem \ref{thm:BLR}.
Indeed the leading coefficient, $\frac{\ell(3^\ell-1)}{(\ell-1)!}$,
in the bound in Theorem \ref{thm:algebraic-total-degree} 
tends to $0$ exponentially fast with $\ell$, in contrast to the leading coefficient in
the bound in Theorem \ref{thm:BLR}.
\end{remark}

\begin{remark}
\label{rem:algebraic-total-degree-2}
Notice that the bound in Theorem  \ref{thm:algebraic-total-degree} is strictly better than the 
Ole{\u\i}nik-Petrovski{\u\i}-Thom-Milnor bound (Theorem \ref{thm:O-P-T-M})  
for all values of $\ell,d$ and $k$, with $d,k>1$, with equality in the case $d=1$ or $k=1$.
Assuming that $d,k>1$, we have that
$$\frac{1}{2}(1+(2d-1)^k)<d(2d-1)^{k-1}$$ since
$$1+(2d-1)^k<(2d-1)^{k-1}+(2d-1)^k=2d(2d-1)^{k-1}.$$
\end{remark}

\begin{remark}
\label{rem:katz}
Even though the results are incomparable, 
it is still interesting to note that using earlier results of Adolphson and Sperber \cite{Adolphson-Sperber} (who used methods involving  bounding exponential sums)
Katz \cite{Katz}  proved a bound of 
\begin{equation}
\label{eqn:katz}
6\cdot 2^r \cdot (2+(1+r d))^{k+1}
\end{equation}
on $\sum_{i \geq 0} \dim_{\Q_\ell} \HH^i_c(V,\Q_\ell)$ (here $\HH_c^*(V,\Q_\ell)$ denotes the $\ell$-adic cohomology groups with compact support), where $V \subset \C^k$ is an affine variety defined by $r$ polynomials in
$\C[X_1,\ldots,X_k]$ of total degrees bounded by $d$. While this result is incomparable with the results proved in this paper, and cannot be derived using our methods, 
notice that the bound in \eqref{eqn:katz} has an exponent of $k+1$ which is worse than the bound in Theorem
\ref{thm:algebraic-total-degree} (in the case $\ell,k$ are fixed and $d$ is large).
\end{remark}

\subsection{Betti numbers of sets defined by polynomials of bounded multi-degrees}
\label{subsec:multi}
We now consider the multi-degree case.

\begin{notation}
\label{not:multi}
Given, $\kk=(k_1,\ldots,k_p),\dd=(d_1,\ldots,d_p) \in \N^p, k=\sum_{i=1}^p k_i$, and $j>0$, we denote:
\[
G_{\gen}(\dd,\kk,j)=
1 + (-1)^{k- j+1} +  (k-j+2)^2\binom{k}{j-1}\binom{k}{\kk}^{-1} \frac{(1+p)^{3k-j+1}}{p(p+2)} d_1^{k_1}\cdots d_p^{k_p}.
\]
\end{notation}

\begin{theorem}
\label{thm:multi-algebraic}
Let $\mathcal{Q} = \{Q_1,\ldots,Q_\ell\} \subset \R[\X^{(1)},\ldots,\X^{(p)}]$ be a finite set of polynomials with $\ell>0$, where for $1 \leq i \leq p$, $\X^{(i)} = (X^{(i)}_1,\ldots,X^{(i)}_{k_i})$,
and $\deg_{\X^{(i)}}(Q) \leq d_i, \; d_i\geq 2$, for  all $Q \in \mathcal{Q}$. Let also $V = \ZZ(\mathcal{Q},\R^k)$, where 
$k = \sum_{i=1}^{p} k_i$.
Denote by $\dd = (d_1,\ldots,d_p)$ and $\kk=(k_1,\ldots,k_p)$.
Then,
\begin{eqnarray*}
\label{eqn:Gmin}
b(V,\Z_2) &\leq& G_{\min}(\dd,\kk,\ell)\\
&\leq& O(1)^kp^{3k}d_1^{k_1}\cdots d_p^{k_p},
\end{eqnarray*}
where $G_{\min}(\dd,\kk,\ell)$ equals
$$
\displaylines{
 \min\left( 3+\sum_{j=1}^{k} \binom{\ell}{j} 2^j (G_{\gen}(\dd',\kk,j)+G_{\gen}(\dd',\kk,j+1 )),
\; \frac{1}{2} G_{\gen}(2\dd,\kk,1)\right),
}
$$
$\dd' = (d_1',\ldots,d_p')$, and for $1 \leq i \leq p$, $d_i'$ is the least even integer $\geq d_i$.
\end{theorem}

\begin{theorem}
\label{thm:multi-semi}
Let $\mathcal{P} = \{P_1,\ldots,P_s\} \subset \R[\X^{(1)},\ldots,\X^{(p)}]$ be a finite set of polynomials with $s>0$, where for $1 \leq i \leq p$, $\X^{(i)} = (X^{(i)}_1,\ldots,X^{(i)}_{k_i})$,
and $\deg_{\X^{(i)}}(P) \leq d_i,\; d_i\geq 2$, for  all $P \in \mathcal{P}$.
Denote by $\dd = (d_1,\ldots,d_p)$ and $\kk=(k_1,\ldots,k_p)$.
Then, for each $i, 0 \leq i \leq k-1$,
\begin{eqnarray*}
\sum_{\sigma \in \{0,1,-1\}^{\mathcal{P}}}  b_i(\RR(\sigma,\R^k),\Z_2) 
&\leq&
\sum_{j=1}^{k-i} \binom{s}{j} 4^{j} G_{\min}(\dd,\kk,j)\\
&\leq& O(1)^k s^{k-i} p^{3k}d_1^{k_1}\cdots d_p^{k_p}.
\end{eqnarray*}
Furthermore, if $S$ is any $\mathcal{P}$-closed semi-algebraic set, then
\begin{eqnarray*}
b(S,\Z_2) &\leq& \sum_{i=0}^{k} \sum_{j=1}^{k-i} \binom{s+1}{j} 6^{j} G_{\min}(\dd,\kk,j)\\
&\leq& O(1)^k s^k p^{3k}d_1^{k_1}\cdots d_p^{k_p}.
\end{eqnarray*}
\end{theorem}

\subsection{Betti numbers of semi-algebraic sets defined by polynomials with different multi-degrees}
\label{subsec:different-boxes}
We now consider the case when different polynomials are allowed to have different multi-degrees.
\begin{notation}
For a matrix $\mathbf{d} \in \Z^{\ell \times k}$ and $I \subset [1,\ldots \ell], J \subset [1,k]$,
denote by $\mathbf{d}_{I,J}$ the sub-matrix extracted from $\mathbf{d}$ by taking the rows
indexed by $I$ and columns indexed by $J$.
We denote
\begin{eqnarray}
\label{eqn:def:K-gen}
K_{\gen} (\mathbf{d}) &=&
\left(  
\sum_{j=\ell}^{k} \sum_{J  \in  \binom{[1,k]}{j}}
(-1)^{k-j}
\sum_{\substack{
\boldsymbol{\alpha} =(\alpha_1,\ldots,\alpha_\ell) \in \Z_{>0}^\ell\\
\alpha_1+\cdots+\alpha_\ell = j}}  N (\mathbf{d}_{[1,\ell], J},\boldsymbol{\alpha})\right) ,
\end{eqnarray}
where the function $N(\cdot,\cdot)$ is defined in Eqn. \eqref{eqn:mixed-volume-of-boxes-refined}.
\end{notation}

\begin{theorem}
\label{thm:different-boxes-algebraic}
Let $\mathbf{d} = \Z_{\geq 2}^{\ell \times k}$.
Let  for $1 \leq i \leq \ell$,
$B_i = [0,d_{i,1}]\ \times \cdots \times [0,d_{i,k}] \subset \Z^k$. Let $\mathcal{Q} =\{Q_1,\ldots,Q_\ell\} \subset \R[X_1,\ldots,X_k]$,
with $\supp(Q_i) \subset B_i, 1 \leq i \leq \ell$, and let
$V = \ZZ(\mathcal{Q},\R^k)$.

Then,
\[
b(V,\Z_2) \leq K(\mathbf{d}),
\]
where 
\begin{eqnarray*}
K(\mathbf{d}) &=&
3+ \sum_{i=1}^{k} \sum_{I \subset [1,\ell], \card(I) = i} 2^{i+1}  K_{\gen}(\mathbf{d}''_{I,[1,k]}),
\end{eqnarray*}
\[
\mathbf{d}'' = \left[\begin{array}{c} \mathbf{d}'\\\mathbf{d}' \end{array}\right],
\] 
and $\dd' = [d'_{i,j}]_{\substack{1\leq i \leq \ell \\1 \leq j \leq k}}$, with $d'_{i,j}$ the least even number $\geq d_{i,j}$ for $1 \leq i \leq \ell, 1 \leq j \leq k$.
\end{theorem}

\begin{theorem}
\label{thm:different-boxes-semi}
Let $\mathbf{d} = \Z_{\geq 2}^{\ell \times k}$.
Let 
for $1 \leq i \leq \ell$,
$B_i = [0,d_{i,1}]\ \times \cdots \times [0,d_{i,k}] \subset \Z^k$.
Let $\mathcal{P} =\{P_1,\ldots,P_s\} \subset \R[X_1,\ldots,X_k]$,
with $\supp(P_i) \subset B_i, 1 \leq i \leq s$, $s >0$. 

Then, for each $i, 0 \leq i \leq k-1$,
\[
\sum_{\sigma \in \{0,1,-1\}^{\mathcal{P}}}  b_i(\RR(\sigma,\R^k),\Z_2) 
\leq
\sum_{j=1}^{k-i} \binom{s}{j} 4^{j} K(\dd).
\]
Furthermore, if $S$ is any $\mathcal{P}$-closed semi-algebraic set, then
\[
b(S,\Z_2) \leq \sum_{i=0}^{k} \sum_{j=1}^{k-i} \binom{s+1}{j} 6^{j} K(\dd).
\]
\end{theorem}

\begin{example}
\label{eg:different-boxes}
We give here an example in which Theorem \ref{thm:different-boxes-algebraic} can be applied.
Let $\mathcal{Q} = \{Q_1,\ldots,Q_\ell\} \subset \R[X_1,\ldots,X_k]$ with $\ell \leq k$.
Suppose that for each $i, 1\leq i \leq \ell$, $\deg_{X_i}(Q_i) \leq d_i, \deg_{X_j}(Q_i) = O(1), j\neq i$.
Moreover, assume that $d_1 \leq d_2 \leq \cdots \leq d_\ell$. 
Then using Theorem  
\ref{thm:different-boxes-algebraic}, one obtains immediately that
\begin{equation}
\label{eg:eqn:different-boxes-algebraic}
b(\ZZ(\mathcal{Q},\R^k),\Z_2) \leq  
O(1)^k \left(\sum_{\substack{
\boldsymbol{\alpha} \in \Z_{>0}^\ell \\ \alpha_1 +\cdots +\alpha_\ell = k}}
\Cont(\ell,k,\boldsymbol{\alpha}, \mathbf{1})
\right) d_1\cdots d_\ell^{k-\ell+1},
\end{equation}
where for any $m,n > 0$, and $\mathbf{r}= (r_1,\ldots,r_m)\in \Z_{>0}^m,
\mathbf{c} = (c_1,\ldots,c_n) \in \Z_{>0}^n, \sum_{i} r_i = \sum_{j} c_j$, 
$\Cont(m,n,\mathbf{r},\mathbf{c})$ denotes the number of matrices in
$\Z_{\geq 0}^{m \times n}$ with the vector of row-sums equal to $\mathbf{r}$, and the vector of
column-sums equal to $\mathbf{c}$ (such matrices are often referred to as \emph{contingency tables}). 
Note that the quantity 
\[
\left(\sum_{\substack{
\boldsymbol{\alpha} \in \Z_{>0}^\ell \\ \alpha_1 +\cdots +\alpha_\ell = k}}
\Cont(\ell,k,\boldsymbol{\alpha}, \mathbf{1})
\right)
\]
appearing in \eqref{eg:eqn:different-boxes-algebraic}
depends only on $k$ and $\ell$,  is independent of the $d_i$'s, and is bounded by $2^{O(k^2)}$ using results
in \cite{Barvinok2009} on the asymptotic number of contingency tables.

Notice that the dependence on the various degrees $d_i$ in the bound above is similar to the bound
proved in \cite{Barone-Basu13} on the number of semi-algebraically connected components of a real variety defined by 
polynomials of increasing \emph{total} degrees, with some added assumptions on the dimensions of
the intermediate varieties defined by some of the subsets of the polynomials. There are no dimensional restrictions 
for the bound in \eqref{eg:eqn:different-boxes-algebraic} to hold, and moreover the inequality in
\eqref{eg:eqn:different-boxes-algebraic} gives a bound on the sum of all Betti numbers not just on the 
zero-th one. However, the degree restrictions in the assumption for  \eqref{eg:eqn:different-boxes-algebraic} 
 is much stronger than just requiring that for each $i$,
the total degree of the polynomial $Q_i$ is bounded by $d_i + O(1)$ which would suffice for the result in   
\cite{Barone-Basu13} to hold.

Finally, note that using Alexander duality, the bound in  \eqref{eg:eqn:different-boxes-algebraic} is also an upper bound on
$b(\R^k \setminus \ZZ(\mathcal{Q},\R^k),\Z_2)$.
\end{example}

\subsection{Betti numbers of sets defined by quadratic and partially  quadratic polynomials}
\label{subsec:partly-quadratic}
In the following theorems  we improve the result in Theorem \ref{thm:BP'R07}. In Theorems
\ref{thm:partly-quadratic-total-algebraic} and 
\ref {thm:partly-quadratic-total-semi}, we assume that the set $\mathcal{P}_1$ is empty, and we are
able to provide more precise bounds in this situation. In Theorem \ref{thm:BP'R-new} the hypothesis is
the same as in Theorem \ref{thm:BP'R07}, and we improve the bound in 
Theorem \ref{thm:BP'R07} a significant way -- namely the dependence of the bound on $m$.

We first introduce the following notation.
\begin{notation}
\label{not:partly-quadratic-total}
In the following theorems, we will denote by $k=k_1+k_2$, and
\[
H_{\gen}(d,k_1,k_2,j)=2+(-1)^{k-j+1}+ j 2^{j} (k_1+k_2)^{j-1} \left( 2d(k_1+k_2)+1 \right)^{k_1}.
\]
\end{notation}

\begin{remark}
\label{rem:partly-quadratic-total}
Notice that for $j,k_1 < k_2$, 
\[
H_{\gen}(d,k_1,k_2,j) \leq (O(k_2))^{j-1} (O(d k_2))^{k_1}.
\]
\end{remark}

\begin{theorem}
\label{thm:partly-quadratic-total-algebraic}
Let $\mathcal{Q} = \{Q_1,\ldots,Q_\ell\} \subset \R[X_1,\ldots,X_{k_1},Y_1,\ldots,Y_{k_2}]$ be a finite set of polynomials with $\ell>0$, $\deg_\X(Q)\leq d, \; d\geq 2$, and $\deg_\Y(Q)\leq 2$ for all $Q \in \mathcal{Q}$. Let $V$ denote $\ZZ(\mathcal{Q},\R^k)$. Then,
\begin{eqnarray*}
b(V,\Z_2) &\leq&  H(d,k_1,k_2,\ell),
\end{eqnarray*}
where 
\begin{equation}
\label{eqn:partly-quadratic-total-algebraic}
H(d,k_1,k_2,\ell) =
3+\sum_{j=1}^{k} \binom{\ell}{j}2^j (H_{\gen}(d',k_1,k_2,j)+H_{\gen}(d',k_1,k_2,j+1)),
\end{equation}
where $d'$ is the least even integer $\geq d$.
In particular, for $\ell,k_1 \leq k_2$,
\begin{eqnarray}
\label{eqn:partly-quadratic-total-algebraic-asymptotic}
b(V,\Z_2)
&\leq & (O(k_2))^{\ell+k_1} d^{k_1}.
\end{eqnarray}
\end{theorem}

\begin{remark}
\label{rem:partly-quadratic-total-algebraic}
Notice that in the case $k_1=0$ (thus the polynomials in $\mathcal{Q}$ are \emph{fully} quadratic),
the bound in inequality \eqref{eqn:partly-quadratic-total-algebraic-asymptotic} 
reduces to $(O(k_2))^{\ell}$ almost recovering (i.e. up to a factor $k$) the bound in 
Theorem \ref{thm:lerario}
(albeit for affine varieties).
\end{remark}

\begin{remark}
\label{rem:Lerario-new}
It might also be possible with more work (using the same ideas as in the proof of Theorem
\ref{thm:quadratic-projective} taking into account signs)
to remove a factor of $k_2$ from the bound in Theorem
\ref{thm:partly-quadratic-total-algebraic}
(cf. Remark \ref{rem:quadric-projective}), and we leave this as an open question. 
\end{remark}

For projective varieties  in $\PP_{\R}^k$ defined by a fixed number of homogeneous quadratic polynomials we have the 
following bound that is asymptotically a slight improvement over the tightest bound known previously \cite[Theorem 15]{Lerario2012} (namely, the bound  $(O(k))^{\ell-1}$).

\begin{theorem}
\label{thm:Lerario-new}
For each fixed $\ell >0$, and for each set 
$\mathcal{P} \subset \R[X_0,\ldots,X_k]$ of homogeneous polynomials of degree $2$ of $\card(\mathcal{P}) \leq \ell$,
\begin{eqnarray}
\label{eqn:Lerario-new}
\nonumber
b(\ZZ(\mathcal{P},\PP_\R^k),\Z_2)  &\leq & 
k+1
\sum_{i=1}^{k} 
\binom{\ell}{i}2^i 
H_{\gen}'(k,i) \\
&\leq & \left(O\left(\frac{k}{\ell}\right)\right)^{\ell-1},
\end{eqnarray}
where
$$\displaylines{
H_\gen'(k,i) = 
(1+(-1)^{k-i+1})(k-i+1)+ \cr
 (-1)^{k-i} \left(\sum_{h=0}^{i-1} (-2)^h\left( \sum_{j=i}^{k} (-1)^{j+1} \binom{j}{h} \right) +  (k-i+1)\right).
}
$$
\end{theorem}

\begin{theorem}
\label{thm:partly-quadratic-total-semi}
Let $\mathcal{P} = \{P_1,\ldots,P_s\} \subset \R[X_1,\ldots,X_{k_1},Y_1,\ldots,Y_{k_2}]$ be a finite set of polynomials with $s>0$, $\deg_\X(P)\leq d, \; d\geq 2$, and $\deg_\Y(P)\leq 2$ for all $P \in \mathcal{P}$.

Then, for each $i, \; 0 \leq i \leq k-1$,
\begin{eqnarray}
\label{eqn:partly-quadratic-total-semi1}
\nonumber
\sum_{\sigma \in \{0,1,-1\}^{\mathcal{P}}}  b_i(\RR(\sigma,\R^k),\Z_2) 
&\leq&
\sum_{j=1}^{k-i} \binom{s}{j} 4^{j} H(d',k_1,k_2,j) \\
&\leq &
(O(k_2))^{s+k_1} d^{k_1} \mbox{ for } s,k_1 < k_2,
\end{eqnarray}
where $d'$ is the least even integer $\geq d$.
Furthermore, if $S$ is any $\mathcal{P}$-closed semi-algebraic set, then

\begin{eqnarray}
\label{eqn:partly-quadratic-total-semi2}
\nonumber
b(S,\Z_2) &\leq & \sum_{i=0}^{k}
  \sum_{j=1}^{k-i} \binom{s+1}{j} 6^{j} H(d',k_1,k_2,j) \\
 &\leq &       (k_1+k_2+1)  (O(k_2))^{s+k_1+1} d^{k_1}, \mbox{ for } k_1 < k_2, \\
 \nonumber
 &=& (O(k_2))^{s + k_1+2} d^{k_1}.
\end{eqnarray}
\end{theorem}

\begin{remark}
\label{rem:partly-quadratic-total-semi}
Notice that the  bound in inequality \eqref{eqn:partly-quadratic-total-semi2} in Theorem 
\ref{thm:partly-quadratic-total-semi} is significantly better than the  previous best bound
known on this quantity (namely,
inequality \eqref{eqn:BP'R07} in Remark \ref{rem:BP'R07-1}).
\end{remark}

\begin{theorem}
\label{thm:BP'R-new}
With the same notation as in Theorem \ref{thm:BP'R07},  
for each $i, \; 0 \leq i \leq k-1$ and assuming $m \leq k_2$,
$\sum_{\sigma \in \{0,1,-1\}^{\mathcal{P}_1 \cup \mathcal{P}_2}}  b_i(\RR(\sigma,\R^k),\Z_2)$
is bounded by  
\begin{eqnarray}
\label{eqn:BP'R-new1}
\nonumber
&&\sum_{j=1}^{k-i} 
\sum_{\substack{0 \leq j_1 \leq \min(s,k_1)\\  0\leq j_2\leq \min(m+1,k_1+k_2-j_1)\\ j_1+j_2 = j}}  
\binom{s}{j_1} \binom{m+1}{j_2}  5^{j} H
(2d,k_1,k_2,j_2+1) \\ 
&\leq&             (O(k_2))^{k_1+m+2} (O(sd))^{k_1}, 
\end{eqnarray}
for $m, k_1 < k_2$.

Furthermore, if $S$ is any $(\mathcal{P}_1 \cup \mathcal{P}_2)$-closed semi-algebraic set, then

\begin{eqnarray}
\label{eqn:BP'R-new2}
\nonumber
b(S,\Z_2) &\leq & \sum_{i=0}^{k}
\sum_{\substack{0 \leq j_1 \leq \min(s,k_1)\\ 0 \leq j_2 \leq \min(m+1,k_1+k_2-j_1)\\j_1+j_2 = j \leq k-i}}  
\binom{s}{j_1}\binom{m+1}{j_2}7^{j} H(2d,k_1,k_2,j_2+1) \\
 &\leq &             (O(k_2))^{k_1+m+3} (O(sd))^{k_1}, \mbox{ for }  m, k_1 < k_2.
\end{eqnarray}
\end{theorem}

\begin{remark}
\label{rem:BP'R-new}
Notice that the  bound in inequalities \eqref{eqn:BP'R-new1} and \eqref{eqn:BP'R-new2}
in  Theorem 
\ref{thm:BP'R-new} is significantly better than the corresponding bounds
in Theorem \ref{thm:BP'R07} (namely, in the dependence on $m$ and the exponent of $k_2$).
\end{remark}

\subsection{Betti numbers of semi-algebraic sets defined by partially quadratic polynomials with several blocks of variables}
\label{subsec:partly-quadratic-multi}
Lastly, we consider the case of partially quadratic polynomials, with the non-quadratically bounded variables
allowed to have different degrees. 
\begin{notation}
\label{not:partly-quadratic-multi}
In the following theorems, we will denote by $k=k_1+k_2, \; \dd=(d_1,\ldots,d_{k_1}) \in \N^{k_1}$ and
\[
M_{\gen}(\dd,k_1,k_2,j)=2 + (-1)^{k-j+1}+j 2^{j} k_1!(k_1+k_2)^{j-1}\left( 2(k_1+k_2)+1 \right)^{k_1}d_1\cdots d_{k_1}.
\]
\end{notation}

\begin{theorem}
\label{thm:partly-quadratic-multi-algebraic}
Let $\mathcal{Q} = \{Q_1,\ldots,Q_\ell\} \subset \R[X_1,\ldots,X_{k_1},Y_1,\ldots,Y_{k_2}]$ be a finite set of polynomials with $\ell>0$, $\deg_{X_i}(Q)\leq d_i, \; d_i\geq 2$, and $\deg_\Y(Q)\leq 2$ for all $Q \in \mathcal{Q}$. Let $V$ denote $\ZZ(\mathcal{Q},\R^k)$. Then,
\[
b(V,\Z_2) \leq 
M(\dd,k_1,k_2,\ell)
\]
where 
$$
M(\dd,k_1,k_2,\ell) = 
3+\sum_{j=1}^{k} \binom{\ell}{j}2^j (M_{\gen}(\dd',k_1,k_2,j)+M_{\gen}(\dd',k_1,k_2,j+1)),
$$
and 
where $\dd' = (d_1',\ldots,d_{k_1}')$ and for $1 \leq i \leq k_1$, $d_i'$ is the least even integer $\geq d_i$.

In particular, for $\ell,k_1 \leq k_2$,
\begin{eqnarray*}
b(V,\Z_2) &\leq& (O(k_2))^{\ell+k_1}d_1\cdots d_{k_1}.
\end{eqnarray*}
\end{theorem}

\begin{theorem}
\label{thm:partly-quadratic-multi-semi}
Let $\mathcal{P} = \{P_1,\ldots,P_s\} \subset \R[X_1,\ldots,X_{k_1},Y_1,\ldots,Y_{k_2}]$ be a finite set of polynomials with $s>0$, $\deg_{X_i}(P)\leq d_i,\; d_i\geq 2$, and $\deg_\Y(P)\leq 2$ for all $P \in \mathcal{P}$.
Then, for each $i, \; 0 \leq i \leq k-1$,
\begin{eqnarray*}
\sum_{\sigma \in \{0,1,-1\}^{\mathcal{P}}}  b_i(\RR(\sigma,\R^k),\Z_2) 
&\leq&
\sum_{j=1}^{k-i} \binom{s}{j} 4^{j} 
M(\dd',k_1,k_2,j)\\
&\leq&
(O(k_2))^{s+k_1} d_1\cdots d_{k_1} \mbox{ for } s,k_1 < k_2,
\end{eqnarray*}
where $\dd' = (d_1',\ldots,d_{k_1}')$ and for $1 \leq i \leq k_1$, $d_i'$ is the least even integer $\geq d_i$.
Furthermore, if $S$ is any $\mathcal{P}$-closed semi-algebraic set, then
\begin{eqnarray*}
b(S,\Z_2) &\leq& \sum_{i=0}^{k} \sum_{j=1}^{k-i} \binom{s+1}{j} 6^{j} M(\dd',k_1,k_2,j)\\
&\leq &(O(k_2))^{s+k_1+2} d_1\cdots d_{k_1}, \mbox{ for } s,k_1 < k_2,
\end{eqnarray*}
where $\dd' = (d_1',\ldots,d_{k_1}')$ and for $1 \leq i \leq k_1$, $d_i'$ is the least even integer $\geq d_i$.
\end{theorem}

\section{Preliminaries}
\label{sec:preliminaries}

We first recall some preliminary results that we will need in the paper.
\subsection{Real algebraic preliminaries}
\begin{notation}
  For $\R$ a real closed field we denote by $\R \left\langle \eps
  \right\rangle$ the real closed field of algebraic Puiseux series in $\eps$
  with coefficients in $\R$. We use the notation $\R \left\langle \eps_{1} ,
  \ldots , \eps_{m} \right\rangle$ to denote the real closed field $\R
  \left\langle \eps_{1} \right\rangle \left\langle \eps_{2} \right\rangle
  \cdots \left\langle \eps_{m} \right\rangle$. Note that in the unique
  ordering of the field $\R \left\langle \eps_{1} , \ldots , \eps_{m}
  \right\rangle$, $0< \eps_{m} \ll \eps_{m-1} \ll \cdots \ll \eps_{1} \ll 1$.
\end{notation}

\begin{notation}
  For elements $x \in \R \left\langle \eps \right\rangle$ which are bounded
  over $\R$ we denote by $\lim_{\eps}  x$ to be the image in $\R$ under the
  usual map that sets $\eps$ to $0$ in the Puiseux series $x$.
\end{notation}

\begin{notation}
\label{not:ext}
  If $\R'$ is a real closed extension of a real closed field $\R$, and $S
  \subset \R^{k}$ is a semi-algebraic set defined by a first-order formula
  with coefficients in $\R$, then we will denote by $\Ext \left( S, \R'
  \right) \subset \R'^{k}$ the semi-algebraic subset of $\R'^{k}$ defined by
  the same formula. It is well-known that $\Ext \left( S, \R' \right)$ does
  not depend on the choice of the formula defining $S$ {\cite{BPRbook2}}.
\end{notation}

\begin{notation}
\label{not:ball-and-sphere}
  For $x \in \R^{k}$ and $r \in \R$, $r>0$, we will denote by $B_{k} ( x,r )$
  the open Euclidean ball centered at $x$ of radius $r$. If $\R'$ is a real
  closed extension of the real closed field $\R$ and when the context is
  clear, we will continue to denote by $B_{k} ( x,r )$ the extension $\Ext
  \left( B_{k} ( x,r ) , \R' \right)$. This should not cause any confusion.
  We also denote by $\Sphere^{k-1}(x,r)$ the $(k-1)$-dimensional sphere in 
  $\R^k$, centered at $x$ and of radius $r$.
\end{notation}

\subsection{Topological preliminaries}

\subsubsection{Mayer-Vietoris inequalities}
Let $S_1,S_2$ be two closed semi-algebraic sets, and $\F$ any field of coefficients. We will use heavily the following
inequalities which are consequences of the exactness of the Mayer-Vietoris sequence. 
\begin{eqnarray} 
b_i(S_1 \cup S_2,\F) &\leq & b_i(S_1,\F)+b_i(S_2,\F)+b_{i-1}(S_1 \cap S_2,\F),  \label{eqn:MV:S1cupS2}\\
b_i(S_1\cap S_2,\F)  &\leq &b_i(S_1,\F)+b_i(S_2,\F)+b_{i+1}(S_1 \cup S_2,\F),  \label{eqn:MV:S1capS_2}\\
b_i(S_1,\F)+b_i(S_2,\F) &\leq & b_i(S_1 \cup S_2,\F)+b_i(S_1 \cap S_2,\F).\label{eqn:MV:S1+S2}
\end{eqnarray}

The following generalization in the case of more than two sets will also be useful for us (see for example
\cite[Proposition 7.33]{BPRbook2}).
\begin{proposition}
\label{prop:unionint}
Let $S_1,\ldots,S_s \subset \R^k, \; s\geq 1$, be closed semi-algebraic sets contained in a closed semi-algebraic set $T$ of dimension $k'$. For $S_{\leq t}=\bigcap_{1 \leq j \leq t}S_j$, and $S^{\leq t}=\bigcup_{1\leq j \leq t} S_j$. Also, for $J\subset \{1,\ldots,s\}, \; J \neq \emptyset$, let $S_J=\bigcap_{j\in J} S_j$, and $S^J=\bigcup_{j\in J} S_j$. Finally, let $S^\emptyset =T$. Then
\begin{enumerate}
\item For $0 \leq i \leq k'$, $$b_i(S^{\leq s},\F)\leq \sum_{j=1}^{i+1} \sum_{\substack{J  \subset \{1,\ldots,s\} \\ \card(J) = j }} b_{i-j+1}(S_J,\F).$$
\item For $0 \leq i \leq k'$, $$b_i(S_{\leq s},\F)\leq \sum_{j=1}^{k'-i} \sum_{\substack{J  \subset \{1,\ldots,s\} \\ \card(J) = j }} b_{i+j-1}(S^J,\F)+\binom{s}{k'-i}b_{k'}(S^\emptyset,\F).$$
\end{enumerate}
\end{proposition}

\subsubsection{Smith inequality}

Let $X$ be a compact space  (or a regular complex)  equipped with an involution map $c:X \rightarrow S$. Let 
$\Fix(c) \subset X$ denote the subspace of fixed points of $X$. The Smith exact sequence
(see for example 
\cite[page 126]{Bredon-book})
then implies that
\begin{eqnarray}
\label{eqn:Smith}
b(\Fix(c),\Z_2) &\leq& b(X,\Z_2).
\end{eqnarray}

Taking the involution $c$ to be the complex conjugation we obtain:

\begin{theorem}[Smith inequality for affine sub-varieties of $\C^k$ defined over $\R$]
\label{thm:Smith}
Let $\mathcal{Q}\subset \mathbb{R}[X_1,\ldots,X_{k}]$ be a finite set of polynomials. Then,
\begin{eqnarray*}
b(\ZZ(\mathcal{Q},\R^k),\Z_2) &\leq & b(\ZZ(\mathcal{Q},\C^k),\Z_2).
\end{eqnarray*}
\end{theorem}

\begin{proof}[Proof of Theorem \ref{thm:Smith}]
See \S \ref{subsec:proofs} in the Appendix.
\end{proof}

\subsubsection{Descent spectral sequence}

The following theorem proved in \cite{GVZ04} allows one to bound the Betti numbers of the image of a
closed and bounded semi-algebraic set $S$ under a polynomial map $\mathbf{F}$ in terms of the Betti numbers of  the iterated fibered product of $S$ over $\mathbf{F}$. More precisely:
 
\begin{theorem} \cite{GVZ04}
\label{thm:descent}
Let $S \subset \R^k$ be a closed and bounded semi-algebraic set, and $\mathbf{F} = (F_1,\ldots,F_m):\R^k \rightarrow \R^m$ be a polynomial map. For for all $p, 0 \leq p \leq m$,
\[
b_p(\mathbf{F}(S),\Z_2) \leq \sum_{\substack{i,j \geq 0\\i+j = p}} b_i(\underbrace{S \times_{\mathbf{F}} \cdots \times_{\mathbf{F}} S}_{(j+1)},\Z_2).
\]
\end{theorem}

\subsection{Mixed volume}
Mixed volumes of (Newton) polytopes of polynomials play a very important role in the role of
toric varieties, and will play an important role in this paper (see Theorem \ref{thm:Khovansky} below).
We recall here the definition and certain elementary property of mixed volumes that we will need later in the paper
referring the reader to \cite[\S A.4]{Oda} for any missing detail.

\subsubsection{Definition of mixed volume}

\begin{definition}[Mixed volume]
\label{def:mixed-volume}
Given compact, convex sets $K_1,\ldots,K_m \subset \re^m$, and $\lambda_1,\ldots,\lambda_m \geq 0$,
$\lambda_1K_1 + \cdots + \lambda_mK_m$ is also a compact, convex subset of $\re^m$, and 
$\frac{1}{m!}\vol_m( \lambda_1K_1 + \cdots + \lambda_mK_m)$ is given by a polynomial in $\lambda_1,\ldots,\lambda_m$. The coefficient of $\lambda_1\cdots \lambda_m$ in the polynomial
$\frac{1}{m!}\vol_m( \lambda_1K_1 + \cdots + \lambda_mK_m)$ is called the \emph{mixed volume}
of $K_1,\ldots,K_m$ and denoted by $\MV(K_1,\ldots,K_m)$.
\end{definition}

We will use a few basic properties of mixed volume that we list below (see \cite[\S A.4]{Oda} for an exposition).

\begin{enumerate}
\item(Linearity)
\label{itemlabel:mixed-volume-1}
$$
\displaylines{
\MV(K_1,\ldots,K_{i-1},\lambda' K_i' + \lambda'' K_i'', K_{i+1},\ldots,K_m) = \cr
\lambda' \MV(K_1,\ldots,K_i',\ldots,K_m) + \lambda'' \MV(K_1,\ldots,K_i'',\ldots,K_m).
}
$$
\item (Monotonicity)
\label{itemlabel:mixed-volume-2}
  $K_i' \subset K_i''$ implies that
\[
\MV(K_1,\ldots,K_{i-1},K_i', K_{i+1},\ldots,K_m) \leq 
\MV(K_1,\ldots,K_{i-1},K_i'', K_{i+1},\ldots,K_m).
\]
\item
\label{itemlabel:mixed-volume-3}
If $K_1 =\cdots = K_m = K$, then
\[
\MV(K_1,\ldots,K_m) = \vol_m(K). 
\]
\end{enumerate}

Since in many of our applications we will be interested in obtaining upper bounds on mixed-volumes of certain
special polytopes -- namely products of simplices or boxes, the following simple consequences of Properties
\eqref{itemlabel:mixed-volume-1}, \eqref{itemlabel:mixed-volume-2}, and \eqref{itemlabel:mixed-volume-3} will be useful.

\begin{lemma}
\label{lem:mixed-volume-multiplicative}
For $a_i \geq 0, 1\leq i \leq m$, let 
\[
K_i = \underbrace{\{0\} \times \{0\}}_{i-1} \times [0,a_i] \times \underbrace{\{0\} \times \{0\}}_{m-i-1}.
\]
Then,
\[
\MV(K_1,\ldots,K_m) = \frac{a_1 \cdots a_m}{m!}.
\]
\end{lemma}

\begin{proof}
First observe that
\begin{eqnarray*}
\frac{1}{m!}\vol_m(\lambda_1 K_1 +\cdots +\lambda_m K_m) &=& 
\frac{1}{m!}\vol_m([0,\lambda_1 a_1] \times \cdots \times [0,\lambda_m a_m]) \\
&=& \left(\frac{a_1\cdots a_m}{m!}\right) \lambda_1\cdots\lambda_m.
\end{eqnarray*}

It then follows from Definition \ref{def:mixed-volume} that 
\[
\MV(K_1,\ldots,K_m) = \frac{a_1 \cdots a_m}{m!}.
\]
\end{proof}

\begin{proposition}
\label{prop:mixed-volumes-of-boxes}
Let $B_1,\ldots,B_\ell \subset \re^k, 1 \leq \ell \leq k$, and for $1 \leq i \leq \ell$,
$B_i = [0,d_{i,1}]\ \times \cdots \times [0,d_{i,k}]$. We denote by 
$\mathbf{d} = \left(d_{i,j}\right)_{\substack{1 \leq i \leq\ell \\ 1 \leq j \leq k}}$.
Let $\alpha_1,\ldots,\alpha_\ell  \in \Z_{> 0}$, with $\sum_{i=1}^{\ell} \alpha_i = k$, and
denote $\boldsymbol{\alpha} = (\alpha_1,\ldots,\alpha_\ell)$.
Let 
\[
N(\mathbf{d},\boldsymbol{\alpha}) = \MV(\underbrace{B_1,\ldots,B_1}_{\alpha_1},\ldots, \underbrace{B_\ell,\ldots,B_\ell}_{\alpha_\ell}).
\]
Then, 
\begin{eqnarray}
\label{eqn:mixed-volume-of-boxes-refined}
N(\mathbf{d},\boldsymbol{\alpha})  &=& 
\sum_{\substack{ \mathbf{A}  =(a_{ij}) \in \{0,1\}^{\ell \times k} \\ \sum_{1 \leq j \leq k} a_{ij} = \alpha_i, 1\leq i \leq \ell\\
\sum_{1 \leq i \leq \ell}a_{ij} = 1, 1\leq j \leq k} } \mathbf{d}^{\mathbf{A}}.
\end{eqnarray}

Denoting by $\binom{[1,k]}{\boldsymbol{\alpha}}$ the set of all partitions of $[1,k]$ into disjoint
subsets $J_1,\ldots,J_\ell$ with $\card(J_i) = \alpha_i, 1\leq i \leq \ell$,
\begin{eqnarray}
\label{eqn:mixed-volume-of-boxes-coarse}
N(\mathbf{d},\boldsymbol{\alpha})  &\leq& 
\max_{(J_1,\ldots,J_\ell) \in \binom{[1,k]}{\boldsymbol{\alpha}}} \left(
\prod_{\substack{1\leq i \leq \ell \\ j \in J_i}} d_{i,j} 
\right).
\end{eqnarray}

In the special case, when for $1 \leq i \leq \ell, 1 \leq j \leq k$, $d_{i,j} = d_i$, 
\begin{eqnarray}
\label{eqn:mixed-volume-of-boxes-simple}
N(\mathbf{d},\boldsymbol{\alpha})  &=& d_1^{\alpha_1} \cdots d_\ell^{\alpha_\ell}.
\end{eqnarray}
\end{proposition}

\begin{proof}
Eqns. \eqref{eqn:mixed-volume-of-boxes-refined} and \eqref{eqn:mixed-volume-of-boxes-coarse}
follow from Lemma \ref{lem:mixed-volume-multiplicative} and properties 
\eqref{itemlabel:mixed-volume-1}, \eqref{itemlabel:mixed-volume-2},  and 
\eqref{itemlabel:mixed-volume-3} of mixed volumes stated previously.
Eqn. \eqref{eqn:mixed-volume-of-boxes-simple} is a simple consequence of 
properties 
\eqref{itemlabel:mixed-volume-1}  and 
\eqref{itemlabel:mixed-volume-3} of mixed volumes stated earlier.
\end{proof}

\begin{corollary}
\label{cor:mixed-volumes-of-boxes}
With the same notation as in Proposition \ref{prop:mixed-volumes-of-boxes}, denoting for
$1 \leq i \leq \ell$, 
\[
d_i = \max_{ J \subset [1,k], \card(J) = \alpha_i} \prod_{j \in J} d_{i,j},
\]
\begin{eqnarray}
\label{eqn:mixed-volume-of-boxes-very-coarse}
N(\mathbf{d},\boldsymbol{\alpha})  &\leq& d_1^{\alpha_1} \cdots d_\ell^{\alpha_\ell}.
\end{eqnarray}
\end{corollary}

\begin{proof}
Immediate from Eqn. \eqref{eqn:mixed-volume-of-boxes-refined} in Proposition
\ref{prop:mixed-volumes-of-boxes}.
\end{proof}

\subsection{Topology of complex varieties}
\label{subsec:topology-of-complex}

\subsubsection{Euler-Poincar\'e characteristics of generic intersections in $\C^n$}
\label{subsubsec:E-P-of-complex}
In this section we recall a fundamental result due to Khovanski{\u\i}  \cite{Khovansky78} that we will
exploit heavily  later in the paper. This result in conjunction with Theorem 
\ref{thm:Smith} (Smith inequality) allows us to bound the Betti numbers of generic algebraic varieties
in $\R^k$ in terms of the Newton polytopes of the defining polynomials (under a weak hypothesis on the
Newton polytopes stated in Property \ref{property:khovansky} below).
 
Before recalling  Khovanski{\u\i}'s result we first introduce some more notation.

\begin{notation}
\label{not:supp}
Let $\kb$ be any field. For $P = \sum_{\pmb{\alpha} \in \N^k} c_{\pmb{\alpha}} \X^{\pmb{\alpha}} \in \kb[X_1,\ldots,X_k]$, we denote by $\supp(P) \subset \Q^k$ the 
convex hull of the set $\{\pmb{\alpha} \in \N^k \;\mid\; c_{\pmb{\alpha}} \neq 0\}$. 
\end{notation}

\begin{property}
\label{property:khovansky}
Given a tuple $\pmb{\Delta}= (\Delta_1,\ldots,\Delta_\ell)$, where for $i=1,\ldots,\ell$, $\Delta_i \subset \Q^k$ is a convex polytope, we say that $\pmb{\Delta}$ satisfies
Property \ref{property:khovansky} if for each non-empty subset $L \subset [1,\ell]$, $\dim(\sum_{i \in L} \Delta_i)$ is at least $k-\ell +\card(L)$. We say that 
a tuple of polynomials $\mathcal{P} = (P_1,\ldots,P_\ell), P_i \in \kb[X_1,\ldots,X_k]$ satisfies the same property if the tuple
$\supp(\mathcal{P}) = (\supp(P_1),\ldots,\supp(P_\ell))$ satisfies the above property. 
\end{property}

The following two special cases where Property \ref{property:khovansky} holds will be important later and 
we record this fact here.

\begin{remark}
\label{rem:khovansky}
Notice that if each $\Delta_i$ is a standard $k$-dimensional simplex in $\re^k$ 
of side length  $d_i$ (i.e. the convex hull of
$\mathbf{0}, (d_i,0\ldots,0),\ldots, (0,\ldots,0,d_i)$, with $d_i>0)$, then the tuple 
$(\Delta_1,\ldots,\Delta_\ell)$ satisfies Property \ref{property:khovansky}. The same is true if each
$\Delta_i = [0,d_{i,1}] \times \cdots \times [0,d_{i,k}]$, where 
$\mathbf{d} = (d_{i,j})_{\substack{1\leq i \leq \ell\\1 \leq j \leq k}} \in \Z_{>0}^{\ell \times k}$.
\end{remark}

We first need a notation.

\begin{notation}
\label{not:khovansky}
For $\boldsymbol{\alpha} \in \Z_{\geq 0}^p$, and polytopes $\Delta_1,\ldots,\Delta_p \subset \re^k$, and
any monomial $M(X_1,\ldots,X_p) = \X^{\boldsymbol{\alpha}}$ of degree $k$, 
\begin{eqnarray*}
M(\Delta_1,\ldots,\Delta_p) = k! \; \MV(\underbrace{\Delta_1,\ldots,\Delta_1}_{\alpha_1},\ldots, \underbrace{\Delta_p,\ldots,\Delta_p}_{\alpha_p}),
\end{eqnarray*}
and the definition is extended to any form $H \in \Q[X_1,\ldots,X_p]$ of degree $k$,
by linearity. Finally,  for any rational function $F(X_1,\ldots,X_p)$, 
we define 
\begin{eqnarray}
\label{eqn:not:khovansky}
F(\Delta_1,\ldots,\Delta_p) &=& F_k(\Delta_1,\ldots,\Delta_p),
\end{eqnarray} where 
$F_k$ is the degree $k$ homogeneous component of the Taylor expansion of $F$ at $\mathbf{0}$.
\end{notation}

\begin{theorem}\cite{Khovansky78}
\label{thm:Khovansky}
Let $\mathcal{P} = (P_1,\ldots,P_\ell)$, where each $P_j \in \C[X_1,\ldots,X_k]$, and such that 
$\mathcal{P}$ satisfies Property \ref{property:khovansky}, and the coefficients of the polynomials $P_j$ are generic.
Let $V = \ZZ(\mathcal{P},\C^k)$. Then,
\begin{eqnarray}
\label{eqn:thm:Khovansky}
\chi(V) =  \sum_{I \subset [1,k]}  \prod_{j=1}^\ell \frac{\Delta_j^I}{1+\Delta_j^I},
\end{eqnarray}
where $\Delta_j^I$ is the face of $\Delta_j$ obtained by setting $X_i=0$ for all $i \in I$
(cf. Notation \ref{not:khovansky} and Remark \ref{rem:thm:Khovansky} below).
\end{theorem}

\begin{remark}
\label{rem:thm:Khovansky}
The right hand side of Eqn. \eqref{eqn:thm:Khovansky} needs some explanation. For $I \subset [1,k]$, each simplex $\Delta_j^I$ is a polytope in 
$\mathbb{R}^{k - \card(I)}$, and the expression
\[
\prod_{j=1}^\ell \frac{\Delta_j^I}{1+\Delta_j^I}
\] 
being a rational function of these polytopes represents  a rational number using Eqn. \eqref{eqn:not:khovansky} in  Notation \ref{not:khovansky} (with $k$ being replaced by $k - \card(I)$).
\end{remark}

\subsubsection{Betti numbers of smooth complete intersections in $\C^k$}
\label{subsubsec:Khovansky-applications}

The following proposition appears in \cite{Benedetti-Loeser}  
relates the Euler-Poincar\'e characteristic of affine varieties defined by generic systems of polynomials whose Newton polytopes satisfies Property \ref{property:khovansky}.

\begin{proposition}\cite[Lemma 3.5]{Benedetti-Loeser}.
\label{prop:chi-to-b}
Let $\mathcal{P} = \{P_1,\ldots,P_\ell\}$, where each 
\[P_j \in \C[X_1,\ldots,X_k],
\] 
and such that 
$\mathcal{P}$ satisfies Property \ref{property:khovansky}, and the coefficients of the polynomials $P_j$ are generic.
Let $V = \ZZ(\mathcal{P},\C^k)$. Then,
\begin{equation}
\label{eqn:chi-to-b}
b(V,\Z_2)=   1+(-1)^{k-\ell+1}+(-1)^{k-\ell}\chi(V,\Z_2).
\end{equation}
\end{proposition}

\begin{proof}
See \cite[Lemma 3.5]{Benedetti-Loeser}.
\end{proof}

A similar result also holds for non-singular projective complete intersection varieties. This is
very well known and we sketch a proof.
\begin{proposition}
\label{prop:chi-to-b-projective}
Let $\mathcal{P} = \{ P_1,\ldots,P_\ell \} \subset \C[X_0,\ldots,X_k]$ be a set of generic 
homogeneous polynomials and $\ell \leq k$.
Let $V = \ZZ(\mathcal{P},\PP_\C^k)$.
Then, 
\begin{eqnarray}
\label{eqn:prop:chi-to-b-projective}
b(V,\Z_2) &=& (1 + (-1)^{k-\ell+1})\cdot(k-\ell+1) +(-1)^{k-\ell} \cdot \chi(V,\Z_2).
\end{eqnarray}
\end{proposition}

\begin{proof}
It follows from the Lefschetz hyperplane section theorem and Poincar\'e duality that
\begin{eqnarray}
\label{eqn:Lefschetz-Poincare}
b_i(V,\Z_2) &=& b_i(\PP^{k-\ell}_\C,\Z_2), i \neq k-\ell.
\end{eqnarray}
The proposition follows immediately from \eqref{eqn:Lefschetz-Poincare}, and the fact that
\begin{eqnarray*}
b_i(\PP^{k-\ell}_\C,\Z_2) &=& 1, \mbox{ if $i$ is even and $0 \leq i \leq 2(k-\ell)$, and} \\
&=& 0,  \mbox{ otherwise}.
\end{eqnarray*}
\end{proof}

We stated previously that one of our main tools that we are going to exploit heavily is Theorem \ref{thm:Khovansky}, and this what we proceed to do in this section.
We  use Theorem \ref{thm:Khovansky} in conjunction with Proposition \ref{prop:chi-to-b} and Theorem \ref{thm:Smith}
(Smith inequalities)  to obtain bounds
on the sum of the ($\mathbb{Z}_2$) Betti numbers of certain generic affine complete intersection 
sub-varieties of  $\C^k$ and $\R^k$ that are of interest to us. These include varieties defined by generic
polynomials having prescribed total degrees, or multi-degrees, or with a fixed block of variables appearing at most quadratically with the remaining having prescribed degrees, etc. These results are stated separately
since the precise calculations and the bounds obtained in each case is different (though the main idea
used to obtain these bounds is the same).
   
\subsection{Some applications of Khovanski{\u\i}'s theorem}
\label{subsec:applications-khovansky}

In this section we use Theorem \ref{thm:Khovansky} to obtain bounds on the Betti numbers of generic 
affine intersections in several cases of interest to us. Since some of the calculations are long and technical,
for the sake of readability, we defer the proofs of some of the propositions to  \S \ref{subsec:proofs} in the Appendix.

We begin by observing that as a special case, when $\ell = k$,  Theorem \ref{thm:Khovansky} gives
us a theorem of Bernstein and Kouchnirenko, namely:

\begin{proposition} \cite{Kouchnirenko}
\label{prop:Kouchnirenko}
Let $\mathcal{P} = \{P_1,\ldots,P_k\}\subset \C[X_1,\ldots,X_k]$ be a finite set such that 
$\mathcal{P}$ satisfies Property \ref{property:khovansky}, and the coefficients of the polynomials $P_i$ are 
generic (cf. Definition \ref{def:generic})
Then, $\ZZ(\mathcal{P},\C^k)$ is a finite set, and 
\[
\card(\ZZ(\mathcal{P},\C^k)) = k! \;  \MV(\Delta_1,\ldots,\Delta_k),
\]
where $\Delta_i = \supp(P_i), 1\leq i \leq k$.
Moreover, if additionally $\mathcal{P} \subset \R[X_1,\ldots,X_k]$, then
\[
\card(\ZZ(\mathcal{P},\R^k)) \leq k! \; \MV(\Delta_1,\ldots,\Delta_k).
\]
\end{proposition}

\begin{proof}
Immediate from Theorem \ref{thm:Khovansky}.
\end{proof}

Proposition \ref{prop:Kouchnirenko} deals with the generic zero-dimensional case. Another result that
follows immediately from Theorem \ref{thm:Khovansky} is the following well known expression giving the sum of the Betti numbers of a generic affine hypersurface in $\C^k$ defined by one polynomial of degree $d$. Note that this
proposition could be also be deduced using an argument involving counting multiplicities of Milnor fibers
(see for example, \cite[page 152]{Dimca-book}) or as a special case of Proposition  \ref{prop:many-total-diff-degrees} proved below.

\begin{proposition}
\label{prop:one-total}
Let $P\in \C[X_1,\ldots,X_k], k>0$ be a generic polynomial of total degree $d$.
Then, 
\begin{eqnarray}
\label{eqn:one-total-complex}
b(\ZZ(P,\C^k),\Z_2) &=&  1 + (d-1)^k.
\end{eqnarray}

If $P \in \R[X_1,\ldots,X_k]$, 
\begin{eqnarray}
\label{eqn:one-total-real}
b(\ZZ(P,\R^k),\Z_2) &\leq&  1 + (d-1)^k.
\end{eqnarray}
\end{proposition}

\begin{proof}
See  \S \ref{subsec:proofs} in the Appendix.
\end{proof}

\begin{remark}
\label{rem:one-total}
Notice that if $P \in \R[X_1,\ldots,X_k]$ and defines a bounded, non-singular hypersurface in $\R^k$,
then Proposition \ref{prop:one-total} gives a better bound than just counting critical points of a linear functional on $\ZZ(P,\R^k)$. The latter gives a bound of $d(d-1)^{k-1} > 1 + (d-1)^{k}$ for all $d,k  >1$,
since 
\begin{eqnarray*}
d(d-1)^{k-1} - (1+ (d-1)^{k}) &=& (d-1)^{k-1} -1 \\
&>& 0, \mbox{ for all } d,k > 1.
\end{eqnarray*}
\end{remark}

We now consider the case of generic affine intersections of hypersufaces defined by generic polynomials of
possibly different degrees.

\begin{proposition}
\label{prop:many-total-diff-degrees}
Let $\mathcal{P} = \{P_1,\ldots,P_\ell\} \subset \C[X_1,\ldots,X_k], k\geq \ell >0$ be  a set of generic polynomials with $\deg(P_i) = d_i$.
Then, 
\begin{equation}
\label{eqn:many-total-diff-degrees-complex-1}
b(\ZZ(\mathcal{P},\C^k),\Z_2) = 1 + (-1)^{k- \ell+1} +  d_1\cdots d_\ell\cdot \left( \sum_{j=0}^{k-\ell } 
(-1)^{k-\ell+j}
\binom{k}{j+\ell} h_{j}(d_1,\ldots,d_\ell) \right),
\end{equation}
where
\begin{eqnarray}
\label{eqn:complete-homogeneous}
h_j(d_1,\ldots,d_\ell) =  \sum_{\substack{\boldsymbol{\alpha} \in \Z_{\geq 0},\\ |\boldsymbol{\alpha}|=j}} \dd^{\boldsymbol{\alpha}}
\end{eqnarray}
is the complete homogeneous symmetric polynomial of degree $j$ in $\dd=(d_1,\ldots, d_\ell)$.

Moreover, if $d_1 = \cdots = d_\ell = d$  then
\begin{equation}
\label{eqn:many-total-diff-degrees-complex-2}
b(\ZZ(\mathcal{P},\C^k),\Z_2) \leq 1+(-1)^{k-\ell+1}+ \binom{k-1}{\ell-1} (d^k  + k - 1) \mbox{ if } \ell <k,
\end{equation}

\begin{equation}
\label{eqn:many-total-diff-degrees-complex-2'}
b(\ZZ(\mathcal{P},\C^k),\Z_2) \leq d^k \mbox{ if } \ell = k.
\end{equation}

Additionally, if $\mathcal{P} \subset \R[X_1,\ldots,X_k]$, then
\begin{equation}
\label{eqn:many-total-diff-degrees-real-1}
b(\ZZ(\mathcal{P},\R^k),\Z_2) \leq  1 + (-1)^{k- \ell+1} +  d_1\cdots d_\ell\cdot \left( \sum_{j=0}^{k-\ell }
(-1)^{k-\ell+j}
\binom{k}{j+\ell} h_{j}(d_1,\ldots,d_\ell) \right).
\end{equation}

In the case $d_1 = \cdots = d_\ell = d$,
\begin{equation}
\label{eqn:many-total-diff-degrees-real-2}
b(\ZZ(\mathcal{P},\R^k),\Z_2) \leq 1+(-1)^{k-\ell+1}+ \binom{k-1}{\ell-1} (d^k  + k - 1),
\end{equation}
if $\ell < k$, and
\begin{equation}
\label{eqn:many-total-diff-degrees-real-2'}
b(\ZZ(\mathcal{P},\R^k),\Z_2) \leq  d^k,
\end{equation}
if $\ell=k$.
\end{proposition}

\begin{proof}
See  \S \ref{subsec:proofs} in the Appendix.
\end{proof}

\begin{remark}
\label{rem:why-affine}
We remark that Eqn. \eqref{eqn:many-total-diff-degrees-complex-1}  in Proposition \ref{prop:many-total-diff-degrees} below can also be obtained by a Chern class computation
and an argument using Lefschetz duality.
(In view of its importance we include this alternate proof in  \S\ref{subsec:Chern} in the Appendix.)
However, this approach is not applicable in
many of the situations considered in the current paper since the generic affine intersections that we
consider might necessarily be singular at infinity.
For example, let $k=3$, and $\Delta = [0,1] \times [0,1] \times [0,1]$. Then, a generic polynomial 
$P  \in \C[X_1,X_2,X_3]$ with
$\supp(P) \subset \Delta$, defines a non-singular hypersurface in $\C^3$, but defines a singular curve in the projective plane at infinity defined by the (homogeneous) equation
$X_1X_2X_3=0$ (with three singular points -- namely, $(1:0:0),(0:1:0),(0:0:1)$). 
One could take more complicated compactifications of affine space -- for example, multi-projective spaces -- which would solve this problem, at the cost of increasing the complexity of the process of subtracting the added part, a process which moreover would be 
different in each of the cases that we consider.
Because of these reasons it is convenient for us to have directly an expression for the Betti numbers of generic complex \emph{affine} intersections -- which is afforded by Theorem \ref{thm:Khovansky} in conjunction with Proposition \ref{prop:chi-to-b}. 
\end{remark}

\begin{remark}
\label{rem:degree2andd}
A special case of Proposition \ref{prop:many-total-diff-degrees} will be used later, and we record it here 
for future use. Let $\mathcal{P} \subset \R[X_1,\ldots,X_k]$ be generic and equal to the disjoint union of
$\mathcal{P}_1$ and $\mathcal{P}_2$, with $\card(\mathcal{P}_i) = \ell_i$, and the total
degrees of the polynomials in $\mathcal{P}_i$ equal $d_i$, $i=1,2$. 
Applying Proposition \ref{prop:many-total-diff-degrees} to this special case we obtain that
$b(\ZZ(\mathcal{P},\R^k),\Z_2)$ is bounded by 
\begin{eqnarray}
\label{eqn:degreed1andd2}
\nonumber
&& 1 + (-1)^{k- \ell+1} + d_1^{\ell_1} d_2^{\ell_2} \cdot\left( \sum_{j=0}^{k-\ell } (-1)^{k-\ell+j}\binom{k}{j+\ell} h_{j}(\underbrace{d_1,\ldots,d_1}_{\ell_1},\underbrace{d_2,\ldots,d_2}_{\ell_2}) \right) \\\nonumber
&=& 1 + (-1)^{k- \ell+1} + d_1^{\ell_1}d_2^{\ell_2} \cdot\left( \sum_{j=0}^{k-\ell } \sum_{i=0}^{j}(-1)^{k-\ell+j}\binom{k}{j+\ell}\binom{\ell+i-2}{i}d_1^{j-i}d_2^i\right) \\\nonumber
& \leq &  \binom{k-2}{\ell -2} d_1^{\ell_1} d_2^{k-\ell_1} + O(1)^k d_1^{\ell_1+1}d_2^{k-\ell_1-1}.
\end{eqnarray}
In particular if $\ell_1 = 1$, and $d_1=2$, 
$b(\ZZ(\mathcal{P},\R^k),\Z_2)$ is bounded by

\begin{eqnarray}
\label{eqn:degree2andd}
&& 1 + (-1)^{k- \ell+1} + 2d^{\ell-1} \cdot\left( \sum_{j=0}^{k-\ell } \sum_{i=0}^{j}(-1)^{k-\ell+j}\binom{k}{j+\ell}\binom{\ell+i-2}{i}2^{j-i}d^i\right) \\
\nonumber
& \leq &  2 \binom{k-2}{\ell -2} d^{k-1} + (O(d)) ^{k-2}.
\end{eqnarray}
\end{remark}

\begin{proposition}
\label{prop:manyblocks-total}
Let $\mathcal{P} \subset \C[\X^{(1)},\ldots,\X^{(p)}]$, 
for $1 \leq i \leq p$, $\X^{(i)} = (X^{(i)}_1,\ldots,X^{(i)}_{k_i})$ and $\deg_{\X^{(i)}}(P) \leq d_i, P \in \mathcal{P}$, with $\ell = \card(\mathcal{P}) > 0$.
Let $\kk = (k_1,\ldots,k_p),k = \sum_{i=1}^p k_i$  and $V_\C =  \ZZ(\mathcal{P},\C^k)$.
Suppose also that the polynomials in $\mathcal{P}$  are generic. 

Then,
\begin{equation}
\label{eqn:manyblocks-total-complex}
b(V_\C,\Z_2) \leq 1 + (-1)^{k- \ell+1} +  (k-\ell+2)^2 \binom{k}{\ell-1}\binom{k}{\kk}^{-1}\frac{(1+p)^{3k-\ell+1}}{p(p+2)} d_1^{k_1}\cdots d_p^{k_p}.
\end{equation}

If additionally, $\mathcal{P} \subset \R[\X^{(1)},\ldots,\X^{(p)}]$, and  $V_\R=  \ZZ(\mathcal{P},\R^k)$,  then
\begin{equation}
\label{eqn:manyblocks-total-real}
b(V_\R,\Z_2) \leq 1 + (-1)^{k- \ell+1} +  (k-\ell+2)^2\binom{k}{\ell-1}\binom{k}{\kk}^{-1} \frac{(1+p)^{3k-\ell+1}}{p(p+2)} d_1^{k_1}\cdots d_p^{k_p}.
\end{equation}
\end{proposition}

\begin{proof}
See  \S \ref{subsec:proofs} in the Appendix.
\end{proof}

In order to investigate the tightness of the inequalities in Proposition \ref{prop:manyblocks-total}, it 
is instructive to consider the 
special case of Proposition \ref{prop:manyblocks-total} when the block sizes are all equal to one (i.e. $p=k$) and $\ell=1$.

\begin{proposition}
\label{prop:one-multi}
 Let $P$ be a generic  polynomial in $\C[X_1,\ldots,X_{k}]$, 
and $\deg_{X_i}(P) \leq d_i$ with $d_1 \geq d_2 \geq \cdots \geq d_k \geq 0$. 
We denote by $\bar{d} = (d_1,\ldots,d_{k})$, and for $J \subset [1,k]$,
we denote $\bar{d}^{J} = \prod_{i \in J} d_i$.  Then,
\begin{equation}
\label{eqn:one-multi-complex}
b(\ZZ(P,\C^k) ,\Z_2) = 1 + (-1)^k + \left(  \sum_{j=1}^{k} (-1)^{k-j} \sum_{\substack{J  \subset [1,k] \\ \card(J) = j \leq k }} j! \bar{d}^{J}  \right),
\end{equation}
and if $P \in \R[X_1,\ldots,X_k]$, then 
\begin{equation}
\label{eqn:one-multi-real}
b(\ZZ(P,\R^k) ,\Z_2) \leq 1 + (-1)^k + \left(  \sum_{j=1}^{k} (-1)^{k-j} \sum_{\substack{J  \subset [1,k] \\ \card(J) = j \leq k }} j! \bar{d}^{J}  \right).
\end{equation}
\end{proposition}

\begin{proof}
See  \S \ref{subsec:proofs} in the Appendix.
\end{proof}

\begin{remark}[Comparison with the 
Ole{\u\i}nik-Petrovski{\u\i}-Thom-Milnor bound]
\label{rem:one-multi1}
Theorem \ref{thm:O-P-T-M} gives a
bound of $D(2D-1)^{k-1}$ in the context of Proposition \ref{prop:one-multi} with the total degree $D = d_1 + \cdots + d_k$. This bound is in general much worse than the bound in inequality
\eqref{eqn:one-multi-real} in Proposition \ref{prop:one-multi}.
For example, take $k=2$, 
$\bar{d} = (d,d)$. Then the 
bound from Theorem \ref{thm:O-P-T-M} (i.e. the Ole{\u\i}nik-Petrovski{\u\i}-Thom-Milnor bound)
is $2d(4d-1)$, while Proposition \ref{prop:one-multi} yields a bound of 
\[
1 + 1 - 2d + 2 d^2 = 2d^2 -2d+ 2< 2d(4d-1)
\] 
for all $d>0$. 
\end{remark}

\begin{remark}
\label{rem:one-multi2}
Proposition \ref{prop:one-multi} is tight when $k_1 = k_2 =1$ and $\bar{d} = (2,2)$. Then, the bound in Proposition \ref{prop:one-multi}  is (using the formula in Remark \ref{rem:one-multi1})
\[
2\cdot 2^2 - 2\cdot 2 + 2 = 6.
\]
Consider the polynomial
\[
P_\eps = (X_1^2 -1)(X_2^2 -1) - \eps,
\]
and let $V = \ZZ(P_\eps,\R^2)$.
Then, for all sufficiently small $\eps > 0$, $b_0(V,\Z_2) = 5$, and $b_1(V,\Z_2) = 1$ (see Figure \ref{fig:one-multi}), so that
$b(V,\Z_2)=6$.

Notice that the Ole{\u\i}nik-Petrovski{\u\i}-Thom-Milnor bound of $d(2d-1)^{k-1}$, where $d$ is the total degree, yields in this case
$4\cdot(8-1) =  28$,  which is much worse than the bound in Proposition \ref{prop:one-multi}.

\begin{figure}
\includegraphics[scale=0.3]{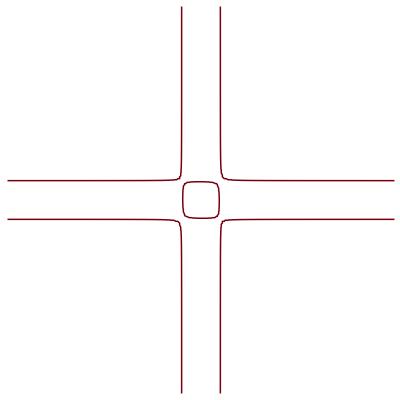}
\caption{The set of real zeros in $\R^2$ of  $P_\eps=(X_1^2 -1)(X_2^2 -1) - \eps$ for small $\eps$.}
\label{fig:one-multi}
\end{figure}
\end{remark}

\begin{proposition}
\label{prop:different-boxes}
Let $B_1,\ldots,B_\ell \subset \re^k, 1 \leq \ell \leq k$, and for $1 \leq i \leq \ell$,
$B_i = [0,d_{i,1}]\ \times \cdots \times [0,d_{i,k}]$. We denote by 
$\mathbf{d} = \left(d_{i,j}\right)_{\substack{1 \leq i \leq\ell \\ 1 \leq j \leq k}}$.
Let $\mathcal{P} =\{P_1,\ldots,P_\ell\} \subset \C[X_1,\ldots,X_k$, with $\supp(P_i) = B_i$,
$V_\C =  \ZZ(\mathcal{P},\C^k)$,
with $\mathcal{P}$ assumed to be generic.
\begin{eqnarray}
\label{eqn:different-boxes-complex}
\nonumber
b(V_\C ,\Z_2) &\leq& 1 + (-1)^{k-\ell+1} + 
\sum_{j=\ell}^{k} \sum_{J  \in  \binom{[1,k]}{j}}
(-1)^{k-j}
\sum_{\substack{
\boldsymbol{\alpha} =(\alpha_1,\ldots,\alpha_\ell) \in \Z_{>0}^\ell\\
\alpha_1+\cdots+\alpha_\ell = j}}  N (\mathbf{d}_J,\boldsymbol{\alpha})
\\
&\leq&
O(\ell)^k\cdot
\max_{\substack{\boldsymbol{\alpha} = (\alpha_1,\ldots,\alpha_\ell) \in \Z_{>0} \\ \alpha_1+\cdots+\alpha_\ell = k}} 
\max_{(J_1,\ldots,J_\ell) \in \binom{[1,k]}{\boldsymbol{\alpha}}} \left(
\prod_{\substack{1\leq i \leq \ell \\ j \in J_i}} d_{i,j} 
\right),
\end{eqnarray}
where for $J \subset [1,k]$, $\mathbf{d}_J$ is the $\ell \times \card(J)$ sub-matrix obtained by extracting the columns corresponding to $J$ in $\mathbf{d}$, and 
$N(\mathbf{d}_J,\boldsymbol{\alpha})$ is defined in 
\eqref{eqn:mixed-volume-of-boxes-refined}.

The same bound as in \eqref{eqn:different-boxes-complex} 
also hold for  $b(\ZZ(\mathcal{P},\R^k) ,\Z_2)$,  if additionally $\mathcal{P}$ has coefficients in $\R$.
\end{proposition}

\begin{proof}
See  \S \ref{subsec:proofs} in the Appendix.
\end{proof}

\subsubsection{Quadratic and partially quadratic case}
\label{subsec:quadratic}
We now use Theorem \ref{thm:Khovansky} to obtain bounds on the Betti numbers of generic
intersections of quadratic and partially quadratic polynomials. 
Since the dependence of the bounds on the different parameters in this case are rather different
from the previous cases,
we start by explaining the most simple case in detail.

\begin{proposition}
\label{eg:two-quadratic}
Let $P_1,P_2$ be two generic quadratic polynomials in $\C[X_1,\ldots,X_k]$.
Then,
 \begin{eqnarray}
\label{eqn:two-quadratic-complex}
b(\ZZ(\{P_1,P_2\},\C^k),\Z_2) &=& 2k.
\end{eqnarray}
\end{proposition}

\begin{proof}
First note that $\supp(P_1),\supp(P_2)$ are both equal to the 
convex hull of $\mathbf{0}, (2,0,\ldots,0),\ldots, (0,\ldots,0,2)$. 
It follows applying Theorem \ref{thm:Khovansky} that for $k \geq 2$,
\begin{eqnarray*}
\chi(\ZZ(\{P_1,P_2\},\C^k),\Z_2) &=& \sum_{j=2}^{k} \binom{k}{j} (-1)^{j} (j-1)j! \frac{2^j}{j!} \\
            &=& 1+ \sum_{j=0}^{k} \binom{k}{j} (-1)^{j} (j-1){2^j} \\
            &=& 1+ \sum_{j=0}^{k} \binom{k}{j} (-1)^{j} j2^j -\sum_{j=0}^{k} \binom{k}{j} (-1)^{j} 2^j \\
            &=&  1+ 2\sum_{j=0}^{k} \binom{k}{j} (-1)^{j} j2^{j-1} - \sum_{j=0}^{k} \binom{k}{j} (-1)^{j} 2^j \\
             &=& 1 + 2k(1-2)^{k-1}(-1) - (1-2)^k \\
             &=&  1 + (-1)^k2k - (-1)^k \\
             &=&  1 + (-1)^k (2k-1). 
\end{eqnarray*}
This implies (using Eqn. \eqref{eqn:chi-to-b})  that 
\begin{eqnarray*}
b(V,\Z_2) 
&=& 1 + (-1)^{k}(\chi(V_k,\Z_2) -1) \\
&=& 1 + (-1)^{k+1} + (-1)^{k} \chi(V_k,\Z_2) \\
&=& 1 + (-1)^{k+1} + (-1)^{k} \left( 1 + (-1)^k (2k-1) \right)\\
&=& 2k.
\end{eqnarray*}
\end{proof}

\begin{remark}
\label{rem:two-quadratic}
In particular, when $k=2$, $b(\ZZ(\{P_1,P_2\},\C^k),\Z_2) = 4$, agreeing with the fact that two generic quadratic polynomials
in two variables will have $4$ points in their intersection. 

In the case $k=3$,  notice that intersection, $W$  of two generic quadric surfaces in $\PP_\C^3$ is topologically a torus $\mathbf{S}^1 \times \mathbf{S}^{1}$. The intersection of $W$ with the plane at infinity is 
$4$ points. Hence, $\ZZ(\{P_1,P_2\},\C^k)$ in this case is homeomorphic to $\mathbf{S}^1 \times \mathbf{S}^{1}$ minus $4$ points.
This gives,
\begin{eqnarray*}
b_0(\ZZ(\{P_1,P_2\},\C^k),\Z_2) &=& 1, \\
b_1(\ZZ(\{P_1,P_2\},\C^k),\Z_2) &=& 5.
\end{eqnarray*}
This gives, $\chi(V,\Z_2) = -4$ which agrees with the formula above.
\end{remark}

We now consider a more general situation.

\begin{proposition}
\label{prop:many-total-mixed}
 Let $\mathcal{P} = \{P_1,\ldots,P_\ell\}$ be a finite set of generic  polynomials in $\C[X_1,\ldots,X_{k_1},Y_1,\ldots,Y_{k_2}]$ with $0<\ell \leq k = k_1+k_2$, 
and $\deg_{\X}(P_i) \leq d,  \deg_\Y(P_i) \leq 2$. 

Then,
\begin{equation}
\label{eqn:many-total-mixed-complex}
b(\ZZ(\mathcal{P},\C^k) ,\Z_2) \leq 2 + (-1)^{k-\ell +1} + \ell 2^\ell (k_1+k_2)^{\ell-1} \left( 2d(k_1+k_2)+1 \right)^{k_1}.
\end{equation}
The same bound as in \eqref{eqn:many-total-mixed-complex} also holds for $b(\ZZ(\mathcal{P},\R^k) ,\Z_2)$,
if additionally $\mathcal{P}$ has coefficients in $\R$.
\end{proposition}

\begin{proof}
See  \S \ref{subsec:proofs} in the Appendix.
\end{proof}

\subsubsection{Generic intersections of quadrics in affine and projective spaces}
Since the intersections of quadrics is a very well studied topic \cite{Agrachev, Lerario2012} we investigate the special case of
 Proposition \ref{prop:many-total-mixed} where $k_1=0$.
 In particular, we calculate the leading coefficient of the polynomial in $k$ giving the sum of the Betti numbers 
 of the intersection of $\ell$ generic quadrics in $\PP_\C^k$ for every fixed $\ell$, thus solving a problem
 posed in \cite{Lerario2012} (see Eqn. \ref{eqn:quadratic-projective-betti-exact} below).
 
Setting $k_1=0$ and $k_2=k$, in the above calculation and keeping the same notation, we obtain that
\begin{eqnarray}
\label{eqn:quadratic-exact-affine}
\chi(\ZZ(\mathcal{P},\C^k),\Z_2) &=&  1+ (-1)^{k+1}
 \left( 
\sum_{h=0}^{\ell-1}
\binom{k}{h} (-2)^h
\right).
\end{eqnarray}

In order to calculate the Euler-Poincar\'e characteristic of a generic complete intersection of dimension $k-\ell$ in $\PP_\C^k$ in terms of a fixed degree sequence $(d_1,\ldots,d_\ell)$, it suffices to take the sum of the Euler-Poincar\'e characteristics of the corresponding
affine varieties in $\C^k,\C^{k-1},\ldots,\C^{\ell}$ (with the same degree sequence). Applying this in our situation we obtain 
that if $\mathcal{P} = \{P_1,\ldots,P_\ell\}$ are generic homogeneous quadrics in $\C[X_0,\ldots,X_k]$ then it follows from the above and
\eqref{eqn:quadratic-exact-affine} that

\begin{eqnarray}
\label{eqn:quadratic-exact-projective}
\nonumber
\chi(\ZZ(\mathcal{P},\PP_\C^k),\Z_2) &=& 
\sum_{j=\ell}^{k}
\left( 1+
 (-1)^{j+1}
\left( 
\sum_{h=0}^{\ell-1}
\binom{j}{h} (-2)^h
\right) 
\right) \\
&=&
\sum_{h=0}^{\ell-1} (-2)^h\left( \sum_{j=\ell}^{k} (-1)^{j+1} \binom{j}{h} \right) +  (k-\ell+1),
\end{eqnarray}
and using Proposition \ref{prop:chi-to-b-projective}
\begin{equation}
\label{eqn:quadratic-projective-betti-exact}
b (\ZZ(\mathcal{P},\PP_\C^k),\Z_2) = 
(1+(-1)^{k-\ell+1})(k-\ell+1) + (-1)^{k-\ell}\chi(\ZZ(\mathcal{P},\PP_\C^k),\Z_2),
\end{equation}
for $\ell < k$.

It  is easy to deduce directly from \eqref{eqn:quadratic-exact-projective} and
\eqref{eqn:quadratic-projective-betti-exact} that for $\ell \geq 3$
\begin{eqnarray}
\label{eqn:quadratic-exact-projective-inequality3} 
b(\ZZ(\mathcal{P},\PP_\C^k),\Z_2) &= & \frac{2^{\ell-2}}{(\ell-1)!} {k}^{\ell-1} +O(k^{\ell-2}), 
\end{eqnarray}
for fixed $\ell$ and $k$ large.

To see this define,
\[
B(h,k,\ell) = 2^h\left( \sum_{j=\ell}^{k} (-1)^{j+1} \binom{j}{h} \right).
\]

$B(\ell-1,k,\ell)$ is the absolute value of the term corresponding to $h=\ell-1$  in the expression in
\eqref{eqn:quadratic-exact-projective}.
We have
\begin{eqnarray}
\label{eqn:Lerario-inequality1}
\nonumber
B(\ell-1,k,\ell) 
&=&2^{\ell-1}\left( \binom{k}{\ell-1} - \binom{k-1}{\ell-1} + \cdots + (-1)^{k-\ell-1}\binom{\ell-1}{\ell-1}\right)\\
&=&
2^{\ell-1}\left( \binom{k-1}{\ell-2} + \binom{k-3}{\ell-2} + \cdots \right).
\end{eqnarray}

\begin{lemma}
\label{lem:Lerario-inequality}
For $p \geq 0$,  and all large $n$,
\begin{eqnarray*}
\sum_{i=0}^{\lfloor(n-p)/2\rfloor } \binom{n -2i}{p} &=& \frac{1}{2} \binom{n+1}{p+1} + O(n^{p}).  
\end{eqnarray*}
\end{lemma}

\begin{proof}
Let 
\[
A(n,p) =  \sum_{i=0}^{\lfloor(n-p)/2\rfloor } \binom{n -2i}{p}.
\]
From standard binomial identities we deduce
\[
A(n,p) + A(n-1,p) = \sum_{i=0}^{\lfloor(n-p)\rfloor } \binom{n -i}{p} = \binom{n+1}{p+1},
\]
and 
\[
A(n,p) - A(n-1,p) = A(n-1,p-1) \leq A(n,p-1).
\]
The lemma now follows by induction on $p$, the case $p=0$ being trivial.
\end{proof}

It follows from 
Lemma \ref{lem:Lerario-inequality} and \eqref{eqn:Lerario-inequality1} that
\begin{eqnarray}
\label{eqn:Lerario-inequality2}
B(\ell-1,k,\ell) &=& 2^{\ell-2}\binom{k}{\ell-1} + O(k^{\ell-2}).
\end{eqnarray}

Moreover,
\begin{eqnarray}
\label{eqn:Lerario-inequality3}
\sum_{h=0}^{\ell-2} B(h,k,\ell) &= &  O(k^{\ell-2}).
\end{eqnarray}

It follows from \eqref{eqn:Lerario-inequality2}, \eqref{eqn:Lerario-inequality3}, Proposition
\ref{prop:chi-to-b-projective} that
\begin{eqnarray}
\label{eqn:Lerario}
b(\ZZ(\mathcal{P},\PP_\C^k),\Z_2)
&=& 
2^{\ell-2}\binom{k}{\ell-1} + O(k^{\ell-2}),
\end{eqnarray}
which implies inequality \eqref{eqn:quadratic-exact-projective-inequality3}.

This answers a question raised in \cite[page 4]{Lerario2012}, where the first few values of the leading coefficients
were given, and the problem of calculating it exactly was posed. 
Notice that this coefficient, $\frac{2^{\ell-2}}{(\ell-1)!}$,  goes to zero exponentially  fast  with $\ell$.

Using  Theorem \ref{thm:Smith} (Smith inequality)  and  \eqref{eqn:quadratic-exact-projective-inequality3}
we obtain the following theorem  (cf. \cite{Lerario2012}).
\begin{theorem}
\label{thm:quadratic-projective}
Let $\mathcal{P} \subset \R[X_0,\ldots,X_k]$ be a set of $\ell$ generic homogeneous polynomials of degree $2$.
Then for every fixed $\ell \geq 3$,  
\begin{eqnarray*}
b(\ZZ(\mathcal{P},\PP_\R^k),\Z_2) &\leq& 2^{\ell-2} \binom{k}{\ell-1} + O(k^{\ell-2}).
\end{eqnarray*}
\end{theorem}

\begin{remark}
\label{rem:quadric-projective}
Note that a more naive approach using a bound of $(O(k))^{\ell-1}$ on the Betti numbers of generic intersections of 
$\ell$ affine quadrics in $\C^k$, the fact that
$\PP_\C^k$ is the disjoint union of $\C^k,\C^{k-1},\ldots,\C^0$, and the additivity property of the
Euler-Poincar\'e characteristics, 
yields a slightly coarser bound of $(O(k))^{\ell}$. The signs thus play an important role in the
proof of Theorem \ref{thm:quadratic-projective}.
Finally, note that the more general version of Theorem \ref{thm:quadratic-projective}, with a slightly worse bound but without any assumption of genericity appears in
Theorem \ref{thm:Lerario-new}.
\end{remark}

\subsubsection{One quadratic block and multi-degree case}
We now consider the case of generic intersections of polynomials having one block of variables of
degree $2$, while the other variables are allowed to have different degrees. More precisely, we prove:

\begin{proposition}
\label{prop:several-blocks-mixed}
Let $\mathcal{P}=\{P_1,\ldots,P_\ell \}\subset \C[X_1,\ldots,X_{k_1},Y_1,\ldots,Y_{k_2}]$ be a finite set of generic polynomials with $\deg_{X_i}(P)\leq d_i$ and $\deg_\Y(P)\leq 2$ for all $P\in \mathcal{P}$ with $d_1 \geq d_2 \geq \cdots \geq d_{k_1}$ and $0<\ell \leq k=k_1+k_2$. 
Then,
\begin{equation}
\label{eqn:several-blocks-mixed-complex}
b(\ZZ(\mathcal{P},\C^k),\Z_2) \leq 
2 + (-1)^{k-\ell+1}+\ell 2^\ell k_1!(k_1+k_2)^{\ell-1}\left( 2(k_1+k_2)+1 \right)^{k_1}d_1\cdots d_{k_1}.
\end{equation}
The same bound as in \eqref{eqn:many-total-mixed-complex} holds for $b(\ZZ(\mathcal{P},\R^k) ,\Z_2)$,
if additionally $\mathcal{P}$ has coefficients in $\R$.
\end{proposition}

\begin{proof}
See  \S \ref{subsec:proofs} in the Appendix.
\end{proof}

\section{Proofs of the main theorems}
\label{sec:main-results}

\subsection{Summary of the methods}
Our main tools are the bounds on the Betti numbers of generic  intersections proved in \S
\ref{subsec:applications-khovansky} above, which are all consequences of Theorem \ref{thm:Khovansky} and Theorem
\ref{thm:Smith} (Smith inequality), the techniques of infinitesimal perturbations (\cite[Chapter 7]{BPRbook2}), and the inequalities derived from the Mayer-Vietoris exact sequence (Proposition
\ref{prop:unionint}). Using the techniques of infinitesimal perturbations, and the inequalities in
Proposition \ref{prop:unionint}, we reduce the problem of bounding the Betti numbers of semi-algebraic sets defined
by general (non-generic) polynomials $P \in \mathcal{P} \subset \R[X_1,\ldots,X_k]$ with support contained in given Newton polytopes $\Delta_P, P \in \mathcal{P}$ (the tuple $(\Delta_P)_{P \in \mathcal{P}}$
satisfying Property \ref{property:khovansky}), to bounding the Betti numbers of a collection of real affine algebraic 
varieties defined by generic polynomials with (nearly) the same support. The proofs of Theorems
\ref{thm:algebraic-total-degree}, \ref{thm:multi-algebraic}, \ref{thm:different-boxes-algebraic},
\ref{thm:partly-quadratic-total-algebraic}, and \ref{thm:partly-quadratic-multi-algebraic} 
(i.e. the cases of different classes of algebraic sets)  are very similar to each other, 
differing only in the application of the appropriate generic bounds. Because of this reason we explain only the proof of Theorem \ref{thm:algebraic-total-degree} in full detail. Similarly, the proofs of Theorems 
\ref{thm:multi-semi}, \ref{thm:different-boxes-semi}, 
\ref{thm:partly-quadratic-total-semi}, and \ref{thm:partly-quadratic-multi-semi} 
(the semi-algebraic cases) are all similar in structure
to the proof of  
\cite[Theorem 7.30]{BPRbook2} and \cite[Theorem 7.38]{BPRbook2}, again 
differing only in the application of the appropriate generic bounds. We refer the reader to
\cite{BPRbook2} for any missing detail.

\subsection{Deformation to generic}
As mentioned above we will  reduce the problem of bounding the Betti numbers of semi-algebraic sets defined
by general (non-generic) polynomials $P \in \mathcal{P} \subset \R[X_1,\ldots,X_k]$  to bounding the Betti numbers of a collection of real affine algebraic 
varieties defined by generic polynomials with  the same support. For this we will need to following technical lemma.

\begin{lemma}
\label{lem:generic}
Let $V \subset \R^N$ be a Zariski closed subset with $\dim_\R(V) < N$. Let $\mathbf{a} \in \R^N \setminus V$. Then, for all $\mathbf{b}\in \R^N$, $(1-\delta)\cdot\mathbf{b} + \delta \cdot\mathbf{a} \not\in \Ext(V,\R\la\delta\ra)$.
\end{lemma}

 \begin{proof}
 Since, $V$ is Zariski-closed, $V = \ZZ(P,\R^N)$ for some $P \in R[X_1,\ldots,X_N]$.
 The polynomial $F(T) := P((1-T)\mathbf{b} + T \mathbf{a}) \in \R[T]$ is not identically $0$, since $F(1) \neq 0$. Hence, there exists $t_0 \in \R$ with $ t_0 > 0$, such that $F$ does not vanish in the interval $(0,t_0)$.
 This implies that $P((1-\delta)\cdot\mathbf{b} + \delta\cdot \mathbf{a}) \neq 0$, and hence
 $(1-\delta)\cdot\mathbf{b} + \delta \cdot\mathbf{a} \not\in \Ext(V,\R\la\delta\ra)$.
 \end{proof} 
 
 \begin{corollary}
 \label{cor:generic}
 For $1 \leq i \leq \ell$, let $\Delta_i \subset \Q^k$ be a convex polytope,
 and let $H_i = \sum_{\pmb{\alpha} \in \Delta_i \cap \N^k} a_{i, \pmb{\alpha}} \X^{\pmb{\alpha}}, P_i = \sum_{\pmb{\alpha} \in \Delta_i \cap \N^k} b_{i, \pmb{\alpha}} \X^{\pmb{\alpha}}, \in \R[X_1,\ldots,X_k], 1  \leq i \leq \ell$.
 Let $V \subset \R^N$ be a real variety, with 
 $N = \sum_{i=1}^{\ell} \card(\Delta_i \cap \N^k)$.
 Let 
 $\mathbf{a} = (\ldots, a_{i,\pmb{\alpha}},\ldots),
  \mathbf{b} = (\ldots,b_{i,\pmb{\alpha}},\ldots) \in \R^N$ denote the vectors of coefficients of $(H_1,\ldots,H_\ell)$ and $(P_1,\ldots,P_\ell)$ respectively. Suppose also that $\mathbf{a} \not\in V$. Then, if $\mathbf{c} \in \R\la\delta\ra^N$ is the vector of coeffcients of $(1-\delta)\cdot P_1 + \delta \cdot H_1,\ldots,(1-\delta)\cdot P_\ell +\delta\cdot H_\ell$, then $\mathbf{c} \not\in \Ext(V,\R\la\delta\ra)$. 
 (In other words, using the language introduced in Definition \ref{def:generic}, if $(H_1,\ldots,H_\ell)$ is a generic tuple of polynomials, so is 
 $(1-\delta)\cdot P_1+\delta \cdot H_1,\ldots, (1-\delta)\cdot P_\ell + \delta\cdot H_\ell$ for all tuples of polynomials $(P_1,\ldots,P_\ell)$ with $\supp(P_i) \subset \Delta_i, 1 \leq i \leq \ell$.)
 \end{corollary}
 
 \begin{proof}
 Apply Lemma \ref{lem:generic}.
 \end{proof}

\subsection{Proof of  Theorem \ref{thm:algebraic-total-degree}}
\label{subsec:proof-total-degree}

\begin{proof}[Proof of Theorem \ref{thm:algebraic-total-degree}]
We first prove 
\begin{equation}
\label{eqn:algebraic-total-degree-tilde-bound-2}
b(V,\Z_2) \leq  \frac{1}{2}(1 + (2d-1)^k).
\end{equation}

Let $r > 0$, and
let 
\[
F(X_1,\ldots,X_k)=(Q_1^2+\cdots+Q_\ell^2)/(r^2-\| \X\|^2).
\] 
The set of critical values of $F$ is finite, so there exists 
$
c_0\in \R,
c_0>0$ so that 
$\ZZ(\tilde{Q},\R^k)$ is a non-singular hypersurface in $\R^k$, where
\begin{equation}
\label{eqn:Qtilde}
\tilde{Q}_r=Q_1^2+\cdots+Q_\ell^2+c (\| \X\|^2-r^2)=0,
\end{equation}
for all $c \in (0,c_0)$.

Denote by $\tilde{V}=\ZZ (\tilde{Q},\R^k)$. 
Since $\deg(\tilde{Q}) \leq 2d$, we have by Proposition \ref{prop:one-total} that 
\begin{equation}
\label{eqn:algebraic-total-degree-tilde-V}
b(\tilde{V}_r,\Z_2) \leq 1+(2d-1)^k.
\end{equation}

Notice that the 
the closed and bounded semi-algebraic set, $S_r$ defined by $\tilde{Q}_r \leq 0$,
is 
semi-algebraically
homotopy equivalent to $V_r = V \cap \overline{B_k(0,r)}$. 

Since $S_r$ is bounded by $\tilde{V}_r$, 
\begin{equation}
\label{eqn:algebraic-total-degree-tilde-V-and-S}
b(S_r,\Z_2) \leq \frac{1}{2}b(\tilde{V}_r,\Z_2)
\end{equation}
(using for example \cite[Proposition 7.27]{BPRbook2}). 

Combining 
\eqref{eqn:algebraic-total-degree-tilde-V} and \eqref{eqn:algebraic-total-degree-tilde-V-and-S} we obtain
\begin{equation*}
b(V_r,\Z_2)=b(S_r,\Z_2) \leq \frac{1}{2} b(\tilde{V}_r,\Z_2) \leq \frac{1}{2} (1+(2d-1)^k),
\end{equation*}
which proves \eqref{eqn:algebraic-total-degree-tilde-bound-2}.

Finally, using conical structure theroem at infinity for semi-algebraic sets, we have that for all $r$ large enough,
$V_r$ is semi-algebraically homeomorphic to $V$, from whence it follows that
\[
b(V,\Z_2) \leq \frac{1}{2} b(\tilde{V}_r,\Z_2) \leq \frac{1}{2} (1+(2d-1)^k).
\]

We now prove that
\begin{equation}
\label{eqn:algebraic-total-degree-tilde-bound-1}
 b(V,\Z_2) \leq \binom{\ell}{k}2^kd'^k + \sum_{j=1}^{k-1}\binom{\ell}{j}2^j (F_1(d',k,j) + F_2(d',k,j)) +3,
 \end{equation}
 where $d'$ is the least even integer $\geq 0$.
 
 Denote for every \emph{even}  $d'' \geq 0$, by  $H_{d'',k} \subset \R\la\eps\ra[X_1,\ldots,X_k]$ the subspace of
 polynomials of degree $\leq d''$, and observe that 
 there is a non-empty, open semi-algebraic subset $U_{d'',k} \subset H_{d'',k}$ 
 such that for every $H \in U_{d'',k}$, 
 $H$ is strictly positive on 
 $\R\la\eps_0\ra^k$.

Denote by $W \subset H_{2,k} \times \underbrace{H_{d'',k} \times \cdots \times H_{d'',k}}_{2\ell}$ the non-empty Zariski open subset of 
generic tuples of polynomials (cf. Definition \ref{def:generic}). Then, $W \cap  U_{2,k} \times U_{d',k} \times \cdots \times U_{d',k}$ is clearly non-empty.
Choose polynomials $H_0,H_1,\ldots,H_{2\ell}$ such that the tuple 
$(H_0,H_1,\ldots,H_{2\ell}) \in W \cap  U_{2,k} \times U_{d',k} \times \cdots \times U_{d',k}$.
 
Let  
\[
\tilde{Q}_0=(1- \delta)\cdot(\|\X\|^2-1/\eps^2) - \delta \cdot H_0,
\]
and for $1 \leq i \leq \ell$, let 
\begin{eqnarray*}
\tilde{Q}_{i,+}  &=& (1-\delta)\cdot Q_i+\delta \cdot H_{2i-1}, \\
\tilde{Q}_{i,-}   &=& (1-\delta)\cdot Q_i-\delta \cdot H_{2i}.
\end{eqnarray*}
Also, let $\R' = \R\la \eps, \delta\ra$.

Note that using Corollary \ref{cor:generic}, the polynomials
$Q_0,Q_{1,\pm},\ldots,Q_{\ell,\pm}$ are generic.

We need the following lemma.
\begin{lemma}
\label{lem:algebraic-total-degree}
The real algebraic variety $\Ext(V,\R')$ (cf. Notation \ref{not:ext}) is semi-algebraically homotopy equivalent to the semi-algebraic set $\tilde{S} \subset \R'^k$ defined by 
\[
(\tilde{Q}_0 < 0) \wedge \bigwedge_{1 \leq i \leq \ell}( (\tilde{Q}_{i,+} > 0) \wedge 
(\tilde{Q}_{i,-} < 0)).
\] 
\end{lemma}
 
 \begin{proof}
 Follows from \cite[Lemma 16.17]{BPRbook2}. 
 \end{proof}
 
Let 
%%sb \epsilon ---> \sigma
\[
\tilde{W} = \ZZ(\tilde{Q}_0,\R'^k) \cup \bigcup_{\substack{1 \leq i \leq \ell\\ \sigma\in\{+,-\}}} \ZZ(\tilde{Q}_{i,\sigma},\R'^k).
\]

Note that $\R'^k \setminus \tilde{W}$ is an open semi-algebraic set, and $\tilde{S}$ is the union of
a subset of the semi-algebraically connected components of $\R'^k \setminus \tilde{W}$ which are bounded. This implies that
\[
b(\tilde{S},\Z_2) \leq b(\overset{\circ}{\R'^k}\setminus \tilde{W},\Z_2),
\]
where $\overset{\circ}{\R'^k}$ is the one-point compactification of $\R'^k$, and is
semi-algebraically homeomorphic to the sphere $\Sphere^{k}$ defined over $\R'$. 
We will now use Alexander duality theorem for semi-algebraic subsets defined over $\R'$.
Note that the statement of this theorem over $\mathbb{R}$ restricted to semi-algebraic sets
(proved for example in \cite[page 296]{Spanier}) implies the corresponding statement over $\R'$ by a standard application of the Tarski-Seidenberg transfer principle. 
Using Alexander duality 
and the fact that $\tilde{W}$ is non-empty we obtain
 \[
 b(\tilde{S},\Z_2)  \leq b(\overset{\circ}{\R'^k} \setminus \tilde{W},\Z_2) \leq b(\tilde{W},\Z_2) + 1.
 \] 

We now bound $b(\tilde{W},\Z_2)$ using Proposition \ref{prop:unionint} noting that
$\tilde{W}$ is the union of $2\ell+1$ real algebraic sets.

Let 
%%sb epsilon --> sigma
\[
\tilde{\mathcal{Q}} = \{\tilde{Q}_{i,\sigma} \;\mid\; 1 \leq i \leq \ell, \sigma \in\{+,-\} \}.
\]

Proposition \ref{prop:unionint} implies that 
\begin{equation}
\label{eqn:algebraic-total-degree-MV1}
b(\tilde{W},\Z_2) \leq  
\sum_{\substack{\tilde{\mathcal{Q}}' \subset \tilde{\mathcal{Q}}, \\ \card(\tilde{\mathcal{Q}}') \leq k-1}}  b(\ZZ(\{\tilde{Q}_0\}\cup \tilde{\mathcal{Q}}',\R'^k),\Z_2) + \sum_{\substack{\tilde{\mathcal{Q}}' \subset \tilde{\mathcal{Q}}, \\ \card(\tilde{\mathcal{Q}}') \leq k}}
b(\ZZ(\tilde{\mathcal{Q}}',\R'^k),\Z_2).
\end{equation}

Notice that for $\tilde{\mathcal{Q}}' \subset \tilde{\mathcal{Q}}$ with $\card(\tilde{\mathcal{Q}}') = k$, we have by Bezout's theorem that
\begin{equation}
\label{eqn:algebraic-total-degree-Bezout}
b(\ZZ(\tilde{\mathcal{Q}}',\R'^k),\Z_2) \leq d'^k.
\end{equation}

Inequalities \eqref{eqn:algebraic-total-degree-Bezout} and \eqref{eqn:algebraic-total-degree-MV1}
imply
\begin{equation}
\label{eqn:algebraic-total-degree-MV2}
b(\tilde{W},\Z_2) \leq  
\binom{\ell}{k} 2^k d'^k + \sum_{\substack{\tilde{\mathcal{Q}}' \subset \tilde{\mathcal{Q}}, \\ \card(\tilde{\mathcal{Q}}') \leq k-1}}  (b(\ZZ(\{\tilde{Q}_0\}\cup \tilde{\mathcal{Q}}',\R'^k),\Z_2) + b(\ZZ(\tilde{\mathcal{Q}}',\R'^k),\Z_2)).
\end{equation}

Finally for $\tilde{\mathcal{Q}}' \subset \tilde{\mathcal{Q}}$ with $\card(\tilde{\mathcal{Q}}') = j, 1 \leq j \leq \min(\ell,k-1)$,
we have,
\begin{equation}
\label{eqn:algebraic-total-degree-MV3}
b(\ZZ(\{\tilde{Q}\}\cup \tilde{\mathcal{Q}}',\R'^k),\Z_2)  \leq  F_1(d',k,j),
\end{equation}
using inequality \eqref{eqn:degree2andd} in Remark \ref{rem:degree2andd},
and 
\begin{equation}
\label{eqn:algebraic-total-degree-MV4}
b(\ZZ(\tilde{\mathcal{Q}}',\R'^k),\Z_2)) \leq 
F_2(d',k,j)
\end{equation}
using  inequality \eqref{eqn:many-total-diff-degrees-real-2} in 
Proposition \ref{prop:many-total-diff-degrees}. Note that if $j=0$, then 
\[
b(\ZZ(\{\tilde{Q_0}\},\R'^k),\Z_2)=2,
\] 
and 
\[
b(\ZZ(\tilde{\mathcal{Q}}',\R'^k),\Z_2)=b(\ZZ(\emptyset,\R'^k),\Z_2)=b(\R'^k,\Z_2)= 1.
\]
 
Inequality \eqref{eqn:algebraic-total-degree-tilde-bound-1} now follows from
Lemma \ref{lem:algebraic-total-degree}, and inequalities \eqref{eqn:algebraic-total-degree-MV1},  
\eqref{eqn:algebraic-total-degree-MV2},  
\eqref{eqn:algebraic-total-degree-MV3}, 
and \eqref{eqn:algebraic-total-degree-MV4}.

Finally, inequality \eqref{eqn:algebraic-total-degree} follows from 
inequalities \eqref{eqn:algebraic-total-degree-tilde-bound-2}
and \eqref{eqn:algebraic-total-degree-tilde-bound-1}.
\end{proof}

\subsection{Proofs of Theorems \ref{thm:multi-algebraic} and \ref{thm:multi-semi}}
\label{subsec:proofs-multi}

\begin{proof}[Proof of Theorem \ref{thm:multi-algebraic}]
We first prove
\begin{eqnarray}
\label{eqn:multi-algebraic1}
b(V,\Z_2) &\leq& \frac{1}{2}G_{\gen}(2\dd,\kk,1).
\end{eqnarray}

The proof is similar to that of Theorem \ref{thm:algebraic-total-degree}, but we note that the polynomial (\ref{eqn:Qtilde}) $\tilde{Q}=Q_1^2+\cdots+Q_\ell^2+c (\| \X\|^2-r^2)\in \R[\X^{(1)},\ldots,\X^{(p)}]$ has multi-degree bounded by $2\dd$. Therefore, 
$b(\ZZ(\tilde{Q},\R^k),\Z_2) \leq G_{\gen}(2\dd,\kk,1)$
using Proposition \ref{prop:manyblocks-total}.

We now prove
\begin{eqnarray}
\label{eqn:multi-algebraic2}
b(V,\Z_2) &\leq&  3+\sum_{j=1}^{k} \binom{\ell}{j}2^j (G_{\gen}(\dd',\kk,j) + G_{\gen}(\dd',\kk,j+1)).
\end{eqnarray}

We proceed in the same manner as in the proof of Theorem \ref{thm:algebraic-total-degree}. We note that the sphere can also be viewed as a polynomial in  $\R[\X^{(1)},\ldots,\X^{(p)}]$, where each block has degree equal to $2$. Notice that we assume all $d_i \geq 2$, so we can view the polynomial $Q_0$ in the proof of Theorem 
\ref{thm:algebraic-total-degree}
as another polynomial with the same block structure and degree bounds as each polynomial in $\mathcal{Q}$. Therefore, we can replace both $F_1(d',k,j)$ (resp. $F_2(d',k,j)$) with $G_{\gen}(\dd',\kk,j+1)$ (resp. $G_{\gen}(\dd',\kk,j)$).
The theorem follows from inequalities \eqref{eqn:multi-algebraic1} and \eqref{eqn:multi-algebraic2}.
\end{proof}

\begin{proof}[Proof of Theorem \ref{thm:multi-semi}]
The proof is similar to those of \cite[Theorem 7.30]{BPRbook2} and \cite[Theorem 7.38]{BPRbook2}
%%ssb
%%. From Proposition \ref{prop:many-total-diff-degrees} 
with the following modification that instead of using Theorem \ref{thm:O-P-T-M}
for bounding the sum of the Betti numbers of various algebraic sets that occur, we use the bound in Theorem \ref{thm:multi-algebraic}.
\end{proof}

\subsection{Proofs of Theorems \ref{thm:different-boxes-algebraic} and \ref{thm:different-boxes-semi}}
\label{subsec:proofs-different-boxes}

\begin{proof}[Proof of Theorem \ref{thm:different-boxes-algebraic}]
The proof is similar to the proof of Theorem \ref{thm:multi-algebraic} using Proposition
\ref{prop:different-boxes} instead of Proposition \ref{prop:manyblocks-total}.
\end{proof}

\begin{proof}[Proof of Theorem \ref{thm:different-boxes-semi}]
The proof is similar to the proof of Theorem \ref{thm:multi-semi} using 
%%sb
%%Proposition \ref{prop:different-boxes} 
Theorem \ref{thm:different-boxes-algebraic}
instead of 
%%sb
%%Proposition \ref{prop:manyblocks-total}.
Theorem \ref{thm:multi-algebraic}.
\end{proof}

\subsection{Proofs of Theorems \ref{thm:partly-quadratic-total-algebraic}, \ref{thm:Lerario-new}, \ref{thm:partly-quadratic-total-semi}, and \ref{thm:BP'R-new}}

\begin{proof}[Proof of Theorem \ref{thm:partly-quadratic-total-algebraic}]
The proof is similar to that of Theorem \ref{thm:multi-algebraic}. 
Since we have $\ell$ partially quadratic polynomials, we use $H_{\gen}(d',k_1,k_2,j)$ in place of $G_{\gen}(\dd',\kk,j)$, where $H_{\gen}(d',k_1,k_2,j)$ is the bound from Proposition \ref{prop:many-total-mixed}, noting that we assume $d\geq 2$.
\end{proof}

\begin{proof}[Proof of Theorem \ref{thm:Lerario-new}]
The proof is similar to the proof of Theorem \ref{thm:algebraic-total-degree} above using generic \emph{positive} quadrics to perturb the given polynomials (in lieu of the polynomials $H_i$ in the proof of
Theorem \ref{thm:algebraic-total-degree}).
More precisely,
for $1 \leq i \leq \ell$, let 
\begin{eqnarray*}
\tilde{Q}_{i,+}  &=& (1 -\delta)\cdot Q_i+\delta \cdot H_{2i-1}, \\
\tilde{Q}_{i,-}   &=& (1 - \delta)\cdot Q_i-\delta \cdot  H_{2i}.
\end{eqnarray*}
where the polynomials
$H_1,\ldots,H_{2\ell}$ are chosen to be
generic positive quadrics.
Using Corollary \ref{cor:generic}, the polynomials
$Q_0,Q_{1,\pm},\ldots,Q_{\ell,\pm}$ are generic.

Also, let $\R' = \R\la\delta\ra$.

Using Lemma
\ref{lem:algebraic-total-degree}
The real algebraic variety $\Ext(V,\R')$ (cf. Notation \ref{not:ext}) is semi-algebraically homotopy equivalent to the semi-algebraic set $\tilde{S} \subset \PP_{\R'}^k$ defined by 
\[
\bigwedge_{1 \leq i \leq \ell}( (\tilde{Q}_{i,+} > 0) \wedge 
(\tilde{Q}_{i,-} < 0)).
\]

Let
%%sb epsilon ---> sigma 
\[
\tilde{W} =  \bigcup_{\substack{1 \leq i \leq \ell\\ \sigma \in\{+,-\}}} \ZZ(\tilde{Q}_{i,\sigma},\R'^k).
\]

Note that $\PP_{\R'}^k \setminus \tilde{W}$ is an open semi-algebraic set, and $\tilde{S}$ is the union of
a subset of the semi-algebraically connected components of $\PP_{\R'}^k \setminus \tilde{W}$. This implies that
\[
b(\tilde{S},\Z_2) \leq b(\PP_{\R'}^k\setminus \tilde{W},\Z_2).
\]

We will now use Lefschetz duality theorem for semi-algebraic subsets defined over $\R'$.

Using Lefschetz duality we obtain that
 \begin{eqnarray}
 \label{eqn:Lerario-new1}
 \nonumber
 b(\tilde{S},\Z_2)  &\leq&  b(\PP_{\R'}^k \setminus \tilde{W},\Z_2) \\ 
 \nonumber
  &\leq&  b(\tilde{W},\Z_2) + b(\PP_{\R'}^k,\Z_2) \\ 
  &=& b(\tilde{W},\Z_2) + k+1.
 \end{eqnarray}

We now bound $b(\tilde{W},\Z_2)$ using Proposition \ref{prop:unionint} noting that
$\tilde{W}$ is the union of $2\ell$ real algebraic sets.

Let 
%%sb epsilon --> sigma
\[
\tilde{\mathcal{Q}} = \{\tilde{Q}_{i,\sigma} \;\mid\; 1 \leq i \leq \ell, \sigma \in\{+,-\} \}.
\]

Proposition \ref{prop:unionint} implies that 
\begin{eqnarray}
 \label{eqn:Lerario-new2}
b(\tilde{W},\Z_2) &\leq&  
\sum_{\tilde{\mathcal{Q}}' \subset \tilde{\mathcal{Q}}}
b(\ZZ(\tilde{\mathcal{Q}}',\PP_{\R'}^k),\Z_2).
\end{eqnarray}

Now use inequalities \eqref{eqn:Lerario-new1},  \eqref{eqn:Lerario-new2}, and 
\eqref{eqn:quadratic-projective-betti-exact}  to obtain the bound in the theorem.
\end{proof}

\begin{proof}[Proof of Theorem \ref{thm:partly-quadratic-total-semi}]
The proof is similar to those of \cite[Theorem 7.30]{BPRbook2} and \cite[Theorem 7.38]{BPRbook2} 
 with the modification that instead of using Theorem \ref{thm:O-P-T-M}
for bounding the sum of the Betti numbers of various algebraic sets that occur, we use the bound in Theorem \ref{thm:partly-quadratic-total-algebraic}.
\end{proof}

\begin{proof}[Proof of Theorem \ref{thm:BP'R-new}]
The proof is again similar to the proofs of Theorem 7.30 and Theorem 7.38 in \cite{BPRbook2} with several modifications. 
In the proof of Theorem 7.30 in \cite{BPRbook2}, we let $\mathcal{Q} = \{0\}$, and 
$\mathcal{P} = \mathcal{P}_1 \cup \mathcal{P}_2$, and
we let $\mathcal{P}_1 = \{P_1,\ldots,P_s\}$ and $\mathcal{P}_2 =\{P_{s+1},\ldots,P_{s+m}\}$.
In Proposition 7.34, for
each sign condition $\sigma \in \{0,1,-1\}^{\mathcal{P}}$, we redefine the basic closed semi-algebraic set
\[
\overline{\RR(\sigma)} \subset \R\la\delta,\delta_1,\ldots,\delta_{s+m},\eps_1,\ldots,\eps_{s+m}\ra^k
\] 
in the following way. 
Without loss of generality assume

$$
\displaylines{
\sigma(P_h)=0  \mbox{ if  } h \in I, \cr
\sigma(P_h)=1 \mbox{ if }  h \in J, \cr
\sigma(P_h)=-1  \mbox{ if } h \in \{1,\ldots,s+m\}\setminus(I \cup J),
}
$$

and denote by $\overline{\RR(\sigma)}$ the subset of $\R\la\delta,\delta_1,\ldots,\delta_{s+m},\eps_1,\ldots,\eps_{s+m}\ra^k$
defined by
$$
\displaylines{
\delta(|\X|^2 + |\Y|^2) \leq 1, \cr
-\eps_h \leq P_h \leq \eps_h, \mbox{ if } h \in I, \cr
P_{h}  \geq \delta_h, \mbox{ if } h \in J, \cr
P_{h} \leq -\delta_h, \mbox{ if } h \in \{1,\ldots,s+m\}\setminus(I \cup J).
}
$$

It is easy to verify that Proposition 7.34 in \cite{BPRbook2} remains true with this new definition of $\overline{\RR(\sigma)}$.

Now we observe that since $\mathcal{P}_1 \subset \R[X_1,\ldots, X_{k_1}]$, no more than $k_1$ polynomials
amongst the set $\{P_1 \pm \delta_1, P_1 \pm \eps_1,\ldots,P_s \pm \delta_s, P_s \pm \eps_s\}$ can have a common zero. 
Each non-empty real algebraic set $V$ defined by some subset of the polynomials  
$\{P_1 \pm \delta_1, P_1 \pm \eps_1,\ldots,P_s \pm \delta_s, P_s \pm \eps_s\} \cup \{P_0\}$,
where $P_0 = \delta(|\X|^2+ |\Y|^2) -1$,
is the set of zeros of  two sets of polynomials, namely
%%sb changes - proof-corections
\hide{
$$
\displaylines{
(P_h  + \epsilon_{h} \eta_{h})_{h\in J_1},  \eta_h \in \{\eps_h,\delta_h\}, \epsilon_h \in \{\pm 1, \pm 2\}, J_1 \subset [1,s], \cr
(P_h  + \epsilon_{h} \eta_{h})_{h\in J_2},  \eta_h \in \{\eps_h,\delta_h\}, \epsilon_h \in \{\pm 1, \pm 2\}, J_2 \subset [s+1,s+m],
}
$$
}
%%sb end of hide
$$
\displaylines{
(P_h  + \sigma_{h} \eta_{h})_{h\in J_1},  \eta_h \in \{\eps_h,\delta_h\}, \sigma_h \in \{\pm 1\}, J_1 \subset [1,s], \cr
(P_h  + \sigma_{h} \eta_{h})_{h\in J_2},  \eta_h \in \{\eps_h,\delta_h\}, \sigma_h \in \{\pm 1\}, J_2 \subset [s+1,s+m],
}
$$

and  possibly $P_0$,
with  $j_1= \card(J_1) \leq k_1$,
and $j_2 = \card(J_2) \leq \min(m+1, k_1+k_2 - j_1-i)$.

We also note that $V$ is defined by the  $(\card(J_2) + 1)$ or $(\card(J_2) +2)$ (depending on whether 
$P_0$ is included or not) polynomials
%%sb changes -- proof-corrections
\hide{
\[
\sum_{h \in J_1} (P_h  + \epsilon_{h} \eta_{h})^2, (P_h  + \epsilon_{h} \eta_{h})_{h\in J_2}
\]
}
%%sb end of hide
\[
\sum_{h \in J_1} (P_h  + \sigma_{h} \eta_{h})^2, (P_h  + \sigma_{h} \eta_{h})_{h\in J_2}
\]
(and possibly $P_0$).
The degrees of these polynomials are at most $2d$ in $\X$, and at most $2$ in $\Y$.
 We can use Theorem \ref{thm:partly-quadratic-total-algebraic} to bound  $b(V,\Z_2)$ by $H(2d,k_1,k_2,\card(J_2)+1)$
 or $H(2d,k_1,k_2,\card(J_2)+2)$ (depending on whether 
$P_0$ is included or not) 
(cf. Eqn. \eqref{eqn:partly-quadratic-total-algebraic}).

Moreover, the total
number of non-empty real algebraic sets $V$ that occur in the proof is bounded by 
\[
\sum_{\substack{0\leq j_1\leq k_1\\ 0 \leq j_2 \leq \min(m+1,k_1+k_2-j_1 - i)}} \binom{s}{j_1}\binom{m+1}{j_2} 4^{j_1+j_2}.
\]
One now obtains inequality \eqref{eqn:BP'R-new1} by following the rest of the argument in the proof of 
Theorem 7.30. Note that we needed to increase the number of polynomials by one by including the polynomial
$P_0$. This accounts for the $m+1$ in the subscript of the second sum in the bound.  

The proof of inequality \eqref{eqn:BP'R-new2} is by a similar modification of the proof of Theorem 7.38 in \cite{BPRbook2}  and is omitted.
\end{proof}

\subsection{Proofs of Theorems \ref{thm:partly-quadratic-multi-algebraic} and  \ref{thm:partly-quadratic-multi-semi}}
\begin{proof}[Proof of Theorem \ref{thm:partly-quadratic-multi-algebraic}]
The proof is similar to that of Theorem \ref{thm:multi-algebraic}. 
Since we have $\ell$ partially quadratic polynomials with several blocks, we use the bound 
\[M_{\gen}(\dd',k_1,k_2,j)
\] 
in place of 
$G_{\gen}(\dd',\kk, j)$, where $M_{\gen}(\dd',k_1,k_2,j)$ is the bound from Proposition \ref{prop:several-blocks-mixed}, noting that we assume each $d_i\geq 2$.
\end{proof}

\begin{proof}[Proof of Theorem \ref{thm:partly-quadratic-multi-semi}]
The proof is similar to that of Theorem \ref{thm:multi-semi}, except that the bound for the sum of the Betti numbers corresponding to a subset of indices of cardinality  $j$ is now given by $M(\dd',k_1,k_2,j)$.
\end{proof}

\section{A few applications}
\label{sec:applications}
In this section we give a few applications of the results proved in the last section.
\subsection{Bounding Betti numbers of pull-backs and direct images under polynomial maps}
\label{subsec:applications}
We discuss a few immediate applications of the multi-degree bounds proved in
\S \ref{sec:main-results}.

\begin{theorem}[Bound on pull-back]
\label{thm:pull-back}
Let $\mathcal{F} = \{F_1,\ldots,F_m\} \subset \R[X_1,\ldots,X_k]$ and  $\mathcal{G}\subset
\R[Y_1,\ldots,Y_m]$, with $\deg(F) \leq d, F \in \mathcal{F}$, and
$\deg(G) \leq D, G \in \mathcal{G}$, and let $\card(\mathcal{G}) = s$. Let $\mathbf{F}: \R^k \rightarrow \R^m$ denote the polynomial
map $x \mapsto (F_1(x),\ldots,F_m(x))$, and let $S \subset \R^m$ be a $\mathcal{G}$-closed semi-algebraic set.
Then,
\begin{eqnarray*}
b(\mathbf{F}^{-1}(S),\Z_2) &\leq& \sum_{i=0}^{k+m} \sum_{j=1}^{k+m-i} \binom{m+s+1}{j} 6^{j} 
G_{\min}(\dd, \kk,j) \\
&\leq&  O(1)^{k+m}
(m+s)^{k+m}
d^kD^m.
\end{eqnarray*}
\end{theorem}

 \begin{proof}
 Suppose that $\Phi(Y_1,\ldots,Y_m)$ is a $\mathcal{G}$-closed formula defining $S$.
 Notice that $\mathbf{F}^{-1}(S)$ is semi-algebraically homeomorphic to the semi-algebraic subset of
 $\R^{k+m}$ defined by the formula
 \[
 \Psi(\X,\Y) := \bigwedge_{i=1}^{m} (Y_i - F_i = 0) \wedge \Phi(Y_1,\ldots,Y_m).
 \]
 
 The number of polynomials appearing in $\Psi$ is bounded by $m + s$. The degrees in $\Y$ of the polynomials appearing in $\Psi$ are bounded by $D$, while the degrees in $\X$ are bounded by $d$.
Applying Theorem \ref{thm:multi-semi}  with $p=2$, $\kk = (k,m)$, and $\dd= (d,D)$,  we obtain
\begin{eqnarray*}
b(\RR(\Psi,\R^{k+m}), \Z_2) &\leq& \sum_{i=0}^{k+m} \sum_{j=1}^{k+m-i} \binom{m+s+1}{j} 6^{j} 
G_{\min}(\dd, \kk,j) \\
&\leq& \sum_{i=0}^{k+m} \sum_{j=1}^{k+m-i} \binom{m+s+1}{j} 6^{j} O(1)^{k+m}
d^kD^m  \mbox{ (using (\ref{eqn:Gmin}))} \\
&\leq&  O(1)^{k+m}(m+s)^{k+m}d^kD^m.
\end{eqnarray*}
\end{proof} 

\begin{theorem}[Bound on image]
\label{thm:push-forward}
Let $\mathcal{F} = \{F_1,\ldots,F_m\}, \mathcal{G} \subset \R[X_1,\ldots,X_k]$, with $\deg(F) \leq d, F \in \mathcal{F}$, and $\deg(G) \leq D, G \in \mathcal{G}$, and let $\card(\mathcal{G}) = s$. Let $\mathbf{F}: \R^k \rightarrow \R^m$ denote the polynomial
map $x \mapsto (F_1(x),\ldots,F_m(x))$, and let $T \subset \R^k$ be a bounded $\mathcal{G}$-closed semi-algebraic set.
Suppose also that $d \geq D$.

Then, for $0 \leq i \leq m$, 
\begin{eqnarray*}
b_i(\mathbf{F}(T),\Z_2)  &\leq& 
\sum_{j=0}^{i} 
\sum_{h=0}^{\alpha_j}
\sum_{\ell=1}^{\alpha_j-h} \binom{(j+1)(m+s)+1}{\ell}6^\ell G_{\min}(\dd,\kk,\ell) \\
 &\leq&
 O(i)^{\alpha_i} (m+s)^{\alpha_i} d^{(i+1)k}D^m
 \end{eqnarray*}
 where $\alpha_i=(i+1)k+m$.
\end{theorem}

\begin{proof}
Using the descent spectral sequence we have that 
\begin{eqnarray}
\label{eqn:push-forward1}
b_i(\mathbf{F}(T),\Z_2) &\leq& \sum_{j=0}^{i} b_{i-j}(\underbrace{T \times_{\mathbf{F}} \cdots \times_{\mathbf{F}} T}_{(j+1)}, \Z_2) \nonumber \\
&\leq& \sum_{j=0}^{i} b(\underbrace{T \times_{\mathbf{F}} \cdots \times_{\mathbf{F}} T}_{(j+1)}, \Z_2).
\end{eqnarray}
Suppose that $T$ is defined by a $\mathcal{G}$-closed formula $\Psi$. Notice that for
all $j \geq 0$, $\underbrace{T \times_{\mathbf{F}} \cdots \times_{\mathbf{F}} T}_{(j+1)}$ is 
defined by the formula 
\[
\Psi^{(j)}(\X^{(0)},\ldots,\X^{(j)},\Y) := \bigwedge_{i=1}^{m} \bigwedge_{h=0}^{j} \Psi(\X^{(h)},\Y)\wedge(Y_i - F_i(\X^{(h)}) = 0),
\]
where $\Y= (Y_1,\ldots,Y_m), \X^{(h)} = (X^{(h)}_1,\ldots,X^{(h)}_k), 0 \leq h \leq j$.

The cardinality of the set of polynomials appearing in $\Psi^{(j)}$ is $(j+1)(m+s)$, the degree in each
block $\X^{(h)}$ is bounded by $d$, and that in $\Y$ is bounded by $D$.

Denote by $\alpha_j=(j+1)k+m$. Now apply Theorem \ref{thm:multi-semi} with $p = j+2$, $\kk = (\underbrace{k,\ldots,k}_{j+1},m)$,
$\dd = (\underbrace{d,\ldots,d}_{j+1},D)$ to obtain
\begin{eqnarray}
\label{eqn:push-forward2}
b(\underbrace{T \times_{\mathbf{F}} \cdots \times_{\mathbf{F}} T}_{(j+1)}, \Z_2)
&\leq& 
\sum_{h=0}^{\alpha_j}
\sum_{\ell=1}^{\alpha_j-h} \binom{(j+1)(m+s)+1}{\ell}6^\ell G_{\min}(\dd,\kk,\ell) \nonumber \\
&\leq&
O(j)^{\alpha_j} (m+s)^{\alpha_j} d^{(j+1)k}D^m
\end{eqnarray}

The theorem now follows from Eqns. \eqref{eqn:push-forward1} and \eqref{eqn:push-forward2}.
\end{proof}

\begin{remark}
\label{rem:push-forward-and-pull-back}
Note that versions of Theorem \ref{thm:pull-back} and Theorem \ref{thm:push-forward} without the
distinction between the two degrees $d$ and $D$ were known before (see \cite{GaV}). The novel
aspect of these two theorems is the different dependence of the bounds proved on the two degrees 
$d,D$. In certain applications, this distinction is important. 
\end{remark}

We record the following result similar to that of Theorem  \ref{thm:push-forward}, which is also useful
in practice.
The following situation occurs very frequently in semi-algebraic geometry.

 Let $S_1 \subset \R^{k}$, and $S_2 \subset \R^{k} \times \R^{m}$ be 
  semi-algebraic subsets, and 
 $\pi_\X:\R^{k} \times \R^{m} \rightarrow \R^{k},\pi_\Y: \R^{k} \times \R^{m} \rightarrow \R^m$ be the two projection maps on the first and second factors resp.. 
 Let $T = \pi_\Y(\pi_\X^{-1}(S_{1}) \cap S_2)$ (see figure below).
\[
\xymatrix{
& \pi_\X^{-1}(S_{1}) \cap S_2 \ar[ld]_{\pi_\X} \ar[rd]^{\pi_\Y} \\
S_{1}  && T.
}
\]

\begin{theorem}[Set-theoretic Fourier-Mukai transform]
\label{thm:Fourier-Mukai}
With the same notation as above, 
let $\mathcal{P}_{1} \subset \R[X_1,\ldots,X_k]$ and $\mathcal{P}_2 \subset \R[X_1,\ldots,X_k,Y_1,\ldots,Y_m]$
be finite sets such  that  $S_{1} \subset \R^k$ is a $\mathcal{P}_{1}$-closed semi-algebraic set and 
$S_2 \subset \R^{k+m}$ is a bounded $\mathcal{P}_2$-closed semi-algebraic set.
Suppose that  $\deg_\X(\mathcal{P}_{1}), \deg_\X(\mathcal{P}_2) \leq d$ and $\deg_\Y(\mathcal{P}_2)\leq D$.
Suppose also that $d \geq D$.

Then for $0 \leq i \leq m$,
\begin{eqnarray*}
b_i(T,\Z_2) &\leq& \sum_{j=0}^{i} \sum_{h=0}^{\alpha_j} \sum_{\ell=1}^{\alpha_j-h} \binom{(j+1)(s_1+s_2)+1}{\ell}6^\ell G_{\min}(\dd,\kk,\ell) \\
&\leq& O(i)^{\alpha_i} (s_1+s_2)^{\alpha_i} d^{(i+1)k}D^{m},
 \end{eqnarray*}
 where $\alpha_i=(i+1)(k+m)+k$, and $s_1 = \card(\mathcal{P}_{1}) , s_2 = \card(\mathcal{P}_2)$.
\end{theorem}

\begin{proof}
Note that the semi-algebraic set $\pi_\X^{-1}(S_{1}) \cap S_2$ is a bounded $(\mathcal{P}_{1} \cup \mathcal{P}_2)$-closed semi-algebraic set, with degrees in $\X$ bounded by $d$ and in $\Y$ bounded by $D$. Note that if 
$S_{1}$ is defined by the formula $\Phi(X_1,\ldots,X_k)$ and $S_2$ is defined by the formula $\Psi(X_1,\ldots,X_k,Y_1,\ldots,Y_m)$, then the set $\pi_\X^{-1}(S_1)\cap S_2$ is defined by $\Phi(X_1,\ldots,X_k)\wedge \Psi(X_1,\ldots,X_k,Y_1,\ldots,Y_m)$. 
Note that with the above notation, for all $j\geq 0$, $\underbrace{T \times_{\pi_\Y} \cdots \times_{\pi_\Y} T}_{(j+1)}$ is defined by the formula
\[
\Theta^{(j)}(\X^{(0)},\ldots,\X^{(j)},\Y) := \bigwedge_{h=0}^{j} \Phi(\X^{(h)},\Y)\wedge \Psi(\X^{(h)},\Y),
\]
where $\Y= (Y_1,\ldots,Y_m), \X^{(h)} = (X^{(h)}_1,\ldots,X^{(h)}_k), 0 \leq h \leq j$.

The cardinality of the set of polynomials appearing in $\Theta^{(j)}$ is $(j+1)(s_1+s_2)$, the degree in each
block $\X^{(h)}$ is bounded by $d$, and that in $\Y$ is bounded by $D$.

Applying Theorem \ref{thm:push-forward} with $\alpha_j=(j+1)k+m$ we get
\begin{eqnarray}
b(\underbrace{T \times_{\pi_\Y} \cdots \times_{\pi_\Y} T}_{(j+1)}, \Z_2)
&\leq& \sum_{h=0}^{\alpha_j}
\sum_{\ell=1}^{\alpha_j-h} \binom{(j+1)(s_1+s_2)+1}{\ell}6^\ell G_{\min}(\dd,\kk,\ell) \nonumber \\
&\leq&
O(j)^{\alpha_j} (s_1+s_2)^{\alpha_j} d^{(j+1)k}D^{m}.
\end{eqnarray}
The theorem follows from the inequality \eqref{eqn:push-forward1}.
\end{proof}

\subsection{An application to discrete geometry}
\label{subsec:transversal}
The theory of transversals is a very well-studied topic in discrete geometry with many applications.
Suppose that $S \subset \R^k$ is a closed and bounded semi-algebraic set.
We define the space $\transversal_{k'}(S) \subset \AGr_{k,k'}(\R)$ to be the set of $k'$-dimensional
affine subspaces $\ell$ of $\R^k$ such that $\ell \cap S \neq \emptyset$ (where we denote by 
$\AGr_{k,k'}(\R)$ the space (the affine Grassmannian) 
of $k'$-dimensional affine subspaces of $\R^k$ . Upper bounds on the topology
of such spaces of transversals are important in discrete geometry (see for example \cite{GPW96}).

We prove the following theorem which improves the bound that one obtains using previously known methods
by exploiting the multi-degree bounds proved in the current paper (see Remark \ref{rem:transversals}).

\begin{theorem}
\label{thm:transversals}
Let $S \subset \R^k$ be a bounded $\mathcal{P}$-closed semi-algebraic set, where 
$\mathcal{P} \subset \R[X_1,\ldots,X_k]$ with $\deg(P) \leq d, P \in \mathcal{P}$, and 
$\card(\mathcal{P}) = s$. Then, for all $k', 0 \leq k'\leq k$,
$b_i(\transversal_{k'}(S),\Z_2) $ is bounded by
\begin{eqnarray*}
&& \sum_{j=0}^{i} \sum_{h=0}^{\alpha_j} \sum_{\ell=1}^{\alpha_j-h} \binom{(j+1)(s+m+ 2(k+1))+1}{\ell}6^\ell G_{\min}(\dd,\kk,\ell) \\
&\leq& O(i)^{\alpha_i} (s+m+2(k+1))^{\alpha_i} d^{(i+1)k},
\end{eqnarray*}
where $\kk=(\underbrace{k,\ldots,k}_{j+1},m)$, $\dd=(\underbrace{d,\ldots,d}_{j+1},2)$, and $\alpha_i=(i+1)k+m$, with $m=(k+1)(k+2)/2-1$.
\end{theorem}

\begin{proof}
We first identify $\AGr_{k,k'}(\R)$ with 
%%sb adds
an open dense semi-algebraic  subset of the 
real Grassmannian $\Gr_{k+1,k'+1}(\R)$ of $(k'+1)$-dimensional
subspaces of $\R^{k+1}$ in the standard way, identifying $\ell \in \AGr(k,k')$ with the linear hull of
$\ell' = \{(x,1) \mid x\in \ell' \subset \R^k\} \subset \R^{k+1}$. Similarly, let
$S_1 = \{(x,1) \mid x \in S\} \subset \R^k$. The set $\transversal_{k'}(S)$ can then be identified with
the space (which we also denote by $\transversal_{k'}(S)$) 
\[
\{ \ell' \in \Gr_{k+1,k'+1}(\R) \mid \ell \cap S_1 \neq \emptyset\}.
\]
Now  $\Gr_{k+1,k'+1}(\R)$ is semi-algebraically homeomorphic to the real affine variety defined by 
\begin{eqnarray}
\label{eqn:grassmannian}
\{ A \in \R^{(k+1) \times (k+1)} \mid A^t = A, A^2 = A, \trace(A) = k'+1\}.
\end{eqnarray}
(see for example \cite[Theorem 3.4.4]{BCR}).
 
We identify $\Gr_{k+1,k'+1}(\R)\subset \R^{(k+1)(k+2)/2-1}$ with the subset of the linear
subspace of the space of $(k+1) \times (k+1)$ symmetric matrices with entries in $\R$ having trace
$k'+1$ (notice that the subspace containing $\Gr_{k+1,k'+1}(\R)$  has dimension $(k+1)(k+2)/2 -1$ 
and that the degrees of the polynomials in $(k+1)(k+2)/2 -1$ variables  defining $\Gr_{k+1,k'+1}$ are all bounded by $2$).

Let $S_2 \subset \R^{k} \times \R^{(k+1)(k+1)/2-1}$ be the semi-algebraic set 
(the total space of the tautological bundle over $\Gr_{k+1,k'+1}(\R)$)
defined by 
\begin{eqnarray}
\label{eqn:transversals-S2}
S_2 = \{(x,A) \mid x \in \R^k, A \in \Gr_{k+1,k'+1}(\R),  A x' = x', x' = (x,1) \}.
\end{eqnarray}

Let $\pi_1, \pi_2$ be the
projection maps as depicted in the following figure.
\[
\xymatrix{
 & \R^{k+1} \times \Gr_{k+1,k'+1}(\R) \ar[ld]_{\pi_1} \ar[rd]^{\pi_2} & \\
 \R^{k+1} & & \Gr_{k+1,k'+1}(\R).
 }
 \]
 
Observe that
\[
\transversal_{k'}(S) =  \pi_2(\pi_1^{-1}(S_1) \cap S_2).
\]
Now apply Theorem \ref{thm:Fourier-Mukai} 
noting that the number of polynomial equations (each of degree $2$) used to define
$\Gr_{k+1,k'+1}(\R)$ in Eqn. \eqref{eqn:grassmannian} is equal to $k+m+1$, where
$m = (k+1)(k+2)/2-1$, and hence the number of equations (each of degree at most $2$) used in the definition
of $S_2$ in Eqn. \eqref{eqn:transversals-S2} is equal to
$m+ 2(k+1)$.
\end{proof}
 
\begin{remark}
\label{rem:transversals}
Note that if we used the more standard Pl\"{u}cker 
embedding of the Grassmannian $\Gr_{k+1,k'+1}(\R)$ in the
projective space $\PP(\bigwedge^{k'+1}\R^{k+1})$ of dimension $\binom{k+1}{k'+1}-1$, we would obtain a
bound which is doubly exponential in $k$ in the worst case. The fact that over a real closed field,
the Grassmannians are semi-algebraically homeomorphic to the real affine variety described in 
Eqn. \eqref{eqn:grassmannian} allows us to obtain a much better bound (which is only singly exponential
in $k$). Secondly, if we used the best known prior results on effective quantifier elimination to estimate
$b_i(\transversal_{k'}(S),\Z_2)$ from above, we would obtain a bound of $(O(k s d))^{km}$ which has a much
worse dependence on $d$ than the bound proved in Theorem \ref{thm:transversals}.
\end{remark}
 
\section{Bound on the Betti numbers of real varieties defined by two polynomials having different 
 degrees}
 \label{sec:refined-betti}
 \subsection{Background}
 It was mentioned in the introduction that quantitative bounds on the Betti numbers (in particular, on the
 $0$-th Betti number) has proved to be important tools in several areas. More recently, triggered by
 the development of a new technique in discrete geometry (namely, 
 the \emph{ polynomial partitioning} method) it 
 became necessary to prove bounds which has a finer dependence on the degree \emph{sequence} of the
 polynomials rather than on the maximum degree (as in Theorem \ref{thm:B99}). The following theorem
(conjectured by J.  Matou{\v{s}}ek \cite{Matousek_private}) 
was proved in \cite{Barone-Basu11a} to meet the needs of discrete geometry
and has already found several applications.
 
\begin{theorem} \cite{Barone-Basu11a}
\label{thm:B-B}
Let 
$\mathcal{Q},\mathcal{P} \subset \R[X_1,\dots,X_k]$ be
finite subsets of non-zero polynomials
such that $\deg(Q)\leq d_1$ for all $Q\in \mathcal{Q}$,
$\deg P =d_2$ for all $P\in \mathcal{P}$,
and suppose that  $d_1 \leq d_2$.
Suppose that the real dimension of $\ZZ(\mathcal{Q},\R^k)$ is
$k' \leq k$, and that $\card(\mathcal{P}) = s$.

Then,
$$
\displaylines{
\sum_{\sigma \in \{0,1,-1\}^{\mathcal{P}}}
b_0(\RR(\sigma,\ZZ(\mathcal{Q},\R^k)),\Z_2)
}
$$
is at most
$$
\displaylines{
\sum_{j=0}^{k'}4^j{s +1\choose j}\left(
\textstyle\binom{k+1}{k-k'+j+1} \;(2d_1)^{k-k'}d^j\;\max\{2d_1,d_2 \}^{k'-j}
+2(k-j+1)\right).
}
$$
In particular,
\begin{eqnarray}
\label{eqn:B-B}
\sum_{\sigma \in \{0,1,-1\}^{\mathcal{P}}}
b_0(\RR(\sigma,\ZZ(\mathcal{Q},\R^k)),\Z_2) &\leq& O(1)^k (sd_2)^{k'} d_1^{k-k'}. 
\end{eqnarray}
\end{theorem}

Theorem \ref{thm:B-B} has proved to be important in incidence questions in discrete geometry
\cite{Solymosi-Tao,Zahl2013, SSZ2015, Basu-Sombra}. Even though in these applications it is usually a bound on the number of semi-algebraically connected components of semi-algebraic sets defined by polynomials of possibly  different degrees
that is important, it is a very interesting mathematical question (asked already in
\cite{Barone-Basu11a})  if the inequality 
\eqref{eqn:B-B}  in Theorem \ref{thm:B-B} can be
extended to  a bound on the higher Betti numbers.
We formulate below a more precise conjecture.

\begin{conjecture}
\label{conj:B-B}
With the same notation and hypothesis as in Theorem \ref{thm:B-B}, for all $i, 0 \leq i \leq k'$,
 \begin{eqnarray}
\label{eqn:B-B-conj}
\sum_{\sigma \in \{0,1,-1\}^{\mathcal{P}}}
b_i(\RR(\sigma,\ZZ(\mathcal{Q},\R^k)),\Z_2) &\leq& O(1)^k  s^{k'-i} d_1^{k-k'} d_2^{k'}. 
\end{eqnarray}
\end{conjecture}

At present we do not know how to prove Conjecture \ref{conj:B-B} except in the case $i=0$,  which is
Theorem \ref{thm:B-B},  and
the techniques used in proving Theorem \ref{thm:B-B} do not easily extend to the case of $i >0$.
In this paper,  we make some progress on this problem by proving Conjecture \ref{conj:B-B} for all 
$i \geq 0$, but only in the special case when $k'= k-1$.
In fact we prove the following slightly stronger theorem.

Unlike in the previous sections the bounds stated in this section will be valid for Betti numbers with coefficients in an arbitrary
field $\F$ rather than just $\Z_2$. This is because we do not use Smith inequality in our proofs.

\begin{theorem}
\label{thm:B-B-new}
With the same notation and hypothesis as in Theorem \ref{thm:B-B}, for all $i, 0 \leq i \leq k'<k$, 
and any field of coefficients $\F$,
\[\sum_{\sigma \in \{0,1,-1\}^{\mathcal{P}}}
b_i(\RR(\sigma,\ZZ(\mathcal{Q},\R^k)),\F)
\] 
is bounded by 
 \begin{equation}
\label{eqn:B-B-new}
\sum_{j=1}^{k'-i} \binom{s}{j}  4^j (F(2d_1,2 d_2,k) + F(2d_1,2d_2,k-1) +1) 
\leq O(1)^k  s^{k'-i} d_1 d_2^{k-1},
\end{equation}
where 
\[
F(d_1,d_2,k) =
\binom{k+1}{2} d_1 
\left( (d_1-1)^{k-1}+ \frac{4(k-1)}{3}d_2(d_2-1)^{k-2} \right).
\] 
\end{theorem}

The rest of this section is devoted to the proof Theorem \ref{thm:B-B-new}.
We begin as usual with the algebraic case.

\subsection{The algebraic case}
In this section we prove a nearly optimal bound on the sum of the Betti numbers of a real variety 
$V \subset \R^k$ defined by two
polynomials of possibly differing degrees $d_1 \leq d_2$. We prove that 
\[
b(V,\F) \leq O(1)^k d_1 d_2^{k-1}.
\]
The above bound follows from the following more precise theorem.

\begin{theorem}
\label{thm:refined-betti-main}
Let $P_1,P_2 \in \R[X_1,\ldots,X_k]$, with $0< \deg(P_1) \leq d_1, \deg(P_2) \leq d_2, 2 \leq d_1 \leq d_2$, and 
$V = \ZZ(\{P_1,P_2\},\R^k)$. Then,
\[
b(V,\F) \leq F(d_1,d_2,k) + F(d_1,d_2,k-1) +1,
\]
where 
\[
F(d_1,d_2,k) =
\binom{k+1}{2} d_1 
\left( (d_1-1)^{k-1}+ \frac{4(k-1)}{3}d_2(d_2-1)^{k-2} \right).
\] 
In particular,
\[
b(V,\F) \leq 8 \binom{k+1}{3}d_1d_2(d_2-1)^{k-2}.
\]
\end{theorem}

\begin{remark}
\label{rem:refined-betti}
Notice that direct application of Theorem \ref{thm:algebraic-total-degree} would yield a bound of 
$O(d_2)^k$ which is not optimal if $d_1 \ll d_2$.
\end{remark}

We now prove Theorem \ref{thm:refined-betti-main}. The proof involves several steps and utilizes
a few results from stratified Morse theory that we recall first.

\subsubsection{Stratified Morse Theory}
\label{subsec:SMT}
We follow the exposition in \cite{GM} (and also \cite{Basu3}.

A Whitney stratification of a space $X$ is a decomposition of $X$  into 
sub-manifolds called strata, which satisfy certain frontier conditions,
(see \cite{GM} page  37). In particular, given a compact set bounded by a 
smooth algebraic hypersurface, the boundary and the interior form a Whitney 
stratification.

Now, let $S$ be a compact Whitney stratified subset of $\re^k,$ and $f$
a restriction to $S$ of a smooth function. A {\em critical point} of $f$ is
defined to be a critical point of the restriction of $f$ to any stratum,
and a critical value of $f$ is the value of $f$ at a critical point.
A function  is called a Morse function if it has only non-degenerate
critical points when restricted to each stratum, and all its critical values
are distinct.
(There is an additional non-degeneracy condition which states
that the differential of $F$ at a critical point $p$ of a strata $S$ should not
annihilate any limit of tangent spaces to a stratum other than $S.$ However,
in our situation this will always be true.)

We now assume that $S \subset \re^k$ is a Whitney-stratified set,
and suppose that $f:S \rightarrow \R$ is a Morse function. We denote $S_x$ (resp. $S_{\leq x}$)
denote $S \cap f^{-1}(x)$ (resp. $S \cap f^{-1}((-\infty,x])$).

The first fundamental result of stratified Morse theory is the following.
\begin{theorem}\cite{GM}
\label{thm:PartA}
As $c$ varies in the open interval between two adjacent critical values,
the topological type of $S \cap \pi^{-1}((-\infty,c])$ remains constant.
\end{theorem}
Stratified Morse theory actually gives a recipe for describing the topological
change in $S_{\leq c}$ as $c$ crosses a critical value of $f$.
This is given in terms of \emph{Morse data}, which consists of a pair
of topological spaces $(A,B), A \supset B,$ with the property that as $c$
crosses the critical value $v = f(p),$ the change in  
$S_{\leq c}$  can be described by gluing in $A$ along $B$.

In stratified Morse theory the Morse data is presented as a product of
two pairs, called the \emph{tangential Morse data} and the \emph{normal Morse data}.
The notion of product of pairs is the
standard one in topology, namely 
\[ (A,B) \times (A',B')= (A\times A',A\times B'\cup B \times A') .\]

\begin{definition}[Tangential Morse data  \cite{GM}] 
\label{def:SMT-tangential-data}
The tangential Morse data at a critical point $p$ 
is then given by
$(B^{\lambda}\times B^{k-\lambda},(\partial B^{\lambda})\times B^{k-\lambda})$
where $B^k$ is the closed $k$-dimensional disk,  $\partial$ is the
boundary map, and $\lambda$ is the index of the Hessian matrix of $f$
(in any local co-ordinate system of the stratum containing $p$
in a neighborhood of $p$)
of $f$  (restricted to the stratum containing $p$)  at $p$.
\end{definition}

\begin{definition}[Normal Morse data \cite{GM}]
\label{def:SMT-normal-data}
Let $p$ 
be a critical point in some $k'$-dimensional stratum $Z$ of
a stratified subset $S$ of $R^k$.

Let $N'$ be any $(k-k')$-dimensional hyperplane passing through
the point $p$ which is transverse to $Z$ which intersects the
stratum $Z$ locally at the single point $p$.

Then, the {\em normal slice},
$N(p)$ at the point $p$ is defined to be,
\[ N(p) = N' \cap S \cap \overline{B_k(p,\delta)}, \]
for sufficiently small $\delta >0$.

Choose $\delta \gg\eps > 0$, and 
let $\ell^{-} = N(p)  \cap f^{-1}(f(p) -\eps). $
The normal Morse data has the homotopy type of the 
pair $({\rm cone}(\ell^-), \ell^-) .$
\end{definition}

The following theorem measures the change in topology as
we cross a critical value.

\begin{theorem}\cite[page 69]{GM}
\label{thm:PartB}
Let $[a,b] \subset \re$ an interval which
contains no critical values except for an isolated critical
value $v \in (a,b)$ which corresponds to a critical point
$p$ of $f$ restricted to some stratum $Z$ of $S$.  Let $\lambda$ be the
Morse index of the critical point $p$, Then, the
space $S_{\leq b}$ has the homotopy type of a space which is obtained from
$S_{\leq a}$ by attaching the pair
$(B^{\lambda}, \partial B^{\lambda}) \times ({\rm cone}(\ell^-), \ell^-) .$
\end{theorem}

We will need to use Theorem \ref{thm:PartB} in the following particularly simple situation. 
Let $S \subset \R^k$ be a closed and bounded semi-algebraic set defined by
$\bigwedge_{P \in \mathcal{P}} (P =0) \wedge Q \geq 0$, where $\mathcal{P}  \cup \{Q\}\subset \R[X_1,\ldots,X_k]$  
such that $\ZZ(P,\R^k), \ZZ(Q,\R^k), P \in \mathcal{P}$ 
are non-singular hypersurfaces  intersecting transversally. Then $S$ is Whitney stratified with two strata -- namely,
$Z= \ZZ(\mathcal{P}\cup \{Q\},\R^k)$ and $Z' = S \setminus Z$. Suppose that $f$ is a Morse function
on the stratified set $S$, and moreover
$f$ restricted to $\ZZ(\mathcal{P},\R^k)$ has no
critical points that belong to $S$. We prove the following theorem as a consequence of Theorems 
\ref{thm:PartA} and \ref{thm:PartB} above.

\begin{theorem}
\label{thm:SMT}
With tfhe assumptions stated above,
$b(S,\F)$ is bounded by the number of critical points of $f$ restricted to $\ZZ(\mathcal{P} \cup \{Q\},\R^k)$.
\end{theorem}

\begin{proof}
We note first that it suffices to prove the theorem for $\R = \re$. The general case then follows after a standard application of the
Tarski-Seidenberg transfer principle. 
Let $p \in \R^k$ be a critical point of $f$ restricted to $\ZZ(\mathcal{P} \cup \{Q\},\R^k)$ and without loss of generality let $p = \mathbf{0}$.
Let $W = T_p \ZZ(\mathcal{P},\R^k)$ and $V = T_p(\ZZ(\mathcal{P} \cup \{Q\},R^k)$, and we have $V$ is subspace of $W$ of codimension one.
Since $p$ is a non-degenerate  critical point of $f$ restricted to $\ZZ(\mathcal{P} \cup \{Q\},\R^k)$,  but not of $\ZZ(\mathcal{P},\R^k)$, the linear form
$df$ vanishes on $W$, but not on $V$. Let $\mathbf{u}$ (resp. $\mathbf{v}$)  denote the orthogonal projection of $\grad(Q)(p)$ (resp. $\grad(f)(p)$) to $V$. Note that $\mathbf{u},\mathbf{v} \neq 0$.
There are two cases to consider. We denote by $(\cdot,\cdot)$ the standard inner product in $\re^k$.
\begin{enumerate}[(a)]
\item
$(\mathbf{u},\mathbf{v}) > 0$: In this case following the notation in Definition \ref{def:SMT-normal-data}, 
$\ell^- = \emptyset$, and  it follows from Definition \ref{def:SMT-normal-data} that the normal Morse data
at $p$ equals $(p,\emptyset)$, and hence the product of the tangential and the normal Morse data  equals the
tangential Morse data in this case. Thus, in this case the change in $b(S_{\leq c})$ as $c$ crosses $f(p)$ is
$\pm 1$ as in ordinary Morse theory.
\item $(\mathbf{u},\mathbf{v}) <  0$.
In this case the normal Morse data is homotopy equivalent to the pair
$([0,1],\{0\})$. 
Since the product $(B^{\lambda}, \partial B^{\lambda}) \times ([0,1] ,\{0\})$ where $\lambda$ is the index of the 
critical point of $p$, the Morse data 
is homotopy equivalent to $(*,*)$. 
Thus in this case there is no change in the homotopy type of the sublevel
set $S_{\leq c}$ as $c$ crosses the critical value $f(p)$ (using Theorem \ref{thm:PartB})
as the pair that is being added is contractible.
\end{enumerate}
The theorem now follows from Theorems \ref{thm:PartA} and \ref{thm:PartB} just as in the case
of usual Morse theory.
\end{proof}

\subsubsection{Summary of the ideas behind the proof of Theorem \ref{thm:refined-betti-main}}
\label{subsubsec:ideas}

\begin{figure}
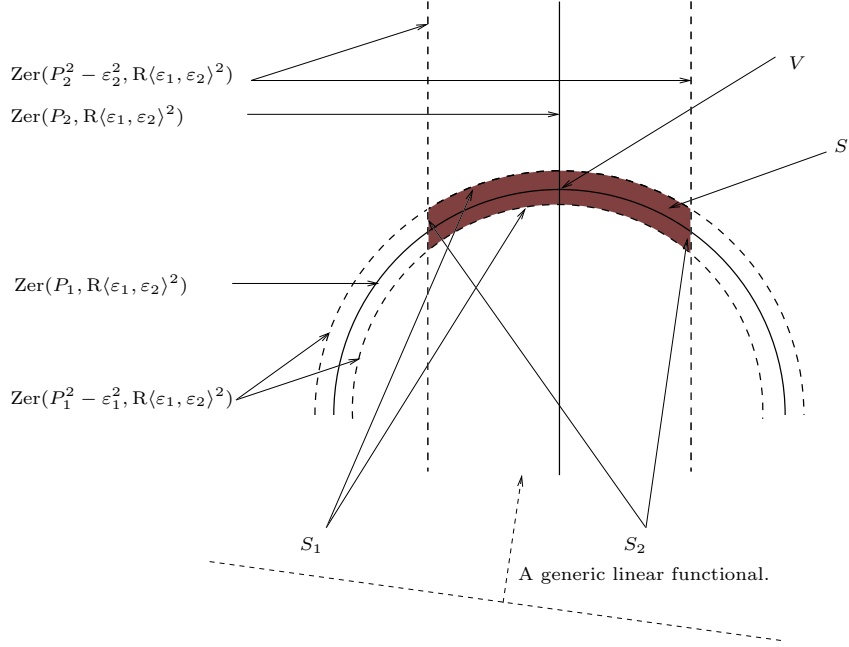
\caption{An Illustrative figure in the plane.}
\label{fig:refined-betti}
\end{figure}

For simplicity let us assume that $V = \ZZ(\{P_1,P_2\},\R^k)$ is bounded.
The case of unbounded $V$ introduces an additional complication which  we ignore in this
informal summary.
We replace $V$ by a closed bounded semi-algebraic subset $S \subset \R\la\eps_1,\eps_2\ra^k$ defined by $-\eps_i \leq P_i \leq \eps_i, i=1,2$. Then, $S$ is semi-algebraically homotopy equivalent to $V$, and moreover $S$ is a topological manifold whose boundary is a union of basic
closed semi-algebraic sets, $S_1,S_2$, where
$S_1$ is defined by $P_1^2 - \eps_1^2 = 0 \wedge -\eps_2 \leq P_2 \leq \eps_2$, and
$S_2$ is defined by $P_2^2 - \eps_2^2 = 0 \wedge -\eps_1 \leq P_1 \leq \eps_1$.
Figure \ref{fig:refined-betti} gives a schematic diagram of all these sets.
Using Alexander duality, in order to bound $b(S,\F)$ it suffices to
bound $b(\partial S,\F)$ (see Lemma \ref{lem:refined-betti-duality} below). 
Now, in order to bound $b(S_1 \cup S_2,\F)$ it suffices to bound (using inequality \eqref{eqn:MV:S1cupS2}) $b(S_1,\F), b(S_2,\F)$ as well as $b(S_1\cap S_2,\F)$
(see Lemma \ref{lem:refined-betti-MV} below).

The techniques used for bounding each of the above quantities are distinct.  
We bound $b(S_1,\F)$ by first reducing the problem to bounding 
$b(\partial S_1,\F)$ and bounding $b(\ZZ(P_1^2 -\eps_1,\R\la\eps_1,\eps_2\ra^k),\F)$
using inequality \eqref{eqn:MV:S1cupS2}, 
and then  using Corollary \ref{cor:refined-betti-kouchnirenko} to bound these quantities
(see Proposition \ref{prop:refined-betti-S1}).

In order to bound $b(S_2,\F)$, we observe that a generic linear functional
has no critical points in the relative interior of $S_2$ (see Lemmas \ref{lem:refined-betti-generic1} and \ref{lem:refined-betti-generic2}, and \ref{lem:refined-betti-generic3}). Note that this fact is not
necessarily true for $S_1$.  
This fact allows us to bound $b(S_2,\F)$
by counting the critical points of the functional on its boundary strata using Theorem \ref{thm:SMT}
(see Lemma \ref{lem:refined-betti-smt}).  Finally, we bound the number of such critical points using Proposition \ref{prop:Kouchnirenko} (see Lemma \ref{lem:refined-betti-kouchnirenko}).
The number of such critical points also gives an upper bound on 
$b(S_1 \cap S_2,\F)$ (Corollary \ref{cor:refined-betti-kouchnirenko}).

\subsubsection{Proof of Theorem \ref{thm:refined-betti-main}}
For the rest of this
section we keep the same notation as in Theorem \ref{thm:refined-betti-main}.
 
Let $\R_0= \R\la \eps_0\ra, \R_1=\la \eps_0,\eps_1\ra, \R_2 = \R\la \eps_0,\eps_1,\eps_2\ra$, and let $W_{k+1} \subset \R_0^{k+1}$ denote the real variety 
defined by the polynomials $P_1,P_2$ and $Q_{k+1}= \eps_0\sum_{i=1}^{k+1} X_i^2 -1$, and 
$W_k \subset \R_0^k$ the real variety defined by  the polynomials $P_1,P_2$ and $Q_k= \eps_0\sum_{i=1}^{k} X_i^2 -1$.

It follows from 
\cite[Corollary 9.3.7]{BCR} (Local Conic Structure at infinity of semi-algebraic sets) that,
the intersection, $W_{k+1}^+$ (resp. $W_{k+1}^-$) of $W_{k+1}$ with the closed half-space defined by $X_{k+1} \geq 0$ 
(resp. $X_{k+1} \leq 0$) are each semi-algebraically homeomorphic to $V$.
\begin{eqnarray*}
W_{k+1}^+ \cup W_{k+1}^- &=& W_{k+1}, \\
W_{k+1}^+ \cap W_{k+1}^- &=& W_k.
\end{eqnarray*}

 It now follows from inequality \eqref{eqn:MV:S1+S2} that:
\begin{proposition}
\label{prop:refined-betti-big-sphere}
\[
b(V,\F) \leq \frac{1}{2}\left(b(W_k,\F) + b(W_{k+1},\F)\right).
\]
\end{proposition}

We now bound $b(W_{k+1},\F)$ (the proof for the bound on $b(W_{k},\F)$ is very similar).
Consider the closed and bounded semi-algebraic set, $\tilde{W}_{k+1} \subset \R_2^{k+1}$ defined by
$ {Q}_{k+1}= 0, -\eps_1 \leq P_1 \leq \eps_1, -\eps_2 \leq P_2 \leq \eps_2$.
Let for $i=1,2$,
\begin{eqnarray*}
Z_i^+ &=& \ZZ(\{Q_{k+1},P_i+\eps_i\},\R_i^{k+1}) ,\\
Z_i^- &=& \ZZ(\{Q_{k+1},P_i-\eps_i\},\R_i^{k+1}), \\
S_1^+ &=& \Ext(Z_1^+,\R_2)  \cap \tilde{W}_{k+1}, \\
S_1^- &=& \Ext(Z_1^-,\R_2)  \cap \tilde{W}_{k+1}, \\
S_2^+ &=& Z_2^+ \cap \tilde{W}_{k+1}, \\
S_2^- &=& Z_2^- \cap \tilde{W}_{k+1}.
\end{eqnarray*}

\begin{lemma}
\label{lem:W_1^+}
The semi-algebraic set $\tilde{W}_{k+1}$ is semi-algebraically homotopy equivalent to $\Ext(W_{k+1},\R')$. In particular,
\[
b(\tilde{W}_{k+1},\F) = b(W_{k+1},\F).
\]
\end{lemma}

\begin{proof}
Clearly $\tilde{W}_{k+1}$ is closed and bounded over $\R'$, and $\lim_{\eps_1} \tilde{W}_{k+1}$. The lemma now follows from \cite[Lemma 16.17]{BPRbook2}.
\end{proof}

\begin{lemma}
\label{lem:refined-betti-duality}
%%sb epsilon ----> sigma
\[
b(\tilde{W}_{k+1},\F) \leq \frac{1}{2}\left( b(\bigcup_{\substack{i=1,2 \\\sigma\in\{+,-\}} }S_i^{\sigma},\F) +1\right).
\]
\end{lemma}

\begin{proof}
Let $\tilde{W}_{k+1}'$ be the closure of the semi-algebraic set $\Sphere^{k}(0,\eps_0^{-1/2}) \setminus
\tilde{W}_{k+1}$.
Then,
%%sn epsilon---> sigma
\begin{eqnarray*}
\tilde{W}_{k+1} \cup \tilde{W}_{k+1}'&=& \Sphere^{k}(0,\eps_0^{-1/2}), \\
\tilde{W}_{k+1} \cap \tilde{W}_{k+1}'&=&  \bigcup_{\substack{i=1,2 \\\sigma\in\{+,-\}} }S_i^{\sigma}.
\end{eqnarray*}
We also have that $\tilde{W}_{k+1}'$ is semi-algebraically homotopy equivalent to
$\Sphere^{k}(0,\eps_0^{-1/2}) \setminus
\tilde{W}_{k+1}$, and hence by \cite[page 296]{Spanier} (Alexander duality)
\begin{eqnarray}
\label{eqn:refined-betti-duality1}
b(\tilde{W}_{k+1},\F) &=& b(\tilde{W}_{k+1}',\F) -1 .
\end{eqnarray}
Also, using inequality \eqref{eqn:MV:S1cupS2} we have
%%sb epsilon ---> sigma
\begin{eqnarray}
\label{eqn:refined-betti-duality2}
b(\tilde{W}_{k+1},\F) + b(\tilde{W}_{k+1}',\F) &\leq& b(\bigcup_{\substack{i=1,2 \\\sigma\in\{+,-\}} }S_i^{\sigma},\F) +2.
\end{eqnarray}
The lemma now follows from \eqref{eqn:refined-betti-duality1} and \eqref{eqn:refined-betti-duality2}.
\end{proof}

\begin{lemma}
\label{lem:refined-betti-MV}
%%sb epsilon ---> sigma
\[
b(\bigcup_{\substack{i=1,2\\ \sigma \in \{+,-\}}} S_i^{\sigma},\F) \leq \sum_{\substack{i=1,2 \\ \sigma \in \{+,-\}}}  b(S_i^{\sigma},\F) +
\sum_{\sigma_1,\sigma_2 \in\{ +,-\}} b(S_1^{\sigma_1} \cap S_2^{\sigma_2},\F).
\]
\end{lemma}

\begin{proof}
Apply 
inequality \eqref{eqn:MV:S1cupS2}.
\end{proof}

We need separate arguments to bound 
%%sb \epsilon --->\sigma
$b(\bigcup_{\sigma_1\in\{+,-\}} S_1^{\sigma_1},\F)$ and $b(\bigcup_{\sigma_2\in\{+,-\}} S_2^{\sigma_2},\F)$.

We first bound $b(\bigcup_{\sigma_1\in\{+,-\}} S_1^{\sigma_1},\F)$.
 
\begin{proposition}
\label{prop:refined-betti-S1}
%%sb epsilon ---> sigma
\begin{eqnarray*}
\label{eqn:refined-beti-S1}
b(\bigcup_{\sigma_1 \in \{+,-\}} S_1^{\sigma_1},\F)  
&\leq &  
4\binom{k+1}{2}d_1\left(\frac{2(k-1)}{3}d_2(d_2-1)^{k-2}+ (d_1-1)^{k-1}\right).
\end{eqnarray*}
\end{proposition}

\begin{proof}
Using inequality \eqref{eqn:MV:S1cupS2}  and Corollary \ref{cor:refined-betti-kouchnirenko} we have
%%sb epsilon --> sigma
\begin{eqnarray*}
\label{eqn:refined-betti-S_1}
b(\bigcup_{\sigma_1\in \{+,-\}} S_1^{\sigma_1},\F) &\leq& 
\sum_{\sigma_1 \in \{+,-\}} \left(b(\partial S_1^{\sigma_1},\F) + b(Z_1^{\sigma_1},\F)\right)\\
&=& \sum_{\sigma_1,\sigma_2 \in \{+,-\} }
b(S_2^{\sigma_2} \cap S_1^{\sigma_1},\F)  + \sum_{\sigma_1 \in \{+,-\}} b(Z_1^{\sigma_1},\F)  \\
&\leq &  8\binom{k+1}{3} d_1 d_2 (d_2-1)^{k-2} + 4\binom{k+1}{2}d_1(d_1-1)^{k-1} \\
&=&  
4\binom{k+1}{2}d_1\left(\frac{2(k-1)}{3}d_2(d_2-1)^{k-2} +(d_1-1)^{k-1}\right).
\end{eqnarray*}
\end{proof}

%%sb epsilon ---> sigma
Next we bound $b(\bigcup_{\sigma_2 \in \{+,-\}} S_2^{\sigma_2},\F)$ as follows.

\begin{lemma}
\label{lem:refined-betti-generic1}
There exists a linear functional $F:\R_1^{k+1} \rightarrow \R_1$, such that the set of critical points of $F$ restricted to
$Z_1^{\pm}$ has an empty intersection with $\ZZ(P_2, \R_1^{k+1})$.
\end{lemma}

\begin{proof}
The semi-algebraic subset $T \subset \Gr_{k+1,k}(\R)$ defined by,
\[
T:= \cup_{x \in Z_1^{\pm} \cap \ZZ(P_2, \R''^{k+1})} \{H \in \Gr_{k+1,k}(\R) \mid H \supset T_x Z_1^{\pm} \}
\]
 is of co-dimension at least $1$ in $\Gr_{k+1,k}(\R)$. Thus, the complement of $T$ in $\Gr_{k+1,k}(\R)$ contains an open dense set.
\end{proof}

\begin{lemma}
\label{lem:refined-betti-generic2}
There exists an open dense subset of linear functionals  $F:\R_2^{k+1} \rightarrow \R_2$, such that the set of critical points of $F$ restricted to
$\Ext(Z_1^{\pm},\R_2)$ has an empty intersection with $S_1^{\pm}$.
\end{lemma}

\begin{proof}
Follows from Lemma \ref{lem:refined-betti-generic1} and the fact that 
$S_1^{\pm}$ is bounded over $\R_1$, and 
$\eps_2$  is infinitesimal with respect to $\R_1$.
\end{proof}

\begin{lemma}
\label{lem:refined-betti-generic3}
There exists an open dense subset of linear functionals $F:\R_2^{k+1} \rightarrow \R_2$, such that the critical points of $F$ restricted to
$\Ext(Z_1^{\pm},\R_2)  \cap Z_2^{\pm} $ are non-degenerate.
\end{lemma}

\begin{proof}
The lemma can be deduced as a special case of \cite[Theorem 2]{BGHMS10}.
\end{proof}

\begin{lemma}
\label{lem:refined-betti-smt}
Let $F$ be a linear functional satisfying the hypothesis of Lemmas \ref{lem:refined-betti-generic2} and
\ref{lem:refined-betti-generic3}.
Then, for 
%%sb epsilon --> sigma
$\sigma \in \{+,-\}$, 
$b(S_1^{\sigma},\F)$ is bounded by the number of critical points of $F$ restricted to 
$\Ext(Z_1^{\sigma},\R_2)  \cap Z_2^{\pm}$.
\end{lemma}

\begin{proof}
Follows from Theorem \ref{thm:SMT}.
\end{proof}

\begin{lemma}
\label{lem:refined-betti-kouchnirenko}
Let $F$ be a linear functional satisfying the hypothesis of Lemmas \ref{lem:refined-betti-generic2} and
\ref{lem:refined-betti-generic3} below. 
%%sb epsilon ---> sigma
Then, for $\sigma \in \{+,-\}$, 
$b(S_1^{\sigma},\F)$  the number of critical points of $F$ restricted to 
$\Ext(Z_1^{\sigma},\R_2)  \cap Z_2^{\pm}$ is bounded by 
\[
8\binom{k+1}{3} d_1 d_2 (d_2-1)^{k-2}.
\]
\end{lemma}

\begin{proof}
For $d,k \geq 0$, let $\Delta_{d,k} \subset \Q^k$ denote the simplex defined as the convex hull of
$(d,0,\ldots,0),\ldots, (0,\ldots,0,d), \mathbf{0}$. 
%%sb epsilon --> sigma
For $\sigma_1,\sigma_2 \in \{+,-\}$, the set of critical points of $F$ restricted to  $\Ext(Z_1^{\sigma},\R_2)  \cap Z_2^{\pm}$ satisfies the following system of equations:

\begin{eqnarray}
\label{eqn:refined-betti-mixed}
Q_{k+1} &=& 0, \\ \nonumber
P_1 - \sigma_1 \eps_1 &=& 0, \\ \nonumber
P_2  - \sigma_2 \eps_2 &=& 0, \\ \nonumber
 \frac{\partial Q_{k+1}}{\partial X_1} + \lambda_1 \frac{\partial Q_{k+1}}{\partial X_1} + \lambda_2 \frac{\partial P_1}{\partial X_1} +  \lambda_3 \frac{\partial P_2}{\partial X_1} &=& 0, \\\nonumber
\vdots & \vdots & \vdots \\
 \frac{\partial Q_{k+1}}{\partial X_{k+1}} + \lambda_1 \frac{\partial Q_{k+1}}{\partial X_{k+1}} + \lambda_2 \frac{\partial P_1}{\partial X_{k+1}} +  \lambda_3 \frac{\partial P_2}{\partial X_{k+1}} &=& 0. \nonumber
\end{eqnarray}

Using Proposition \ref{prop:Kouchnirenko}  we obtain that the number of solutions of the system
\eqref{eqn:refined-betti-mixed} is
bounded by
$$\displaylines{
\MV(\Delta_{2,k+1}, \Delta_{d_1,k+1}, \Delta_{d_2,k+1}, \underbrace{\Delta_{d_2-1,k+1}+\Delta_{1,3},\ldots,\Delta_{d_2-1,k+1}+\Delta_{1,3}}_{k+1}) \cr
\leq \binom{k+1}{3} 
\MV(\Delta_{2,k+1}, \Delta_{d_1,k+1}, \Delta_{d_2,k+1}, 
\underbrace{\Delta_{d_2-1,k+1},\ldots,\Delta_{d_2-1,k+1}}_{k-2},
\underbrace{\Delta_{1,3},\ldots,\Delta_{1,3}}_{3}) \cr
=
\binom{k+1}{3} 2 d_1 d_2 (d_2-1)^{k-2}.
}
$$
Hence, 
 the number of critical points of $F$ restricted to 
 %%sb epsilon ---> sigma
$\Ext(Z_1^{\sigma},\R_2)  \cap Z_2^{\pm}$ is bounded by 
\[
4 \binom{k+1}{3} 2 d_1 d_2 (d_2-1)^{k-2} = 
8\binom{k+1}{3} d_1 d_2 (d_2-1)^{k-2}.
\]
\end{proof}
In particular, we obtain as an immediate corollary that
\begin{corollary}
\label{cor:refined-betti-kouchnirenko}
%%sb epsilon --> sigma
\[
\sum_{\sigma_1,\sigma_2 \in\{ +,-\}} b(S_1^{\sigma_1} \cap S_2^{\sigma_2},\F)
\leq 
8\binom{k+1}{3} d_1 d_2 (d_2-1)^{k-2}.
\]
\end{corollary}

\begin{proposition}
\label{prop:refined-betti-S2}
%%sb epsilon ---> sigma
\begin{eqnarray*}
\label{eqn:refined-betti-S2}
b(\bigcup_{\sigma_2 \in \{+,-\}}S_2^{\sigma_2},\F) \leq 8\binom{k+1}{3} d_1 d_2 (d_2-1)^{k-2}.
\end{eqnarray*}
\end{proposition}

\begin{proof}
Follows from Lemmas \ref{lem:refined-betti-generic1}, \ref{lem:refined-betti-generic2}, \ref{lem:refined-betti-generic3}, \ref{lem:refined-betti-smt} and \ref{lem:refined-betti-kouchnirenko}.
\end{proof}

\begin{proposition}
\label{prop:refined-betti-main}
\begin{eqnarray*}
b(W_{k+1},\F)
&\leq &
2\binom{k+1}{2} d_1 
\left( (d_1-1)^{k-1}+ \frac{4(k-1)}{3}d_2(d_2-1)^{k-2} \right) +1, \\
b(W_{k},\F)
&\leq &
2\binom{k}{2} d_1 
\left( (d_1-1)^{k-2}+ \frac{4(k-2)}{3}d_2(d_2-1)^{k-3} \right) +1.
\end{eqnarray*}
\end{proposition}

\begin{proof}
The inequality for $b(W_{k+1},\F)$ follows from Lemmas \ref{lem:refined-betti-duality}, \ref{lem:refined-betti-MV},  Corollary \ref{cor:refined-betti-kouchnirenko}, and Propositions \ref{prop:refined-betti-S1} and 
\ref{prop:refined-betti-S2}.
The proof of the inequality involving $b(W_{k},\F)$ follows the same steps as in the proof
for $b(W_{k+1},\F)$ replacing $k$ by $k-1$.  We omit the steps.
\end{proof}

\begin{proof}[Proof of Theorem \ref{thm:refined-betti-main}]
The theorem follows from Propositions \ref{prop:refined-betti-big-sphere}
and \ref{prop:refined-betti-main}.
\end{proof}

\subsection{The semi-algebraic case: Proof of Theorem \ref{thm:B-B-new}}
\begin{proof}[Proof of Theorem \ref{thm:B-B-new}]
The proof  is by following the proof of Proposition 7.30 in \cite{BPRbook2} using Theorem \ref{thm:refined-betti-main} to bound the Betti numbers of the algebraic sets that arise instead of Theorem \ref{thm:O-P-T-M}.
\end{proof}

\section{Open problems and future directions}
\label{sec:open}
In this section we list some of the open problems and conjectures that could serve 
as future directions for research in this area.
\begin{enumerate}[1.]
\item
It is an interesting open question to generalize the bounds proved in \S \ref{sec:main} to the case
of complex varieties, as well as constructible subsets of $\C^k$,
defined by polynomial systems having similar degree sequences. One could apply the results
in the current paper, treating constructible subsets of $\C^k$, as special cases of semi-algebraic
subsets of $\R^{2k}$, but this would lead to non-optimal upper bounds. On the other hand,
the perturbation techniques, that allowed us to reduce to the case of non-degenerate complete
intersections, used throughout this paper does not apply to the complex case. 
\item
Another open problem is to prove Conjecture \ref{conj:B-B} in its full generality. 
\end{enumerate}

\section*{Acknowledgments}
The authors thank Antonio Lerario and Marie-Fran{\c c}oise Roy, 
for their useful comments on a preliminary version of this paper. We also thank an anonymous referee for a careful reading of the manuscript, and whose 
comments helped us to substantially improve the paper.
\bibliographystyle{abbrv}
\bibliography{master}

\section{Appendix}
\label{sec:Appendix}

\subsection{Betti numbers of generic complex projective and affine complete intersection varieties via Chern class computation}
\label{subsec:Chern}
The Betti numbers of
\emph{projective} varieties which are
non-singular complete intersections can be expressed in terms of the degree sequence defining them
(see for example \cite[\S 5.7.2-3]{Eisenbud-Harris}) using a Chern class computation and Proposition \ref{prop:chi-to-b-projective}.
More precisely,
suppose that
$V = \ZZ(\mathcal{P}, \PP_\C^k)$ is a non-singular complete intersection, where $\mathcal{P} = 
\{P_1,\ldots, P_\ell\}$, $\deg(P_i) = d_i, 1\leq i \leq \ell$.
We denote by $\mathcal{T}_V$ the tangent bundle of $V$, and by $c(\mathcal{T}_V)$ the Chern class of $\mathcal{T}_V$. 
The following formula which is classical (see for example \cite[\S 5.7.2-3]{Eisenbud-Harris}) expresses  $c(\mathcal{T}_V)$  in terms of the 
of the restriction, $\zeta_V \in \HH^2(V,\Z)$,  of the hyperplane class of $\PP_\C^k$ to $V$.
\begin{equation}
\label{eqn:chern}
c(\mathcal{T}_V) = \frac{(1 + \zeta_V)^{k+1}}{\prod_{j=1}^\ell (1 + d_i \zeta_V)},
\end{equation}

The Euler class of  $\mathcal{T}_V$ is the top Chern class $c_{k -\ell}(\mathcal{T}_V)$, and 
the Euler-Poincar\'e characteristic of $V$ equals the Euler number of  $\mathcal{T}_V$ and  is then given by:

\begin{eqnarray*}
\chi(V,\Z) &=& \int_V c_{k-\ell}(\mathcal{T}_V).
\end{eqnarray*}

Moreover, it follows from the fact that the cohomology of $V$ is torsion-free that

\begin{eqnarray}
\label{eqn:chi-int1}
\nonumber
\chi(V,\Z_2) &=& \chi(V,\Z) \\
&=& \int_V c_{k-\ell}(\mathcal{T}_V).
\end{eqnarray}

It follows from Eqn. \eqref{eqn:chern} that the coefficient, $N_{k,\ell,\dd}$,  of  $\zeta_V^{k -\ell}$ is given by:
\begin{eqnarray*}
N_{k,\ell,\dd} &=& \sum_{i=0}^{k-\ell}(-1)^{k-\ell -i}\cdot \binom{k+1}{i} \cdot h_{k-\ell -i}(d_1,\ldots,d_\ell),
\end{eqnarray*}
where
$h_j(d_1,\ldots,d_\ell)$
is the complete homogeneous symmetric polynomial of degree $j$ in $\dd=(d_1,\ldots, d_\ell)$ (cf. Eqn. \eqref{eqn:complete-homogeneous}).

Since,
\[
c_{k -\ell}(\mathcal{T}_V) = N_{k,\ell,\dd} \cdot \zeta_V^{k -\ell},
\]
and 
\[
\int_V \zeta_V^{k-\ell}  =  \deg(\zeta_V^{k-\ell}) = d_1\cdots d_\ell,
\] 
it follows from Eqn. \eqref{eqn:chi-int1} that
 \begin{eqnarray}
 \label{eqn:chi-int2}
 \nonumber
 \chi(V,\Z_2) &=&  N_{k,\ell,\dd} \cdot \int_V \zeta_V^{k-\ell} \\
 &=& N_{k,\ell,\dd} \cdot d_1\cdot d_2 \cdots d_\ell.
 \end{eqnarray}

It now follows from Eqns. \eqref{eqn:chi-int1} and \eqref{eqn:chi-int2}  that,

\begin{eqnarray*}
\chi(V,\Z_2) &=& d_1\cdot d_2\cdots d_\ell\cdot \left(\sum_{i=0}^{k-\ell}(-1)^{k-\ell -i}\cdot\binom{k+1}{i} \cdot h_{k-\ell -i}(d_1,\ldots,d_\ell)  \right)\\
&=& d_1\cdot d_2\cdots d_\ell\cdot \left(\sum_{j=0}^{k-\ell}(-1)^{j}\cdot\binom{k+1}{j+\ell+1} \cdot h_{j}(d_1,\ldots,d_\ell)  \right).
\end{eqnarray*}
Finally using  Eqn. \eqref{eqn:prop:chi-to-b-projective} we have
\begin{eqnarray*}
b(V,\Z_2) &=& (1 + (-1)^{k-\ell+1})\cdot (k-\ell+1) + \\
&& d_1 \cdot d_2\cdots d_\ell\cdot \left(\sum_{i=0}^{k-\ell}(-1)^{i}\cdot\binom{k+1}{i} \cdot h_{k-\ell -i}(d_1,\ldots,d_\ell)  \right) \\
&=& (1 + (-1)^{k-\ell+1})\cdot (k-\ell+1) + \\
&& d_1 \cdot d_2\cdots d_\ell\cdot \left(\sum_{j=0}^{k-\ell}(-1)^{k-\ell -j}\cdot\binom{k+1}{j+\ell+1} \cdot h_{j}(d_1,\ldots,d_\ell)  \right).
\end{eqnarray*}

The above calculations yield via Lefschetz duality formulas for the Euler-Poincar\'e characteristics (and hence the  Betti numbers)  of generic  complex affine
complex intersection varieties as well.
More precisely, if $V = \ZZ(\mathcal{P}, \C^k)$ is a non-singular complete intersection, where $\mathcal{P} = 
\{P_1,\ldots, P_\ell\}$, $\deg(P_i) = d_i, 1 \leq i \leq \ell$. Let $\mathcal{P}^h = \{P_1^h,\ldots,P_\ell^h\}$  denote the homogenization of
$\mathcal{P}$, and $\overline{V} = \ZZ(\mathcal{P}^h,\PP_\C^k)$, and $\overline{W} = \overline{V} \cap H_0$, where $H_0$ is the hyperplane at infinity defined by $X_0=0$.

Then, $V = \overline{V} - \overline{W}$, and using Lefschetz duality (see for example \cite[Page 297, Theorem 19]{Spanier}), we have
\begin{eqnarray*}
\chi(V,\Z_2) &=& \chi(\overline{V},\Z_2) - \chi(\overline{W},\Z_2) \\
&=& d_1\cdot d_2\cdots d_\ell  \left(\sum_{j=0}^{k-\ell}(-1)^{j}\cdot\binom{k+1}{j+\ell+1} h_{j}(d_1,\ldots,d_\ell)  \right) \\
&& - d_1\cdot d_2\cdots d_\ell \left(\sum_{j=0}^{k-\ell-1}(-1)^{j}\cdot\binom{k}{j+\ell+1}  h_{j}(d_1,\ldots,d_\ell)  \right) \\
&=&
d_1\cdot d_2\cdots d_\ell \left(\sum_{j=0}^{k-\ell}(-1)^{j} \left(\binom{k+1}{j+\ell+1} - \binom{k}{j+\ell+1}\right) h_{j}(d_1,\ldots,d_\ell)  \right) \\
&=&
d_1\cdot d_2\cdots d_\ell \left(\sum_{j=0}^{k-\ell-1}(-1)^{j}\cdot\binom{k}{j+\ell} h_{j}(d_1,\ldots,d_\ell)  \right)
\end{eqnarray*}
(cf. Eqn. \eqref{eqn:many-total-diff-degrees-complex-1}).

\subsection{Proofs omitted in the main text}
\label{subsec:proofs}
We include in this section proofs of several theorems and propositions which were omitted in the main text.

\begin{proof}[Proof of Theorem \ref{thm:Smith}]
The only difficulty in applying inequality \eqref{eqn:Smith} is that in general an affine sub-variety
of $\C^k$ will not bounded.
In order to apply inequality \eqref{eqn:Smith} we need to reduce to the closed and bounded case which
we do as follows. 

Let for $r>0$, $B_\C(r) \subset \C^k$ be defined by 
\[
B_\C(r) = \{ (x_1+i y_1,\ldots,x_k+i y_k) \in \C^k \mid |x_i|, |y_i| \leq r, 1 \leq i \leq k\},
\]
and denote by $B_\R(r) = B_\C(r) \cap \R^k$.

Then, $\ZZ(\mathcal{Q},\C^k) \cap B_\C(r)$ is closed and bounded, and using 
\cite[Corollary 9.3.7]{BCR} (Local Conic Structure at infinity of semi-algebraic sets), we have that for all $r>0$ and  large enough,
$\ZZ(\mathcal{Q},\C^k) \cap B_\C(r)$ is semi-algebraically homeomorphic to $\ZZ(\mathcal{Q},\C^k)$,
and $\ZZ(\mathcal{Q},\R^k) \cap B_\R(r)$ is semi-algebraically homeomorphic to $\ZZ(\mathcal{Q},\R^k)$.

The complex conjugation restricts to an involution of $B_\C(r)$ with fixed points $B_\R(r)$.
Now apply inequality \eqref{eqn:Smith}.
\end{proof}
\begin{proof}[Proof of Proposition \ref{prop:one-total}]
First observe that $\ZZ(P,\C^k)$ is either $0$-dimensional, or is smooth and connected in case $k >1$ (since $\ZZ(P,\C^k)$ is a non-singular
projective hypersurface of dimension $k-1$ minus a sub-variety of strictly smaller dimension).
Using Theorem  \ref{thm:Khovansky} we obtain
\begin{eqnarray*}
\chi(\ZZ(P,\C^k),\Z_2) &=& \sum_{j=1}^{k} \binom{k}{j} (-1)^{j+1} j! \frac{d^j}{j!} \\
&=& 1+\sum_{j=0}^{k} \binom{k}{j} (-1)^{j+1} j! \frac{d^j}{j!} \\
            &=& 1-  (1 - d)^k \\
             &=& 1+ (-1)^{k-1} (d-1)^k.
\end{eqnarray*}
This implies using Proposition \ref{prop:chi-to-b}  
that 

\begin{eqnarray*}
b(\ZZ(P,\C^k),\Z_2) 
&=& 1 + (-1)^{k-1}(\chi(\ZZ(P,\C^k),\Z_2) -1) \\
&=& 1 + (-1)^k +  (-1)^{k-1} \chi(\ZZ(P,\C^k),\Z_2) \\
&=& 1 + (d-1)^k.
\end{eqnarray*} 

This proves Eqn.\eqref{eqn:one-total-complex}.

Finally, inequality \eqref{eqn:one-total-real} follows from Eqn. \eqref{eqn:one-total-complex} and
Theorem \ref{thm:Smith} (Smith inequality).
\end{proof}

\begin{proof}[Proof of  Proposition \ref{prop:many-total-diff-degrees}]
Let $V = \ZZ(\mathcal{P},\C^k)$.
First observe that  either $V$ is $0$-dimensional (in case $k=\ell$) or is non-singular and connected (if $k > \ell$)
since in the latter case  $V$ is equal to 
a non-singular complete intersection variety in a product of projective varieties minus a sub-variety
of strictly smaller dimension. 

Using Theorem  \ref{thm:Khovansky} we have
\begin{eqnarray}
\label{eqn:many-total-diff-degrees}
\nonumber
\chi(V,\Z_2) &=& \sum_{j=\ell}^{k} (-1)^{j+\ell}\binom{k}{j} d_1\cdots d_\ell \sum_{\substack{j_1,\ldots,j_\ell \geq 0\\ j_1+\cdots + j_\ell = j -\ell}} d_1^{j_1}\cdots d_\ell^{j_\ell}  \\
\nonumber            &=&  
            d_1\cdots d_\ell \cdot\left(\sum_{j=\ell}^{k} (-1)^{j+\ell}\binom{k}{j} h_{j - \ell}(d_1,\ldots,d_\ell)\right)  \\
             &=&  d_1\cdots d_\ell \cdot \left( \sum_{j=0}^{k-\ell } (-1)^{j}\binom{k}{j+\ell} h_{j}(d_1,\ldots,d_\ell) \right).
\end{eqnarray}

Eqns. \eqref{eqn:many-total-diff-degrees} and Proposition \ref{prop:chi-to-b} imply that 
\begin{eqnarray*}
b(V,\Z_2)
&=& 1 + (-1)^{k-\ell}(\chi(\ZZ(\mathcal{P},\C^k),\Z_2) -1) \\
&=& 1 + (-1)^{k-\ell+1} + (-1)^{k-\ell} \chi(\ZZ(\mathcal{P},\C^k),\Z_2) \\
&=& 1 + (-1)^{k- \ell+1} +  d_1\cdots d_\ell \left( \sum_{j=0}^{k-\ell } (-1)^{k-\ell+j}\binom{k}{j+\ell} h_{j}(d_1,\ldots,d_\ell) \right)
\end{eqnarray*}
(recovering the same result proved in \S \ref{subsec:Chern} using Chern class computations).

Now assume that $d_1 = \cdots = d_\ell = d$. It follows from \eqref{eqn:many-total-diff-degrees} that
\begin{eqnarray*}
\chi(V,\Z_2) &=& \sum_{j=\ell}^{k} (-1)^{j+\ell}\binom{k}{j}   \binom{j-1}{\ell-1}  d^j \\
            &=&   \sum_{j=\ell}^{k} (-1)^{j+\ell}\frac{k!}{j!(k-j)!} \frac{(j-1)!}{(j-\ell)!(\ell-1)!}  d^j  \\
             &=& \ell\binom{k}{\ell}\sum_{j=\ell}^{k} (-1)^{j+\ell} \binom{k-\ell}{j- \ell} \frac{d^j}{j}  \\
              &=& \ell\binom{k}{\ell}\sum_{j=0}^{k-\ell} (-1)^{j} \binom{k-\ell}{j} \frac{d^{j+\ell}}{j+\ell}  \\
              &=&\ell\binom{k}{\ell}\sum_{j=0}^{k-\ell} (-1)^{j} \binom{k-\ell}{j} \frac{d^{j+\ell}}{j+\ell}  \\
               &=& \ell\binom{k}{\ell}\int_{0}^{d} x^{\ell-1}(1-x)^{k-\ell} \mathrm{d}x.           
\end{eqnarray*}
This implies
\begin{eqnarray*}
|\chi(V,\Z_2)| &\leq& \ell \binom{k}{\ell}\left(\int_{0}^{1} x^{\ell-1}(1-x)^{k-\ell} \mathrm{d}x + \int_{1}^{d} x^{\ell-1}(x-1)^{k-\ell} \mathrm{d}x \right) \\
&\leq& \ell \binom{k}{\ell}\left(1 + \int_{1}^{d} x^{k-1}\mathrm{d}x \right) \\
&=& \ell \binom{k}{\ell}\left(1 + \frac{d^k}{k}  - \frac{1}{k}\right) \\
&=&  \binom{k-1}{\ell-1} (d^k  + k - 1),
\end{eqnarray*}
whence using Proposition  \ref{prop:chi-to-b} 
\begin{eqnarray*}
b(V,\Z_2) 
&=& 1 + (-1)^{k-\ell}(\chi(V,\Z_2) -1) \\
&=& 1 + (-1)^{k-\ell+1} + (-1)^{k-\ell} \chi(V,\Z_2) \\
&\leq& 1 + (-1)^{k- \ell+1} +  \binom{k-1}{\ell-1} (d^k  + k - 1).
\end{eqnarray*}

This proves Eqn. \eqref{eqn:many-total-diff-degrees-complex-1} 
and inequality
\eqref{eqn:many-total-diff-degrees-complex-2}.
Inequalities \eqref{eqn:many-total-diff-degrees-real-1}
and
\eqref{eqn:many-total-diff-degrees-real-2},
follow from 
Eqn. 
\eqref{eqn:many-total-diff-degrees-complex-1}, 
inequality 
\eqref{eqn:many-total-diff-degrees-complex-2}, 
and 
Theorem \ref{thm:Smith} (Smith inequality).
\end{proof}

\begin{proof}[Proof of  Proposition  \ref{prop:manyblocks-total}]
It follows from Theorem \ref{thm:Khovansky} that
\begin{eqnarray*}
\chi(V,\Z_2) &=& 
\sum_{j=\ell}^{k} (-1)^{j+\ell} \sum_{\substack{j_1,\ldots, j_p \\ 0 \leq j_i \leq k_i, 1\leq i \leq p\\ j_1+\cdots+j_p=j}}  \binom{j-1}{\ell-1}\binom{j}{j_1,\ldots,j_p} \left(\prod_{i=1}^{p} \binom{k_i}{j_i} d_i^{j_i}\right) \\
&=& 
\frac{k_1!\cdots k_p!}{(k-\ell)!(\ell-1)!}  F(\kk,\jj,\dd,\ell),
\end{eqnarray*}
where $\kk = (k_1,\ldots,k_p), \; \jj=(j_1,\ldots,j_p),  \; \dd=(d_1,\ldots,d_p)$, and 
$F(\kk,\jj,\dd,\ell)$ is defined as
\[
\sum_{j=\ell}^{k} \frac{(-1)^{j+\ell}}{j} 
\sum_{\substack{j_1,\ldots, j_p \\ 0 \leq j_i \leq k_i, 1\leq i \leq p\\ j_1+\cdots+j_p=j}}  \binom{j}{j_1,\ldots,j_p}^2 \binom{k-\ell}{k_1-j_1,\ldots,k_p-j_p,j-\ell} d_1^{j_1}\cdots d_p^{j_p}.
\]
We now bound $|F(\kk,\jj,\dd,\ell)|$ as follows.
$|F(\kk,\jj,\dd,\ell)|$ is bounded by
\begin{eqnarray*}
&& 
\left(\prod_{i=1}^p d_i^{k_i}\right) \sum_{j=\ell}^{k} 
(1+p)^{2j} 
\sum_{\substack{j_1,\ldots, j_p \\ 0 \leq j_i \leq k_i \\ 1\leq i \leq p}} \binom{k-\ell}{k_1-j_1,\ldots,k_p-j_p,j-\ell} \prod_{i=1}^{p}d_i^{-(k_i-j_i)} \\
&=& 
\left(\prod_{i=1}^p d_i^{k_i}\right)\sum_{j=\ell}^{k} 
(1+p)^{2j}  \left(1 + \frac{1}{d_1} + \cdots + \frac{1}{d_p}\right)^{k -\ell}   \\
&\leq&\left(\prod_{i=1}^p d_i^{k_i}\right)(1+p)^{k+\ell} \sum_{j=0}^{k-\ell} (1+p)^{2j}  \\
&=&  \left(\prod_{i=1}^p d_i^{k_i}\right) (1+p)^{k+\ell} \frac{(1+p)^{2(k-\ell+1)}-1}{(1+p)^2-1} \\
&\leq&  \frac{(1+p)^{3k-\ell+1}}{p(p+2)} d_1^{k_1}\cdots d_p^{k_p}.
\end{eqnarray*}
This implies that
\begin{eqnarray*}
|\chi(V,\Z_2)| &\leq&
\ell(k -\ell+1)\binom{k}{\ell}\binom{k}{\kk}^{-1} \frac{(1+p)^{3k-\ell+1}}{p(p+2)} d_1^{k_1}\cdots d_p^{k_p} \\
&\leq&
(k-\ell+2)^2\binom{k}{\ell-1}\binom{k}{\kk}^{-1} \frac{(1+p)^{3k-\ell+1}}{p(p+2)} d_1^{k_1}\cdots d_p^{k_p},
\end{eqnarray*}
and also (using Eqn. \eqref{eqn:chi-to-b})
\begin{eqnarray*}
b(V,\Z_2) 
&\leq& 1 + (-1)^{k- \ell+1} +  (k-\ell+2)^2\binom{k}{\ell-1}\binom{k}{\kk}^{-1} \frac{(1+p)^{3k-\ell+1}}{p(p+2)} d_1^{k_1}\cdots d_p^{k_p}.
\end{eqnarray*}
This proves inequality \eqref{eqn:manyblocks-total-complex}. Inequality \eqref{eqn:manyblocks-total-real} follows from inequality \eqref{eqn:manyblocks-total-complex} and Theorem \ref{thm:Smith} (Smith inequality).
\end{proof}

\begin{proof}[Proof of Proposition \ref{prop:one-multi}]
First observe that $\ZZ(P,\C^k)$ is non-singular and connected for $k>1$, and is $0$-dimensional if $k=1$.
Using Theorem \ref{thm:Khovansky}
\begin{eqnarray*}
\chi(\ZZ(P,\C^k),\Z_2) &=& \sum_{j=1}^{k} (-1)^{j+1} \sum_{\substack{J  \subset [1,k] \\ \card(J) = j \leq k }} j! \bar{d}^{J}.
\end{eqnarray*}

Now using Eqn. \eqref{eqn:chi-to-b} we get 
\begin{eqnarray*}
b(\ZZ(P,\C^k) ,\Z_2) 
&=& 1 + (-1)^{k-1}(\chi(V_k,\Z_2) -1) \\
&=& 1 + (-1)^k + (-1)^{k-1} \chi(V_k,\Z_2) \\
&=& 1 + (-1)^k + (-1)^{k-1} \left(  \sum_{j=1}^{k} (-1)^{j+1} \sum_{\substack{J  \subset [1,k] \\ \card(J) = j \leq k }} j! \bar{d}^{J}  \right) \\
&=& 1 + (-1)^k + \left(  \sum_{j=1}^{k} (-1)^{k-j} \sum_{\substack{J  \subset [1,k] \\ \card(J) = j \leq k }} j! \bar{d}^{J}  \right).
\end{eqnarray*}
This proves inequality \eqref{eqn:one-multi-complex}. Inequality \eqref{eqn:one-multi-real} now follows
from Eqn. \eqref{eqn:one-multi-complex} and Theorem \ref{thm:Smith} (Smith inequality).
\end{proof}

\begin{proof}[Proof of Proposition \ref{prop:different-boxes}]
It follows from Theorem \ref{thm:Khovansky} that $|\chi(\ZZ(\mathcal{P},\C^k),\Z_2)|$ is 
\begin{eqnarray*}
&\leq& \sum_{j=\ell}^{k} \sum_{J\in \binom{[1,k]}{j}} \binom{j+\ell-1}{\ell-1}
\max_{\substack{\boldsymbol{\alpha} = (\alpha_1,\ldots,\alpha_\ell) \in \Z_{>0}^\ell \\ \alpha_1+\cdots+\alpha_\ell = k}} \max_{(J_1,\ldots,J_\ell) \in \binom{[1,k]}{\boldsymbol{\alpha}}} \left(
\prod_{\substack{1\leq i \leq \ell \\ j \in J_i}} d_{i,j}
\right),
\end{eqnarray*}
So we have that
\begin{equation}
\label{eqn:different-boxes-C}
|\chi(\ZZ(\mathcal{P},\C^k),\Z_2)| \leq O(\ell)^k \max_{\substack{\boldsymbol{\alpha} = (\alpha_1,\ldots,\alpha_\ell) \in \Z_{>0}^\ell \\ \alpha_1+\cdots+\alpha_\ell = k}} \max_{(J_1,\ldots,J_\ell) \in \binom{[1,k]}{\boldsymbol{\alpha}}} \left(
\prod_{\substack{1\leq i \leq \ell \\ j \in J_i}} d_{i,j}
\right).
\end{equation}

Therefore, combining inequality \eqref{eqn:different-boxes-C} and Proposition \ref{prop:chi-to-b}, we have that
\begin{eqnarray}
\nonumber b(\ZZ(\mathcal{P},\C^k),\Z_2) &=& 1+(-1)^{k-\ell+1}+(-1)^{k-\ell}\chi(\ZZ(\mathcal{P},\C^k),\Z_2)\\
&\leq& O(\ell)^k \max_{\substack{\boldsymbol{\alpha} = (\alpha_1,\ldots,\alpha_\ell) \in \Z_{>0}^\ell \\ \alpha_1+\cdots+\alpha_\ell = k}} \max_{(J_1,\ldots,J_\ell) \in \binom{[1,k]}{\boldsymbol{\alpha}}} \left(
\prod_{\substack{1\leq i \leq \ell \\ j \in J_i}} d_{i,j}
\right).
\end{eqnarray}

This proves inequality \eqref{eqn:different-boxes-complex}. 
The inequality in the case $\mathcal{P}$ has coefficients in $\R$ 
follows from inequality \eqref{eqn:different-boxes-complex} and Theorem \ref{thm:Smith}
(Smith inequality).
\end{proof}

\begin{proof}[Proof of  Proposition \ref{prop:many-total-mixed}]
Using Theorem \ref{thm:Khovansky} we obtain that $\chi(\ZZ(\mathcal{P},\C^k),\Z_2) $ equals
\begin{eqnarray*}
&& \sum_{j=0}^{k} (-1)^{j+\ell} \sum_{\substack{0 \leq j_1 \leq k_1, 0\leq j_2 \leq k_2 \\ j_1+j_2 = j}}  \binom{k_1}{j_1}\binom{k_2}{j_2} \binom{j-1}{\ell-1}\binom{j}{j_1}d^{j_1} 2^{j_2}+1\\
&=&  \sum_{j_1=0}^{k_1} (-1)^{j_1+\ell}
\binom{k_1}{j_1}d^{j_1} 
\left(\sum_{j_2=0}^{k_2} \binom{j_1+j_2}{j_2}  \binom{j_1+j_2-1}{\ell-1}\binom{k_2}{j_2} (-2) ^{j_2}
\right) + 1.
\end{eqnarray*}

We now bound from above the quantity  $|F(j_1,k_2)|$, where $F(j_1,k_2)$ is defined by
\[
F(j_1,k_2)  := \sum_{j_2=0}^{k_2} \binom{j_1+j_2}{j_2} \binom{j_1+j_2-1}{\ell-1} \binom{k_2}{j_2} (-2)^{j_2}.
\]
First notice that $F(j_1,k_2)$ equals, setting $\beta = \beta(j_1,\ell)=  \frac{1}{j_1!(\ell-1)!}$,
\begin{eqnarray*}
 && 
  \beta
\sum_{j_2=0}^{k_2} (j_1+j_2)^{\underline{j_1}} (j_1+j_2-1)^{\underline{\ell-1}}
\binom{k_2}{j_2} (-2)^{j_2} \\
&=& 
\beta
\left[
\frac{\mathrm{d} ^{j_1}}{\mathrm{d}{x^{j_1}}}
x^\ell \left( \frac{\mathrm{d} ^{\ell-1}}{\mathrm{d}{x^{\ell-1}}} (x^{j_1-1}(1+x)^{k_2})\right)\right]_{x=-2} \\
&=& 
\beta
\left[ 
\sum_{i=0}^{j_1} \binom{j_1}{i} \ell^{\underline{i}} x^{\ell-i} 
\left( \frac{\mathrm{d} ^{\ell-1+j_1-i}}{\mathrm{d}{x^{\ell-1+j_1-i}}} (x^{j_1-1}(1+x)^{k_2})
\right)
\right]_{x=-2} 
\\
&=& 
\beta
\left[ 
\sum_{i=0}^{j_1} \binom{j_1}{i} \ell^{\underline{i}} x^{\ell-i} 
\left( 
\sum_{h=0}^{\alpha}
\binom{\alpha}{h}
(j_1-1)^{\underline{h}} k_2^{\underline{\alpha-h}}
x^{j_1-1-h}(1+x)^{k_2- (\alpha-h)}
\right)
\right]_{x=-2}\\
\end{eqnarray*}
where 
$\alpha = \alpha(\ell,j_1,i) =  \ell-1+j_1 -i $,  and we have used the ``falling factorial'' notation
\[
t^{\underline{n}} :=t(t-1)\cdots(t-n+1),
\]
for all real $t$ and integer $n$.

Continuing, we have $F(j_1,k_2)$ equals
\begin{eqnarray*}
 &&
\beta
\left[ 
\sum_{i=0}^{j_1} 
\binom{j_1}{i}
\binom{\ell}{i}
\alpha! i!
x^{\alpha}(1+x)^{k_2-\alpha}
\left( 
\sum_{h=0}^{\alpha}
\binom{k_2}{\alpha-h}
\binom{j_1-1}{h} (\omega(x))^h
\right)
\right]_{x=-2},
 \end{eqnarray*}
 where $\omega(x) = 1+ \frac{1}{x}$.

This implies that
\begin{eqnarray*}
|F(j_1,k_2)| &\leq&  
\beta
\sum_{i=0}^{j_1} 
\binom{j_1}{i}
\binom{\ell}{i}
\alpha! i!
2^{\alpha}
\left( 
\sum_{h=0}^{\alpha}
\binom{k_2}{\alpha-h}
\binom{j_1-1}{h}
\right)
\\
&\leq&
\beta
\sum_{i=0}^{j_1} 
\binom{j_1}{i}
\binom{\ell}{i}
\alpha! i!
2^{\alpha}
\binom{k_2 + j_1 -1}{\alpha}.
\end{eqnarray*}

We obtain 
\begin{eqnarray*}
|\chi(V\ZZ(\mathcal{P},\C^k),\Z_2)| &\leq & 
1+
\sum_{j_1=0}^{k_1}  \binom{k_1}{j_1}d^{j_1} |F(j_1,k_2)| \\
&\leq&
1+\sum_{j_1=0}^{k_1} \sum_{i=0}^{j_1}
 \binom{k_1}{j_1}d^{j_1}
\beta(j_1,\ell)
 \binom{j_1}{i}
\binom{\ell}{i}
\alpha! i!
2^{\alpha}
\binom{k_2 + j_1 -1}{\alpha} \\
&\leq&
1+\ell \sum_{j_1=0}^{k_1} \sum_{i=0}^{j_1}
 \binom{k_1}{j_1} d^{j_1} \alpha!2^{\alpha} \binom{k_2 + j_1 -1}{\alpha} \\
&\leq & 1+ \ell \sum_{j_1=0}^{k_1} \sum_{i=0}^{j_1}
 \binom{k_1}{j_1} d^{j_1} 2^{\alpha} (k_2+j_1)^{\alpha} \\
 &\leq&  1+ 2\ell \sum_{j_1=0}^{k_1} \binom{k_1}{j_1}d^{j_1} \left( 2(k_2+k_1) \right)^{\ell-1+j_1} \\
  &=& 1+ \ell 2^\ell (k_1+k_2)^{\ell-1} \sum_{j_1=0}^{k_1} \binom{k_1}{j_1} d^{j_1} \left( 2(k_1+k_2) \right)^{j_1} \\
 &=& 1+ \ell 2^\ell (k_1+k_2)^{\ell-1} \left( 2d(k_1+k_2)+1 \right)^{k_1}.
 \end{eqnarray*}

Now using Eqn. \eqref{eqn:chi-to-b} we get
\begin{eqnarray*}
b(\ZZ(\mathcal{P},\C^k),\Z_2) 
&=& 2 + (-1)^{k-\ell}(\chi(V_k,\Z_2) -1) \\
&\leq & 2 + (-1)^{k-\ell +1} + \ell 2^\ell (k_1+k_2)^{\ell-1} \left( 2d(k_1+k_2)+1 \right)^{k_1}.
\end{eqnarray*}
which proves inequality \eqref{eqn:many-total-mixed-complex}. 
The inequality in the real case 
follows from 
inequality \eqref{eqn:many-total-mixed-complex} and Theorem \ref{thm:Smith} (Smith inequality).
\end{proof}

\begin{proof}[Proof of Proposition \ref{prop:several-blocks-mixed}]
First observe that either $\ZZ(\mathcal{P},\C^k)$ is $0$ dimensional (in case $k = \ell$), or 
$\ZZ(\mathcal{P},\C^k)$ is non-singular and connected (in case $k > \ell$).
We denote by $\overline{d}=\prod_{j=1}^{k_1}d_j$ and for a subset $J_1\subset [1,k_1], \; \overline{d}^{J_1}=\prod_{j\in J_1} d_j$. Then, using Theorem \ref{thm:Khovansky}, $\chi(\ZZ(\mathcal{P},\C^k),\Z_2)$ equals (setting $\epsilon(j,\ell) = (-1)^{j+\ell}$)
\begin{eqnarray*}
 && \sum_{j=\ell}^{k} 
 \epsilon(j,\ell)
  \sum_{\substack{J=J_1 \sqcup J_2 \\ \card(J_1)=j_1\leq k_1, \card(J_2)=j_2 \leq k_2 \\ j=j_1+j_2}} \binom{j-1}{\ell-1} \frac{j!}{j_2!} \overline{d}^{J_1}2^{j_2} \\
&=& \sum_{j_1=0}^{k_1} 
\epsilon(j_1,\ell)
\sum_{\substack{J_1\subset [1,k_1] \\ \card(J_1)=j_1}} \overline{d}^{J_1} \left( \sum_{\substack{J_2\subset [1,k_2] \\ \card(J_2)=j_2}} (-1)^{j_2} \binom{j_1+j_2-1}{\ell-1} \frac{(j_1+j_2)!}{j_2!} 2^{j_2} \right)+1 \\
&=& \sum_{j_1=0}^{k_1} 
\epsilon(j_1,\ell)
 \sum_{\substack{J_1 \subset [1,k_1] \\ \card(J_1)=j_1}} j_1! \overline{d}^{J_1} \left( \sum_{j_2=0}^{k_2} (-1)^{j_2} \binom{j_1+j_2-1}{\ell-1} \binom{j_1+j_2}{j_2} \binom{k_2}{j_2} 2^{j_2} \right)+1.
\end{eqnarray*}
Note that the last sum is the same function $F(j_1,k_2)$ as in Proposition \ref{prop:many-total-mixed}. Applying the same bound, we have $| \chi(\ZZ(\mathcal{P},\C^k),\Z_2) |$ is bounded by 
\begin{eqnarray*}
 && 1+ \sum_{j_1=0}^{k_1} \sum_{\substack{J_1\subset [1,k_1] \\ \card(J_1)=j_1}} j_1!\overline{d}^{J_1} |F(j_1,k_2)| \\
&\leq& 1+ \ell 2^\ell (k_1+k_2)^{\ell-1} \sum_{j_1=0}^{k_1} j_1! (2(k_1+k_2))^{j_1} \sum_{\substack{J_1\subset [1,k_1] \\ \card(J_1)=j_1}} \overline{d}^{J_1} \\
&\leq& 1+ \ell 2^\ell (k_1+k_2)^{\ell-1} \sum_{j_1=0}^{k_1} j_1! (2(k_1+k_2))^{j_1}\binom{k_1}{j_1}d_1\cdots d_{k_1} \\
&\leq& 1+ \ell 2^\ell k_1! (k_1+k_2)^{\ell-1} d_1\cdots d_{k_1} \sum_{j_1=0}^{k_1} \binom{k_1}{j_1} (2(k_1+k_2))^{j_1} \\
&=& 1+ \ell 2^\ell k_1! (k_1+k_2)^{\ell-1}((2(k_1+k_2)+1)^{k_1}d_1\cdots d_{k_1}.
\end{eqnarray*}
Therefore, using Proposition \ref{prop:chi-to-b}, we have
$b(\ZZ(\mathcal{P},\C^k),\Z_2)$ equals
\begin{eqnarray*}
&& 1+(-1)^{k-\ell}(\chi(\ZZ(\mathcal{P},\C^k),\Z_2)-1) \\
&\leq& 2+(-1)^{k-\ell+1}+\ell 2^\ell k_1! (k_1+k_2)^{\ell-1}((2(k_1+k_2)+1)^{k_1}d_1\cdots d_{k_1}.
\end{eqnarray*}
This proves inequality \eqref{eqn:several-blocks-mixed-complex}. 
The inequality in the real case 
follows from inequality \eqref{eqn:several-blocks-mixed-complex}
and Theorem \ref{thm:Smith} (Smith inequality).
\end{proof}

\end{document}